\documentclass{amsart}

\usepackage{amssymb, amsbsy, amsthm, amsmath, amstext, amsopn, verbatim}
\usepackage[all]{xy}
\usepackage{amsfonts}
\usepackage{amscd}
\newtheorem{thm}{Theorem} [section]
\newtheorem{lemma}[thm]{Lemma}

\newtheorem{corollary}[thm]{Corollary}
\newtheorem{prop}[thm]{Proposition}
\newtheorem{notation}[thm]{Notation}

\theoremstyle{definition}

\newtheorem{defn}[thm]{Definition}

\theoremstyle{remark}

\newtheorem{remark}[thm]{Remark}
\newtheorem*{note}{Note}
\newtheorem{claim}[thm]{Claim}

\begin{document}

\numberwithin{equation}{section}

\newcommand{\hs}{\mbox{\hspace{.4em}}}
\newcommand{\ds}{\displaystyle}
\newcommand{\bd}{\begin{displaymath}}
\newcommand{\ed}{\end{displaymath}}
\newcommand{\bcd}{\begin{CD}}
\newcommand{\ecd}{\end{CD}}

\newcommand{\on}{\operatorname}
\newcommand{\proj}{\operatorname{Proj}}
\newcommand{\Proj}{\operatorname{Proj}}
\newcommand{\bproj}{\underline{\operatorname{Proj}}}
\newcommand{\spec}{\operatorname{Spec}}
\newcommand{\bspec}{\underline{\operatorname{Spec}}}
\newcommand{\pline}{{\mathbf P} ^1}
\newcommand{\aline}{{\mathbf A} ^1}
\newcommand{\pplane}{{\mathbf P}^2}
\newcommand{\aplane}{{\mathbf A}^2}
\newcommand{\coker}{{\operatorname{coker}}}
\newcommand{\ldb}{[[}
\newcommand{\rdb}{]]}

\newcommand{\Sym}{\operatorname{Sym}^{\bullet}}
\newcommand{\Symp}{\operatorname{Sym}}
\newcommand{\Pic}{\operatorname{Pic}}
\newcommand{\AAut}{\operatorname{Aut}}
\newcommand{\PAut}{\operatorname{PAut}}

\newcommand{\too}{\twoheadrightarrow}
\newcommand{\Z}{{\mathbf Z}}
\newcommand{\C}{{\mathbf C}}
\newcommand{\Cx}{{\mathbf C}^{\times}}
\newcommand{\cA}{{\mathcal A}}
\newcommand{\cS}{{\mathcal S}}
\newcommand{\cV}{{\mathcal V}}
\newcommand{\cM}{{\mathcal M}}
\newcommand{\bA}{{\mathbf A}}
\newcommand{\cB}{{\mathcal B}}
\newcommand{\cC}{{\mathcal C}}
\newcommand{\cD}{{\mathcal D}}
\newcommand{\D}{{\mathcal D}}
\newcommand{\cs}{{\mathbf C} ^*}
\newcommand{\boldc}{{\mathbf C}}
\newcommand{\cE}{{\mathcal E}}
\newcommand{\cF}{{\mathcal F}}
\newcommand{\bF}{{\mathbf F}}
\newcommand{\cG}{{\mathcal G}}
\newcommand{\G}{{\mathbf G}}
\newcommand{\cH}{{\mathcal H}}
\newcommand{\cJ}{{\mathcal J}}
\newcommand{\cK}{{\mathcal K}}
\newcommand{\cL}{{\mathcal L}}
\newcommand{\baL}{{\overline{\mathcal L}}}
\newcommand{\M}{{\mathcal M}}
\newcommand{\bM}{{\mathbf M}}
\newcommand{\bm}{{\mathbf m}}
\newcommand{\cN}{{\mathcal N}}
\newcommand{\theo}{\mathcal{O}}
\newcommand{\cP}{{\mathcal P}}
\newcommand{\cR}{{\mathcal R}}
\newcommand{\boldp}{{\mathbf P}}
\newcommand{\boldq}{{\mathbf Q}}
\newcommand{\bbL}{{\mathbf L}}
\newcommand{\cQ}{{\mathcal Q}}
\newcommand{\cO}{{\mathcal O}}
\newcommand{\Oo}{{\mathcal O}}
\newcommand{\OX}{{\Oo_X}}
\newcommand{\OY}{{\Oo_Y}}
\newcommand{\otY}{{\underset{\OY}{\ot}}}
\newcommand{\otX}{{\underset{\OX}{\ot}}}
\newcommand{\cT}{{\mathcal T}}
\newcommand{\cU}{{\mathcal U}}\newcommand{\cX}{{\mathcal X}}
\newcommand{\cW}{{\mathcal W}}
\newcommand{\boldz}{{\mathbf Z}}
\newcommand{\qgr}{\operatorname{q-gr}}
\newcommand{\gr}{\operatorname{gr}}
\newcommand{\rk}{\operatorname{rk}}
\newcommand{\coh}{\operatorname{coh}}
\newcommand{\End}{\operatorname{End}}
\newcommand{\Hom}{\operatorname{Hom}}
\newcommand{\uHom}{\underline{\operatorname{Hom}}}
\newcommand{\uHomY}{\uHom_{\OY}}
\newcommand{\uHomX}{\uHom_{\OX}}
\newcommand{\Ext}{\operatorname{Ext}}
\newcommand{\bExt}{\operatorname{\bf{Ext}}}
\newcommand{\Tor}{\operatorname{Tor}}

\newcommand{\inv}{^{-1}}
\newcommand{\airtilde}{\widetilde{\hspace{.5em}}}
\newcommand{\airhat}{\widehat{\hspace{.5em}}}
\newcommand{\nt}{^{\circ}}
\newcommand{\del}{\partial}

\newcommand{\supp}{\operatorname{supp}}
\newcommand{\GK}{\operatorname{GK-dim}}
\newcommand{\hd}{\operatorname{hd}}
\newcommand{\id}{\operatorname{id}}
\newcommand{\res}{\operatorname{res}}
\newcommand{\lrar}{\leadsto}
\newcommand{\im}{\operatorname{Im}}
\newcommand{\HH}{\operatorname{H}}
\newcommand{\TF}{\operatorname{TF}}
\newcommand{\Bun}{\operatorname{Bun}}
\newcommand{\BunD}{\operatorname{Bun}_{\D}}
\newcommand{\PicD}{\operatorname{Pic}_{\D}}
\newcommand{\Hilb}{\operatorname{Hilb}}
\newcommand{\Fact}{\operatorname{Fact}}
\newcommand{\CM}{\mathfrak{CM}}
\newcommand{\ECM}{\mathfrak{ECM}}
\newcommand{\MD}{\mathfrak{M}^{\D}}
\newcommand{\F}{\mathcal{F}}
\newcommand{\Ff}{\mathbb{F}}
\newcommand{\nthord}{^{(n)}}
\newcommand{\Aut}{\underline{\operatorname{Aut}}}
\newcommand{\Gr}{{\on{Gr}}}
\newcommand{\GR}{\operatorname{GR}}
\newcommand{\GRo}{{\operatorname{GR}^{\circ}}}
\newcommand{\GRon}{\operatorname{GR}^{\circ}_n}
\newcommand{\Fr}{\operatorname{Fr}}
\newcommand{\GL}{\operatorname{GL}}
\newcommand{\gl}{\mathfrak{gl}}
\newcommand{\SL}{\operatorname{SL}}
\newcommand{\ff}{\footnote}
\newcommand{\ot}{\otimes}
\def\Ext{\operatorname {Ext}}
\def\Hom{\operatorname {Hom}}
\def\Ind{\operatorname {Ind}}
\def\bbZ{{\mathbb Z}}
\newcommand{\mOp}{\on{MOp}}

\newcommand{\AutO}{\on{Aut}\Oo}
\newcommand{\Der}{{\on{Der}\,}}
\newcommand{\DerO}{{\on{Der}\Oo}}
\newcommand{\AutK}{\on{Aut}\K}
\newcommand{\SShv}{\on{SShv}}

\newcommand{\nc}{\newcommand}

\nc{\cont}{\on{cont}} \nc{\rmod}{\on{mod}} \nc{\Mtil}{\widetilde{M}}
\nc{\ol}{\overline} \nc{\wb}{\overline} \nc{\wt}{\widetilde}
\nc{\wh}{\widehat} \nc{\sm}{\setminus} \nc{\mc}{\mathcal}
\nc{\mbb}{\mathbb} \nc{\Mbar}{\wb{M}} \nc{\Nbar}{\wb{N}}
\nc{\Mhat}{\wh{M}} \nc{\pihat}{\wh{\pi}} \nc{\JYX}{\cJ_{Y\leftarrow
X}} \nc{\phitil}{\wt{\phi}} \nc{\Qbar}{\wb{Q}}
\nc{\DYX}{\D_{Y\leftarrow X}} \nc{\DXY}{\D_{X\to Y}}
\nc{\dR}{\stackrel{\bbL}{\underset{\D_X}{\ot}}}
\nc{\Winfi}{\cW_{1+\infty}} \nc{\K}{{\mc K}} \nc{\Kx}{{\mc
K}^{\times}} \nc{\Ox}{{\mc O}^{\times}} \nc{\unit}{{\bf \on{unit}}}
\nc{\boxt}{\boxtimes} \nc{\xarr}{\stackrel{\rightarrow}{x}}
\nc{\Gamx}{\Gamma_{-}^{\times}} \nc{\Gap}{\Gamma_+}
\nc{\Gamxn}{\Gamma_{-}^{\times}} \nc{\Gapn}{\Gamma_{+,n}}

\nc{\At}{\on{At}} \nc{\Dlog}{{\mathcal D}_{\on{log}}} \nc{\Ereg}{\G}
\nc{\Jac}{\on{Jac}} \nc{\Vplus}{V^{+}} \nc{\Vc}{V^{\vee}}
\nc{\af}{{\mathfrak a}} \nc{\h}{{\mathfrak h}}
\nc{\BunDP}{\on{Bun}_{\boldp(\cD)}} \nc{\SSh}{\CM}
\nc{\Enat}{E^{\natural}} \nc{\Enatbar}{\overline{E}^{\natural}}
\nc{\Id}{\on{Id}} \nc{\Fin}{F^{-1}} \nc{\Dhat}{\wh{D}}
\nc{\Of}{{\mathbb{O}}}\nc{\Pp}{{\mbb P}} \nc{\Hitch}{{\on{Hitchin}}}
\nc{\framing}{T}\nc{\ul}{\underline}

\def\hangfour{\hangindent=.4in\hangafter=1}
\def\hangone{\hangindent=.1in\hangafter=1}
\def\hangtwo{\hangindent=.2in\hangafter=1}

\newcommand{\Volt}{{\mathbf V}_n}
\nc{\volt}{{\mathbf v}}



\title{$\D$-Bundles and Integrable Hierarchies}
\author{David Ben-Zvi}
\address{Department of Mathematics\\University of Texas at Austin\\Austin, TX 78712 USA}
\email{benzvi@math.utexas.edu}
\author{Thomas Nevins}
\address{Department of Mathematics\\University of Illinois at Urbana-Champaign\\Urbana, IL 61801 USA}
\email{nevins@uiuc.edu}

\date{\today}

\maketitle

\tableofcontents

\section{Introduction}

We study the geometry of $\D$-{\em bundles}---locally projective
$\D$-modules---on algebraic curves, and apply them to the study of
integrable hierarchies, specifically the multicomponent
Kadomtsev-Petviashvili (KP) and spin Calogero-Moser (CM)
hierarchies. We show that KP hierarchies have a geometric
description as flows on moduli spaces of $\D$-bundles; in
particular, we prove that the local structure of $\D$-bundles is
captured by the full Sato Grassmannian. The rational, trigonometric,
and elliptic solutions of KP are therefore captured by $\D$-bundles
on cubic curves $E$, that is, irreducible (smooth, nodal, or
cuspidal) curves of arithmetic genus $1$. We develop a Fourier-Mukai
transform describing $\D$-modules on cubic curves $E$ in terms of
(complexes of) sheaves on a twisted cotangent bundle $\Enat$ over
$E$. We then apply this transform to classify $\D$-bundles on cubic
curves, identifying their moduli spaces with phase spaces of general
CM particle systems (realized through the geometry of spectral
curves in $\Enat$). Moreover, it is immediate from the geometric
construction that the flows of the KP and CM hierarchies are thereby
identified and that the poles of the KP solutions are identified
with the positions of the CM particles. This provides a geometric
explanation of a much-explored, puzzling phenomenon of the theory of
integrable systems: the poles of meromorphic solutions to KP soliton
equations move according to CM particle systems.

\subsection{A Brief History of the KP/CM Correspondence}
We begin with a rough sketch of the history of the problem in
integrable systems that motivated the present work---see the review
articles \cite{Benn,BAMS} for more complete histories and
bibliographies.  Also, an exposition
of the present   work appears in \cite{solitons}, together with a
 leisurely historical discussion and overview.

 In the seminal work \cite{AMM}, Airault, McKean, and
Moser wrote down rational, trigonometric, and elliptic solutions of
the Korteweg-deVries equation and discovered that the motion of
their poles is governed by the Calogero-Moser classical many-body
systems of particles on the line, cylinder, and torus (respectively)
with inverse-square potentials. Krichever \cite{Kr1,Kr2} and the
Chudnovskys \cite{CC} extended this correspondence to the
meromorphic solutions of the KP equation, where it becomes an
isomorphism between the phase spaces of generic rational (decaying
at infinity), trigonometric, and elliptic KP solitons and the
corresponding Calogero-Moser systems. In fact, Krichever wrote the
elliptic CM systems in Lax form with elliptic spectral parameter
and showed that the generic elliptic KP solutions are
algebro-geometric solutions associated to line bundles on the
corresponding branched covers of elliptic curves. This approach is
extended in \cite{Kr Lax Vect}, where the elliptic CM systems are
considered as part of a general Hamiltonian theory of Lax operators
on algebraic curves, and \cite{AKV}, where this field analog of the
elliptic CM system is related to the KP hierarchy. A detailed
algebro-geometric study of the CM spectral curves, as geometric
phase spaces for the elliptic KP/CM systems, was undertaken by
Treibich and Verdier \cite{TV,TV2} and lead, in particular, to a
complete classification of elliptic solutions of the Korteweg-deVries (KdV)
 equation.

The rational KP/CM correspondence was explored and deepened by
Shiota \cite{Sh} and, especially, by Wilson \cite{Wilson CM} (see
\cite{Wsurv3} for a review). Shiota identified all the higher flows
of the KP hierarchy on generic rational solutions with the higher
hamiltonians of the rational CM particles. Wilson extended the
correspondence away from generic solutions by allowing collisions
of Calogero-Moser particles. In \cite{Wilson bispectral}, Wilson had
identified the completed phase space of the rational KP hierarchy
with an {\em ad\`elic Grassmannian} (which appears independently in
the work of Cannings and Holland \cite{CH ideals} classifying ideals
in the Weyl algebra of differential operators on the affine line).
In \cite{Wilson CM}, Wilson gives an explicit formula that defines a
point of the adelic Grassmannian from the linear algebra data
describing the rational CM space.  He then proves by direct
calculations that this map extends continuously to the completed
phase spaces and takes the CM flows to the KP flows.

The emergent
relation between CM spaces and the Weyl algebra was explored in the
papers \cite{BW automorphisms, BW ideals} by Berest and Wilson and
extended in \cite{BGK1,BGK2} by Baranovsky, Ginzburg, and Kuznetsov.
Inspired by ideas of Le Bruyn \cite{Le Bruyn},
these authors used noncommutative geometry to identify the rational
CM space (and its spin versions) with the isomorphism classes of
ideals in (and generally torsion-free modules over) over the Weyl
algebra, and thus with Hilbert schemes of points on a noncommutative
surface (see also \cite{NS}). This is closely related to the study
of noncommutative instantons on $\C^2$ (see, for example,
\cite{Nekrasov-Schwarz,KKO,BraNe}).

An independent development of great significance for the present
project was the work
of Nakayashiki and Rothstein \cite{N1,N2,Ro1} and especially
\cite{Ro2}. In these works, the Fourier-Mukai transform on Jacobians
of smooth curves is applied to construct $\D$-modules and, through
Sato's $\D$-module description of the KP hierarchy, to describe the
Krichever construction of KP solutions from line bundles on the
curve.

In this paper, we provide a direct relation between arbitrary
meromorphic (rational, trigonometric, and elliptic) solutions of KP
(and its multicomponent generalization) and the noncommutative
geometry of modules over differential operators. We show that this
geometric description of the KP hierarchy is directly identified,
through a Fourier-Mukai transform, with the geometric description
of the completed (spin generalizations of the) Calogero-Moser
systems. This provides a uniform conceptual description of the KP/CM
correspondence in its most general setting.

We now describe the contents of the present paper in detail.

\subsection{Background} In Section \ref{cubics}, we include a brief
reminder on our backdrop, the family of Weierstrass cubic curves.
These fall into three types, which parallel the three flavors of
many-body systems: smooth elliptic curves (the elliptic case), the
projective line with one node (trigonometric case), and the
projective line with one cusp (rational case). All the constructions
of this paper work over the universal family of cubic curves;
however, for readability, we usually work
 over individual cubic curves. Over each such curve $E$,
we consider a ruled surface $\Enatbar$, with a section $E_{\infty}$
whose complement is the (unique nontrivial) affine bundle $\Enat$
over $E$. In the elliptic case, $\Enat$ can be described
alternatively as the universal additive extension of $E$, as the
space of line bundles with flat connections, or as the home of the
Weierstrass $\zeta$-function.

In Section \ref{CM section}, we discuss the spin CM
 particle systems associated to cubic curves following
\cite{spin} (see also the overview \cite{Nekrasov survey}). These
are Hamiltonian systems describing a system of $n$ particles living
in (the smooth part of) $E$ and completed so as to allow collisions
of particles. The usual spin CM particles carry spins in an
auxiliary vector space $\C^k$. We consider a more general version,
in which the spins take value in a length $k$ torsion sheaf
$\framing$ on (the smooth locus of) $E$; the usual setting
corresponds to $\framing$ being a length $k$ skyscraper at the
marked point of $E$. We present the general spin CM systems in
``geometric action-angle" variables, that is, in terms of line bundles
on spectral curves. More precisely, in \cite{spin} it is shown
(following in the spinless case the works
\cite{Kr2,TV2,DW,GorNe,Nekrasov}, among others) that the phase spaces
of the generalized spin CM systems are identified with spaces of
torsion-free sheaves supported on curves in $\Enatbar$  (the CM
spectral sheaves). The relation between the particle and spectral
curve descriptions is given by a Fourier transform. Specifically,
the positions of CM particles are recovered from the spectral curve
description as the finite support of the Fourier-Mukai transform of
the projection to $E$ of the spectral sheaf. The CM hamiltonian
flows are explicitly identified with a hierarchy of flows that
preserve the curve and ``tweak" the sheaves along the intersection
with the curve $E_{\infty}\subset \Enatbar$.

In fact, we may tweak sheaves by arbitrary meromorphic
endomorphisms near the curve $E_{\infty}$, giving rise to a natural
Lie algebroid acting on CM spectral sheaves, the CM algebroid. For
generic framings---for example, in the spinless case $\framing=\Oo_b$)---these
endomorphisms consist of several copies of the Lie algebra of Laurent
series, so we get several commuting hierarchies of flows (whose
labeling depends on a choice of coordinate on the spectral curve).
In general, the algebroid is isomorphic to a sum of several twisted
loop algebras $L\gl_k$ corresponding to the rank of a sheaf on the
components of its support. To give names to our flows, we may
``Higgs" our spectral sheaves by picking a distinguished such
endomorphism, breaking down the symmetry to give a commuting family
of flows. Higgsed spectral sheaves carry a family of flows labeled
by natural numbers, corresponding to tweaking the sheaf by powers of
the endomorphism. We refer to this enhancement as the ``Higgsed CM
hierarchy.''

In Section \ref{KP section}, we review the multicomponent
generalizations of the KP hierarchy, including their definition as
flows on formal matrix-valued microdifferential Lax operators, and
Sato's description as flows on the big cell of an
infinite-dimensional Grassmannian. In particular, we are interested
in Sato's reinterpretation of KP wave operators in terms of the free
modules they generate over the ring $\D$ of differential operators
with Taylor series coefficients (that is, as $\D$-modules on the disc).
Sato thereby identifies the big cell of the Grassmannian with free
$\D$-modules embedded in the algebra $\cE$ of microdifferential
operators with Taylor coefficients.

\subsection{$\D$-Bundles} In Section \ref{D bundles}, we introduce
$\D$-bundles on smooth algebraic curves $X$ (and their natural
extension to the case of cubic curves). A $\D$-bundle is a
torsion-free module over the sheaf of differential operators on a
curve. Examples include locally free $\D$-modules ($V\ot\D_X$ for a
vector bundle $V$ on $X$), but also ideals in the Weyl algebra
$\D_\aline$, and have been extensively studied in \cite{CH ideals,BW
ideals,BGK1,BGK2,cusps}, among others. It is convenient to consider
$\D$-bundles from the point of view of noncommutative geometry, as
torsion-free coherent sheaves on a noncommutative affine bundle
$T^*_\hbar X$ over $X$. For the point of view of moduli problems, it
is clearly better to consider torsion-free sheaves on a
noncommutative {\em ruled surface} $\ol{T^*_\hbar} X= T^*_\hbar
X\cup X$, which are the {\em framed} $\D$-bundles. For a vector
bundle $V$ on $X$, a $V$-framed $\D$-bundle is a torsion-free sheaf
on $\ol{T^*_\hbar}X$ whose restriction to the curve $X$ at infinity
is identified with $V$. The resulting moduli spaces give
deformations of Hilbert schemes of points on $T^*X$ (in rank $1$)
and of more general moduli of framed torsion-free sheaves; they may
be considered algebraic analogs of spaces of noncommutative
instantons (see, e.g., \cite{KKO,BraNe}). See Section \ref{D
bundles} for precise algebraic definitions in terms of $\D$-modules
equipped with normalized filtrations.

We state our central result on $\D$-bundles on cubic curves
(Theorem \ref{state KP/CM equivalence}) in Section \ref{D bund and
particles} but defer the proof to Section \ref{isomorphism of
moduli}, after we have studied the Fourier-Mukai transform. The theorem
 identifies the moduli spaces
of $\D$-bundles with the phase spaces of the corresponding
Calogero-Moser systems, via a generalized Fourier-Mukai transform.
Fix a cubic curve $E$ and a semistable degree zero vector bundle $V$
on $E$, which corresponds under the Fourier-Mukai transform to a
torsion sheaf $\Vc$ supported on the smooth locus of $E$.

\begin{thm}[Theorem \ref{KP/CM stack equiv}]\label{state isom theorem}
The moduli stack $\BunDP(E,V)$ of $V$-framed $\D$-bundles on a
cubic curve is isomorphic to the $\Vc$-framed spin Calogero-Moser
phase space $\CM_n(E,\Vc)$. The isomorphism identifies the cusps of
a $\D$-bundle with the positions of the corresponding
Calogero-Moser particles.
\end{thm}
\noindent
Note that the usual spin Calogero-Moser spaces are obtained in the
case of trivial framing $V=\Oo_E^k$, $\Vc=\Oo_b^k$.

 The rank $1$
case of this theorem recovers the identification of the space of
ideals in the Weyl algebra with the rational CM space \cite{BW
ideals}.
It also refines the separation
of variables of \cite{GNR}, giving a birational identification
between (spinless) Calogero-Moser spaces and Hilbert schemes of
points, by identifying the full CM spaces with deformed Hilbert
schemes of noncommutative points.  Also, our description of this
identification as a Fourier-Mukai transform establishes the
speculation of \cite{GNR} that separation of variables is a
T-duality.

In Section \ref{adelic Gr}, we relate $\D$-bundles on arbitrary
curves to Wilson's {\em ad\`elic Grassmannian} \cite{Wilson CM}. We
show that, in parallel to the situation for framed torsion-free
sheaves on a ruled surface, framed $\D$-bundles on any curve have a
canonical trivialization away from a finite subset of the curve, the
{\em cusps} of the $\D$-bundle (named following \cite{CH
cusps,cusps}). Wilson's ad\`elic Grassmannian can be identified,
following Cannings and Holland \cite{CH cusps} (see also \cite{BW
automorphisms, BGK1,cusps}), with the moduli space of ({\em
unframed}) $\D$-bundles with generic trivialization (in Wilson's
spinless settings, these are rank 1 $\D$-bundles, while we consider
ones of arbitrary rank). Thus the canonical trivialization off cusps
defines maps from our moduli spaces of framed $\D$-bundles to the
ad\`elic Grassmannians. These maps give a set theoretic
decomposition of the ad\`elic Grassmannians into finite dimensional
moduli spaces. The relevant ``topologies," i.e. notions of {\em
families} of $\D$-modules, are very different, so that the
different moduli of $\D$-bundles become connected in the ad\`elic
Grassmannian. The results of Section \ref{adelic Gr} generalize those
of \cite{Wilson CM,BW
automorphisms,BW ideals,BGK1,BGK2} on ad\`elic Grassmannians and
$\D$-modules while giving precise algebraic meaning to Wilson's
set-theoretic decomposition of the rational rank 1 Grassmannian
into Calogero-Moser spaces. Using the Cannings-Holland--inspired
interpretation of $\D$-bundles on $X$ as coherent sheaves on
cuspidal curves normalized by $X$, we thus also obtain a
construction of solutions of multicomponent KP hierarchies (orbits
in the Sato Grassmannian) from $\D$-bundles.

\subsection{$\D$-Bundles and KP Hierarchies} In Section \ref{Sato},
we relate the local structure of $\D$-bundles to the KP hierarchy.
We generalize Sato's description of the big cell of the Grassmannian
to a description of the full Grassmannian by replacing the free
$\D$-modules in his construction by $\D$-bundles:
\begin{thm}[Theorem \ref{Sato-type thm}]\label{state Sato theorem}
The rank $n$ Sato Grassmannian is isomorphic to the moduli space of
rank $n$ $\D$-bundles $M$ on the disc equipped with an isomorphism
$M\ot_\cD\cE\to \cE^n$ of the microlocalization with the free rank
$n$ module over microdifferential operators $\cE$ on the disc.
\end{thm}
Equivalently, the Grassmannian parametrizes $\D$-{\em lattices},
which are torsion-free finite rank $\D$-submodules of $\cE^n$
which generate $\cE^n$.\footnote{In fact, in Section \ref{D and d} we prove a
precise algebraic statement showing that Sato's Grassmannian represents
an appropriate functor of flat families of $\D$-modules.}
 This
construction is interpreted as a noncommutative version of the
Krichever construction, where we replace line bundles on a curve
trivialized near a point by $\D$-line bundles on the disc with a
microlocal trivialization and obtain in this way the full Sato
Grassmannian.

The $\D$-bundle description of the Sato Grassmannian also gives
rise to a geometric reformulation of KP Lax operators, the {\em
micro-opers}, introduced in Section \ref{microopers}. Micro-opers
are $\D$-bundles equipped with a microlocal endomorphism, which may
be considered as a flat connection, satisfying a strong form of
Griffiths transversality. A micro-oper on a curve $X$ canonically
determines (and is determined by) a matrix Lax operator away from
the cusps of the underlying $\D$-bundle. In other words, the cusps
of a $\D$-bundle provide a natural geometric description of the
poles of matrix KP Lax operators. Thus, micro-opers are perfectly
suited for the study of meromorphic solutions of multicomponent KP
hierarchies, which are our primary motivation. Micro-opers are the
analogs for setting of (multicomponent) KP equations of the {\em
opers} of Beilinson-Drinfeld \cite{BD opers}, or more precisely of
the {\em affine opers} of \cite{BF}, for the setting of
Drinfeld-Sokolov (generalized KdV) equations.

Micro-opers carry a hierarchy of flows that the KP hierarchy
on Lax operators on the disc. The flows are simply given by
modifying (the transition functions of) the underlying $\D$-bundle
by powers of the microlocal endomorphism. These flows are part of a
natural Lie algebroid on the space of $\D$-bundles, the KP
algebroid. This algebroid consists of microlocal deformations of
$\D$-bundles, that is, deformations coming from endomorphisms of their
microlocalizations, acting by deforming ``transition functions along
the curve at infinity.'' Thus a micro-oper structure on a
$\D$-bundle can be considered as a choice of element of this Lie
algebroid, with fixed polar part at $E_\infty$.

It is immediate from the construction of the Fourier-Mukai
identification between $\D$-bundles and CM spectral sheaves that
the corresponding Lie algebroids are identified. Informally, the
Fourier-Mukai transform identifies the (commutative and
noncommutative) ruled surfaces on which CM spectral sheaves and
$\D$-bundles live, and respects the sections at infinity. Since
both hierarchies are given by modifying sheaves along the respective
sections at infinity, the CM and KP flows are intertwined by the
Fourier-Mukai transform.

\begin{thm}[Theorem \ref{compatibility}]
Let $\bF:\BunDP(E,V)\to\CM_n(E,\Vc)$ denote the Fourier-Mukai
isomorphism of framed $\D$-bundles and framed CM spectral sheaves.
\begin{enumerate}
\item $\bF$ identifies the KP Lie algebroid with the CM Lie algebroid.

\item $\bF$ lifts to an isomorphism of the moduli stack of $V$-framed micro-opers
with the moduli stack of $\Vc$-framed Higgsed CM spectral sheaves,
identifying the multicomponent KP hierarchy on micro-opers with the
Higgsed CM hierarchy.
\end{enumerate}
\end{thm}

In the rank one case this theorem gives a strong form of the
correspondence between CM particles and meromorphic KP solutions:

\begin{corollary}[Rank one reformulation]
The completed phase spaces of the rational, trigonometric, and
elliptic (spinless) Calogero-Moser systems are identified with the
moduli spaces of rational, trigonometric ,and elliptic (rank one) KP
Lax operators (taken up to change of coordinate in $\del\inv$).
This isomorphism
identifies poles of Lax operators with positions of Calogero-Moser
particles and identifies the KP and CM hierarchies.
\end{corollary}

To paraphrase, the tweaking flows on CM spectral sheaves are simply
the expression in action-angle coordinates of the matrix KP
hierarchy. The location of the CM particles and the poles of the
matrix Lax operator are both described by the cusps of a
$\D$-bundle, and all multicomponent KP flows are described
spectrally by different tweakings of spectral sheaves near the curve
$E_\infty$. For generic framings, in particular in the rank $1$
case, these flows form a commuting family that is identified with
Laurent series on the spectral curve near its intersection(s) with
$E_\infty$.\footnote{In other words, the choice of micro-oper structure or
Higgsing is unique up to formal changes of coordinates in
$\del\inv$ or equivalently on the spectral curve.} In other
situations, the possible flows form a noncommuting family, although by
choosing a micro-oper or Higgs structure we pick out abelian
subalgebras of flows. This choice parallels the well-known
dependence of matrix KP or Drinfeld-Sokolov hierarchies on choices
of Heisenberg subalgebras of loop groups.

This theorem demystifies and generalizes the results of
\cite{AMM,CC,Kr1,Kr2,TV,TV2,Sh,Wilson CM,BBKT,T matrix} on the
relation between CM particles (and their collisions) and the
meromorphic solutions of KP. In particular, the theorem extends the
results of Wilson\cite{Wilson CM} in two directions: over the
family of cubic curves (that is, to the trigonometric and elliptic
cases), and to higher rank (multicomponent/spin setting, as
predicted in \cite{Wilson CM,Wsurv3} in the rational trivially
framed case). In particular, in the higher rank case we have the
choice of framing $V$, with the trivial case $V=\Oo^k$ corresponding
to the ordinary spin CM system, but other framings correspond to
integrable systems with very different geometry and dynamics---in
fact, generic framings $V$ give rise to a much simpler (abelian)
hierarchy.

It is interesting to note that we now have (even in the rank $1$
case) two independent relations between $\D$-bundles and KP
solutions---our construction involving micro-opers and the
original construction of Wilson \cite{Wilson CM} and \cite{CH
cusps}, whereby the ad\`elic Grassmannian on a curve $X$
parametrizes Krichever data for algebro-geometric solutions of KP
associated to $X$ together with all of its cuspidal quotients
(curves obtained by adding cusps to $X$). While there is no general
relation between the two constructions, in genus zero (where Wilson
was working) they are identified by the geometric Fourier transform
on $\aline$ (the natural auto-equivalence of $\D$-modules on
$\aline$ exchanging multiplication and differentiation). In other
words, the geometric Fourier transform induces a self-map of the
appropriate moduli of $\D$-bundles on $\aline$ that exchanges the KP
algebroid (microlocal deformations of $\D$-bundles, near
``$\del\inv=0$") with the algebroid deforming $\D$-bundles near
$z\inv=0$. This gives a simple geometric explanation of (and
multicomponent generalization of) the bispectral involution on
rational solutions of KP, studied in detail in \cite{Wilson
bispectral, BW automorphisms} (see Section \ref{bispectrality}). We
plan to apply this point of view further in the less-explored
setting of difference modules.

\subsection{Fourier Transform and Moduli of $\D$-Bundles} In
Section \ref{Fourier}, we recall the Fourier-Mukai autoequivalence
of the bounded derived category of coherent sheaves on a cubic curve
\cite{BuK}. In Section \ref{FT for D}, we extend to singular cubic
curves the Fourier-Mukai transform for $\D$-modules, discovered by
Laumon \cite{La} and Rothstein \cite{Ro2} for abelian varieties (and
extended by Polishchuk and Rothstein \cite{PRo} to general
$D$-algebras).  Namely, the Fourier-Mukai transform on a cubic curve
identifies $\D$-modules with (complexes of) coherent sheaves on the
surface $\Enat\to E$, the ``twisted log cotangent bundle" of $E$:
\begin{thm}[Theorem \ref{cubic Laumon-Rothstein}]
The bounded derived category of coherent $\D_{log}$-modules
on a cubic curve $E$ is equivalent to the bounded derived category
of coherent sheaves on $\Enat$.
\end{thm}
We also show, following \cite{PRo}, that this Fourier transform
``compactifies,'' identifying the derived category of modules over
the Rees algebra of $\D$ (coherent sheaves on the noncommutative
ruled surface $\ol{T^*}_\hbar E$) with that of coherent sheaves on
the ruled surface $\Enatbar$, and respects microlocalization
(restriction near $E_\infty$).

In Section \ref{isomorphism of moduli}, the technical heart of the
paper, we apply the Fourier-Mukai transform to prove
Theorem \ref{state isom theorem}, describing the moduli spaces of
$\D$-bundles on $E$ which are framed by a semistable vector bundle
$V$ on $E$ of degree zero. Let $V^\vee$ denote the torsion coherent
sheaf on $E$ Fourier dual to $V$. We show that $V$-framed
$\D$-bundles are sent by the Fourier transform to coherent sheaves
on $\Enatbar$ (rather than complexes), of pure $1$-dimensional
support, whose restriction to the curve $E_\infty\subset \Enatbar$
is identified with $V^{\vee}$. More precisely, we establish an
equivalence of the corresponding moduli stacks. In the cuspidal
case, the moduli spaces of spectral sheaves are the rational spin
Calogero-Moser spaces, which are identified with certain Nakajima
quiver varieties. In this case, our result recovers the theorems of
\cite{BW ideals} classifying ideals in the Weyl algebra (and more
general framed $\D$-bundles on $\aline$) in terms of quiver data.

In Section \ref{D and d}, we prove Theorem \ref{state Sato theorem},
establishing that the Sato Grassmannian represents the functor of
flat families of {\em $\D$-lattices} on the disc, which are
finitely-generated $\D$-submodules of the ring of
microdifferential operators on the disc. This extends Sato's
description of the big cell and refines it from a set-theoretic to
a scheme-theoretic statement (i.e., to families).

Finally, in Section \ref{D algebras}, we carefully explain some technical
aspects of
$\D$-bundles on (possibly singular) cubic curves. More precisely, we
discuss the sheaf $\D_{log}$ of ``logarithmic" differential operators
on a cubic curve $E$, which is generated by the invariant vector
field for the group structure on the smooth part of $E$. We also
discuss filtered $\D_{log}$-modules and their properties.

\subsection{Further Directions} In \cite{solitons2}, we extend the
picture developed in this paper by replacing differential operators
by difference operators. We describe the nonabelian Toda lattice
hierarchies in terms of difference modules, in particular realizing
meromorphic (rational, trigonometric, and elliptic) solutions in
terms of difference modules on cubic curves. A modified version of
the Fourier-Mukai transform then identifies the moduli spaces of
difference modules with spaces of spectral sheaves on a ruled
surface over the corresponding cubic curves. These spectral sheaves
realize the (generalized spin) Ruijsenaars-Schneiders (RS)
relativistic many-body systems, which are deformations of the
corresponding CM systems. This gives a general geometric picture of
the Toda/RS correspondence \cite{KrZab}. This point of view also has
applications to the difference version of bispectrality.

In \cite{W}, we study the factorization (or ``vertex algebra space")
structure \cite{chiral} on the ad\`elic Grassmannian. This structure
is shown to encapsule both (infinitsimally) the
$\cW_{1+\infty}$-symmetry of the KP hierarchy and (globally) the
B\"acklund transformations.

In \cite{BGN} (joint with V. Ginzburg), we develop a simple
algebraic description of $\D$-bundles on arbitrary curves, which
generalizes the quiver description of Calogero-Moser spaces (the
$\aline$ case). We use Koszul duality to describe moduli spaces of
$\D$-modules on curves as twisted cotangent bundles to moduli of
complexes of sheaves on the curve, extending the relation between
vector bundles with connections and Higgs bundles. The resulting
picture of $\D$-bundles has the advantage of being concrete, local,
and functorial.

\subsection{Contents}
We will now describe the contents of the paper section by section.
Section \ref{background section} recollects some necessary
background material. Section \ref{D bundles} is the heart of the
paper, containing the main results on $\D$-bundles and their
application to the KP/CM correspondence. Sections \ref{Fourier}
through \ref{D algebras} contain the technical tools used in Section
\ref{D bundles}. Specifically, Section \ref{Fourier} discusses the
main tool, the extended Fourier-Mukai transform for $\D$-modules on
cubic curves. Section \ref{isomorphism of moduli} applies this
transform to prove an isomorphism of moduli stacks of $\D$-bundles
and spectral sheaves. Section \ref{D and d} provides the proof for
the $\D$-bundle description of the Sato Grassmannian. Finally
Section \ref{D algebras} is an appendix containing needed material
about differential operators and $D$-algebras on cubic curves.

\subsection{Acknowledgments}
The authors are grateful to R. Donagi, E.  Frenkel, D.
Gaitsgory, V. Ginzburg, T. Pantev, and M. Rothstein for helpful
conversations.

The first author was supported in part by an
 NSF postdoctoral fellowship at the University of Chicago and an
MSRI postdoctoral fellowship, as well as by NSF CAREER grant
DMS-0449830 at the University of Texas. The second author was
supported in part by an NSF postdoctoral fellowship at the
University of Michigan and an MSRI postdoctoral fellowship, as well
as by NSF grant DMS-0500221.

\section{Background Material}\label{background section} In this
section, we review some features of the geometry of cubic curves,
describe the spin Calogero-Moser systems, and review the
multicomponent KP hierarchies.

\subsection{Background on Cubic Curves}\label{cubics}

\subsubsection{One-Dimensional Groups and Cubic Curves}
Connected one-dimensional complex groups $\G$ fall into three
classes: the additive group $\C$, the multiplicative group $\Cx$,
and the one-parameter family of elliptic curves. These cases fall
under the monikers {\em rational, trigonometric,} and {\em elliptic}
according to the type of functions on the universal cover $\C$ which
correspond to meromorphic functions on $\G$.

A parallel classification applies to reduced and irreducible cubic
plane curves, which we will also call (slightly abusively) {\em
Weierstrass cubic curves}: reduced and irreducible complex
projective curves $E$ of arithmetic genus one with a nonsingular
marked point $b\in E$. (See \cite{FriedMorgIII} for a detailed
discussion of cubic curves and bundles on them.)  The map $q \mapsto
\theo(q-b)$ defines an isomorphism from the smooth locus of $E$ to
the Jacobian $\on{Pic}^0(E)$; in particular, the smooth locus, which
we denote by $\G$, is equipped with a structure of one-dimensional
group, with $b$ as the identity element. The three types are then:
\vspace{.4em}

\hangtwo\noindent {\bf Rational:} $E$ is a cuspidal cubic, and is
isomorphic to the curve $y^2=x^3$. Its normalization $\pline\to E$
collapses $2\cdot\infty$ to a cusp on $E$, and defines a group
structure $\C= \G\subset E$ on the smooth locus.

\vspace{.2em}

\hangtwo\noindent {\bf Trigonometric:} $E$ is a nodal cubic, and is
isomorphic to the curve $y^2=x^2(x-1)$.  Its normalization
$\pline\to E$ identifies two points $0$ and $\infty$ to a node on
$E$, and defines a group structure $\Cx=\G\subset E$ on the smooth
locus.

\vspace{.2em}

\hangtwo\noindent
{\bf Elliptic:} $E$ is a smooth elliptic curve (in
particular a group), and may be described by an equation of the form
$y^2=x^3+a x+b$ with $\Delta=-16(4a^3+27b^2)\neq 0$.

\vspace{.4em}

\noindent Let $\infty$ denote the singular point in the rational and
trigonometric cases.  Then $E$ is identified with its own
compactified Jacobian, the moduli space of torsion free sheaves of
rank one and degree zero on $E$. The singular point corresponds to
the unique rank 1 degree $0$ torsion-free sheaf which is not
locally free, namely the modification ${\mathfrak m}_\infty(b)$ of
the ideal sheaf of $\infty$.

\subsubsection{Differential Operators on Cubic Curves}\label{Dlog introduced}
The group variety $\G$ acts on $E\cong \overline{\on{Jac}}(E)$, and
admits a unique nonzero invariant vector field $\del$ on $E$ up to a
scalar. If $E$ is smooth this is the usual translation-invariant
vector field, which is constant in a global analytic coordinate.
Writing the singular cubics in terms of their normalization
$\pline$, a choice of $\del$ is represented by
$\frac{\partial}{\partial z}$ in the cuspidal case (vanishing to
order $2$ at $\infty\in\pline$) and by $z\frac{\partial}{\partial
z}$ in the nodal case (vanishing to order one at $0,\infty$).

We will abuse notation to denote the sheaf $\Oo_E\cdot\del$ (i.e.
the action algebroid of $\G$) by $\cT_E$.  Note that in the nodal
case this is the log tangent bundle of $E$; the dual sheaf will be
denoted by $\Omega_E$. The total space of $\Omega_E$ (which is
isomorphic to $E\times \C$) will be denoted $T^*E$. Both $\cT_E$ and
$\Omega_E$ are trivial line bundles on $E$; however, they are not
trivial over families of cubic curves, i.e. one can't fix a $\del$
(or differential) uniformly for all cubic curves, hence we will
maintain this distinction when necessary. A cubic curve with a
choice of $\del$ gives a {\em differential curve} (curve with choice
of global vector field, in analogy with differential fields). In
fact the constructions of this paper will not require a differential
structure (we only use the canonical line $\C\del$).

The sheaf $\cT_E$ acts by derivations of $\Oo_E$, hence it embeds
in the sheaf of (Grothendieck) differential operators on $E$.
\begin{defn}
$\Dlog$ is the sheaf of $\theo_E$-subalgebras of $\D_E$ generated by
$\C\del$.
\end{defn}
We will usually abuse notation and terminology by denoting $\D_E$
the sheaf of differential operators generated by $\cT_E$ and
$\Oo_E$, and referring to it as the sheaf of {\em log differential
operators} on $E$ (which it is in the smooth and nodal cases). See
Section \ref{D algebras} for more on $\Dlog$.

\subsubsection{The Surface $\Enat$}\label{about Enat}
In this section we discuss a surface $\Enat$ which figures
prominently in the study of both Calogero-Moser systems and the
Fourier-Mukai transform.

We first consider the smooth case. Thus fix an elliptic curve $E$.
Let $A$ denote the {\em Atiyah bundle} on $E$, that is, the unique
nontrivial extension of $\Oo_E$ by itself (up to isomorphism),
\begin{equation}\label{Atiyah extension}
0\rightarrow \theo \rightarrow A \rightarrow \theo\rightarrow 0.
\end{equation}

\begin{defn}
We set $\Enatbar={\mathbf P}(A)$.  
The algebraic surface $\Enat$ is the
complement of the section $E_\infty={\mathbf P}(\Oo)\cong E$ of the
projectivization of the Atiyah bundle,
$\Enat = {\mathbf P}(A)\setminus E_\infty$.
\end{defn}
 The resulting
surface $\Enat$ (which is a Stein manifold, but not an affine algebraic
manifold, when $E$ is an elliptic curve) is the unique (up to isomorphism)
nontrivial torsor over $\theo_E$.
In classical terms, the
surface $\Enat$ may be viewed as the receptacle for the Weierstrass
$\zeta$-function of $E$, see \cite{spin} for a discussion.

To fix $A$ and $\Enat$ canonically we set $\Enat=\on{Conn}\Oo(b)$,
the sheaf of connections on the line bundle $\theo(b)$;
 this is a {\em twisted cotangent bundle} of $E$ \cite{BB}
(i.e. an affine bundle for $\Omega_E\simeq \Oo_E$ with compatible
symplectic structure). Let $\cA$ denote the pushforward of
$\Oo_{\Enat}$ to $E$, i.e. the algebra of functions on the fibers of
$E$. Thus $A=\cA_{\leq 1}$, the subsheaf of affine functions on
$\Enat$, is isomorphic to the Atiyah bundle $A$. The sheaf $A$ is
canonically an extension of $\cT_E$ (which is isomorphic to $\Oo_E$)
by $\Oo_E$, and is isomorphic as $\Oo_E$-module to $A=\D_{\leq
1}(\Oo_E(b))$, the sheaf of differential operators of order at most
one acting on the line bundle $\Oo(b)$. Concretely, the sheaf $A$
lies in between
$$\Oo_E\oplus\cT_E(-b)\subset A \subset \Oo_E\oplus \cT_E(b),$$ and
$A$ is generated (in the canonical local coordinate near $b$) by
$\Oo_E\oplus\cT_E(-b)$ and the section $\del-\frac{1}{z}$. It is
useful to note that the fiber $F_b\subset \Enat$ over the basepoint
$b$ is canonically identified with the cotangent fiber to $E$ at
$b$.

It is also well-known that the twisted cotangent bundle $\Enat$ is
canonically identified with the universal additive extension of $E$,
which is identified with the moduli space of line bundles with flat
connection on $E$. This isomorphism is uniquely characterized as an
isomorphism of torsors over the cotangent bundle preserving
basepoints in the fiber over the identity $b$ (for flat connections
the basepoint is the trivial connection).

\subsubsection{$\Enat$ for Singular Cubics}

The definition and properties of $\Enat$ extend naturally to general
cubics $E$. For any Weierstrass cubic $E$ we have \bd
\Ext^1_E(\theo,\theo)= H^1(E,\theo) \cong \C. \ed So $E$ has a
unique nontrivial extension $A$ of $\theo_E\simeq \cT_E$ by
$\theo_E$, up to isomorphism. (Recall that $\cT_E$ is the subsheaf
of the tangent sheaf generated by the $\G$-action.) We again fix
$A=\D_{\leq 1}(\Oo_E(b))$ (here as elsewhere $\D$ denotes the sheaf
of log differential operators). Let $\Enatbar \overset{\on{def}}{=}
\bproj(\Sym A)$ denote the associated ruled surface, and $p:\Enatbar
\rightarrow E$ denote the projection map.

The quotient map $A\twoheadrightarrow \theo_E$ defines a section $s:
E\rightarrow \Enatbar$; we write $E_\infty = s(E)$ and refer to it
as the {\em section at infinity}.  Note that every other section of
$\Enatbar$ has nonempty intersection with $E_\infty$ since the
sequence \eqref{Atiyah extension} is nonsplit.  However, unlike the case of
smooth $E$, if $E$ is singular there are curves in $\Enatbar$ that fail to
intersect $E_\infty$: for example, the normalization of $E$ embeds in
$\Enat$ since the extension $A$ splits when pulled back to $\pline$.
 The surface $\Enat
\overset{\on{def}}{=} \Enatbar\smallsetminus E_\infty$ is called the
{\em twisted (log) cotangent bundle} of $E$; it is the nontrivial
torsor over $\Omega_E$ given by the nonzero class (up to scale) in
$H^1(\Omega_E)$. Canonically, it is given by the space of
log-connections (liftings of $\cT_E$) on the line bundle $\Oo(b)$.
(See \cite{spin} for the relation of $\Enat$ to the Weierstrass
$\zeta$-function of $E$.) In the context of the Fourier-Mukai
transform, we will identify $\Enat$ in Section \ref{twisted log}
with the moduli space of rank one torsion-free sheaves on $E$ with
log connection. Again this isomorphism is uniquely characterized as
an isomorphism of torsors over the (log) cotangent bundle preserving
basepoints in the fiber over the identity $b$.

\subsection{Calogero-Moser Systems and Spectral Curves}\label{CM
section}

In this section we discuss the spin Calogero-Moser system
\cite{GH}, summarizing the detailed treatment in \cite{spin} (to
which we refer for details and more complete references).
 See also \cite{solitons} and \cite{Nekrasov survey} for reviews of the usual (spinless)
complex Calogero-Moser system following \cite{KKS,Wilson
CM,Nekrasov}.

As in the previous section, we let $\G$ denote a one-dimensional
complex group and $E$ the corresponding cubic curve. The
$k$-spin $n$-particle Calogero-Moser system is a Hamiltonian dynamical
system describing $n$ identical particles moving
on $\G$, each equipped with a vector in
 the auxiliary $k$-dimensional vector space $\C^k$; this vector is known as
the {\em spin} of the particle. The system is then naturally
described in terms of the positions $q_i\neq q_j$ in $\G$ of $n$
distinct particles, momenta $p_i\in\C$ of the particles, spin
vectors $v_j\in \C^k$, and a collection of covectors $u_i\in
(\C^k)^*$. Let $f_{ij}=u_i(v_j)\in\C$ be the contraction of the
$i$th covector with the $j$th vector. The Hamiltonian for the spin
Calogero-Moser system is given by
$$H=\frac{1}{2}\sum_{i=1}^n p_i^2 + \sum_{i<j} f_{ij}f_{ji}U(q_i-q_j).$$
Here the potential energy function $U$ on $\G$
has a double pole at the origin, and is given (as a function on
$\C$, the universal cover of $\G$) by
\begin{equation*}
{\bf Rat: \;}U(q)=\frac{1}{q^2},\hspace{3em} {\bf
Trig:\; }U(q)=\frac{1}{\sin^2(q)}, \hspace{3em} {\bf
Ell:\; }U(q)=\wp(q)
\end{equation*}
where $\wp(q)$ is the Weierstrass $\wp$-function attached to the
elliptic curve $E$. The usual (spinless) Calogero-Moser particle
system is recovered in the case $k=1$ with all $f_{ij}=1$.

The (spin) Calogero-Moser systems have a variety of different
group-theoretic and geometric descriptions, all of which have the
feature that they incorporate collisions of the particles, in other
words a locus where the $q_i$ are no longer distinct. Moreover the
Hamiltonian $H=H_2$ is but one of a family $H_i$ of Poisson
commuting Hamiltonian functions on these completed phase spaces.
These are conveniently realized by writing the system in Lax form
(with spectral parameter) or as a Hitchin system, and taking
(residues of) traces of powers of the Lax operator (respectively,
Higgs field). These descriptions are explained in detail in
\cite{spin}; in this paper, we rely on the description of CM
particles by {\em spectral sheaves}, which are geometric
action-angle variables for the system, and which we will now
review.

\subsubsection{Calogero-Moser Spectral Sheaves}\label{CM spectral
sheaves} Moduli spaces of spectral sheaves (specifically, of line
bundles on curves in a Poisson surface) give a wide class of
examples of integrable systems (see e.g. \cite{DM,Hu}). The
prototypical example of such a setting is the ($GL_n$) Hitchin
system on the moduli space $T^*Bun_n(X)$ of Higgs bundles on a curve
$X$, which can be described as a moduli of torsion free sheaves on
curves in $T^*X$ finite of degree $n$ over $X$. We similarly realize
the spin CM systems in terms of spectral curves on the twisted
cotangent bundle $\Enat$ of $E$, or rather its completion
$\Enatbar$. Specifically, we consider sheaves on $\Enatbar$ with
pure $1$-dimensional support, with a framing condition along the
curve $E_\infty=\Enatbar\sm \Enat$: the restriction of the sheaf to
$E_\infty$ is identified with a $k$-fold skyscraper $\Oo_b^{\oplus
k}$ at the basepoint.

We highlight two aspects of this geometric translation: the
positions of the particles and the role of the spins. We encode the
positions $q_i\in\Ereg$ of the Calogero-Moser particles as follows:
$n$ distinct points in $\Ereg$ define a rank $n$ vector bundle
$W=\bigoplus \Oo(q_i-b)$ on $E$, which is in fact semistable of
degree zero. Conversely, a generic degree $0$ semistable bundle on
$E$ is of the form $W=\bigoplus \Oo(q_i-b)$ for $n$ distinct points
$q_i\in\G\subset E$ (determined up to permutation). More generally,
there is an equivalence between semistable degree zero vector
bundles $W$ on an elliptic curve $E$ and length $n$ torsion coherent
sheaves $W^\vee$ on the smooth locus $\G$ of $E$: $W^\vee$ is given
by the Fourier-Mukai transform of $W$ (see Section \ref{Fourier}
and \cite{FriedMorgminuscule}). The same holds for singular cubics
$E$ if we additionally require that the pullback of $W$ to the
normalization of $E$ is a trivial vector bundle. In the generic case
we have $W^{\vee}=\bigoplus \Oo_{q_i}$. Thus the space of such
bundles $W$ provides a partial completion of the configuration space
of points in $\G$.

The identification between spin CM particles and spectral sheaves
has the feature of identifying the auxiliary space $\C^k$ in which
the spins live with the space of sections of the restriction of the
spectral sheaves to $E_\infty$, which we have normalized to be
$\C^k=\Gamma(\Oo_b^{\oplus k})$. This leads to a natural
generalization of the spin CM system, in which we allow the spins to
live in a general length $k$ coherent sheaf $\framing$ on $E$: we
simply consider sheaves whose restriction to $E_\infty$ is
identified with $\framing$.

\begin{defn}\label{CM spectral sheaves def} Fix a finite length coherent sheaf $\framing$ on $\G\subset E$.
A $\framing$-{\em framed CM spectral sheaf} is a pair $(\cF, \phi)$
consisting of a coherent sheaf $\cF$ on $\Enatbar$ of pure dimension
one, together with an isomorphism $\phi:\cF|_{E_{\infty}}\to
\framing$, satisfying the following two normalization conditions:
\begin{enumerate}
\item[(i)] $W=p_*\cF(-E_\infty)$ is a semistable vector bundle of
degree $0$; if $E$ is singular, we also require that the pullback of
$W$ to the normalization of $E$ is a trivial vector bundle.
\item[(ii)] $\on{deg}(p_*\cF(kE_\infty))=(k+1)\on{deg}(T)$ for all
$k\geq -1$.
\end{enumerate}
The {\em $\framing$-Calogero-Moser space $\CM(E,\framing)=\bigcup
\CM_n(E,\framing)$} is the moduli space of $\framing$-framed CM
spectral sheaves $(\cF,\phi)$ (for which the rank of the vector
bundle $W$ is $n$).

 For the (usual) spin CM system described above,
one chooses $\framing = \theo_b^k$, and we will call objects of
$\CM_n^k(E)=\CM_n(E,\Oo_b^k)$ simply {\em CM spectral sheaves}.
\end{defn}
\begin{remark}
The normalization conditions (i) and (ii) are open conditions on
coherent sheaves of pure dimension one. Condition (i) on the vector
bundle $W$ was discussed above.  Regarding Condition (ii), note that
we have the following:
\end{remark}
\begin{lemma}\label{CM equiv cond}
For a coherent sheaf $\cF$ of pure dimension one on $\Enatbar$ such
that $\cF|_{E_\infty} \cong T$, $\cF$ satisfies Condition (ii) in
Definition \ref{CM spectral sheaves def} if and only if \bd
\on{deg}(p_*\cF(kE_\infty))=(k+1)\on{deg}(T) \; \text{for all $k\gg
0$.} \ed
\end{lemma}

We will see later (Section \ref{Fourier})
 that the full phase space $\CM_n(E,T)$ is
identified by an extended Fourier-Mukai transform with a moduli
space of objects (framed $\D$-bundles) on $E$ whose singularities
are precisely at the points $q_i$, and which are framed by the
semistable vector bundle $T^\vee$ on $E$ Fourier dual to $T$. This
moduli spaces carries a natural integrable system, the
multicomponent KP system.

\subsubsection{Flows on Spectral Sheaves}\label{CM flows}
In \cite{spin} we discuss under the name ``tweaking" a simple
construction of flows on moduli spaces of sheaves from germs of
meromorphic functions. This is a generalization of (the
infinitesimal version of) the action of the Picard group on sheaves
by tensor product. Specifically, let $\cF$ be a sheaf on a variety
$Y$ and $\Sigma\subset Y$ the support of $\cF$. First note that any
class in $c\in \HH^1(\Sigma,\Oo)$ defines a first order deformation
of any sheaf $\cF$---for example, we consider $c$ as a sheaf on
$\Sigma$ over the dual numbers and tensor with $\cF$ to define a
deformation class in $\Ext^1(\cF,\cF)$. If $f$ is the germ near a
point $s\in\Sigma$ of a meromorphic function on $\Sigma$, it defines
a class $[f]\in \HH^1(\Sigma,\Oo)$ (as the connecting map coming
from the inclusion of $\Oo$ into meromorphic germs at $s$), which in
turn defines first-order deformations as above. If we wish to
deform sheaves framed along a divisor $D\subset \Sigma$ (i.e. deform
sheaves fixing their restriction to $D$) we project $f$ to a class
in $\HH^1(\Sigma,\Oo(-D))$ covering $[f]$. In \cite{spin} the flows
of all meromorphic $GL_n$ Hitchin systems are uniformly described in
this way.

The spin CM flows are defined in \cite{spin} in this fashion
(following \cite{TV2} in the rank one case). Specifically, the
choice of a global vector field $\del$ on $E$ (equivalently, a
trivialization of the tangent space $\cT_b$ at $b$) defines a
canonical function $\underline{t}$ on the completion of $\Enat$
along the fiber $F_b$ with first order pole at $E_\infty$. Namely
$\underline{t}=\underline{k}+\pi^* \zeta$, where
$\underline{k}:\Enat\to\pline$ is the composition of the canonical
map from $\Enat=\on{Conn}\Oo(b)$ to $T^*E$ with pole along the zero
fiber, and $\zeta$ is (the Laurent expansion at $b$ of) the
Weierstrass $\zeta$-function. The ring of polynomials in
$\underline{t}$ is canonically identified (independently of choice
of $\del$) with $\Sym\cT_b=\C[\del]$. We thus have the following
definition from \cite{spin}:

\begin{defn}\label{CM flows definition} The Calogero-Moser
hierarchy is the action of the polynomial ring $\Sym\cT_b=\C[\del]$
on $\CM_n^k(E)$ by tweaking of spectral sheaves by powers of
$\underline{t}$ at $b_\infty$.
\end{defn}

\begin{prop} The CM flows are the flows associated to the
CM hamiltonians: the vector field given by the action of $\del^i$ is
given by Poisson bracket with the $i$th  CM hamiltonian $H_i$.
\end{prop}

\subsubsection{The Lie algebroid of Tweakings}\label{Lie algebroids}

Rather than privileging the meromorphic germ $\underline{t}$, we may
consider more general deformations of CM spectral sheaves. Let
$\cA_\cE$ denote the sheaf of functions on the punctured formal
neighborhood of $E_\infty\subset\Enatbar$, i.e. the sheaf of Laurent
series along $E_\infty$. For a spectral sheaf $\cF$ we may consider
its ``microlocalization" $\cF_\cE=\cF\ot_{\Oo_{\Enatbar}}\cA_\cE$,
i.e. the localization of $\cF$ to this deleted formal neighborhood
of $E_\infty$. Since $\cF$ has pure $1$-dimensional support, it
embeds (as $\Oo_{\Enatbar}$-module) into $\cF_\cE$. Given any
endomorphism $\xi\in \End_{\cA_\cE}(\cF_\cE)$, we obtain a canonical
first-order deformation of $\cF$, deforming $\cF\subset\cF_\cE$ to
$(1+\epsilon\xi)\cF\subset \cF_\cE$ (over the dual numbers). More
formally, construct a first-order deformation of $\cF$, $[\xi]\in
\Ext^1(\cF,\cF)$, as the image (under a connecting homomorphism) of
the operation of restricting sections of $\cF$ to $\cF_\cE$:
$$\{s\mapsto \xi(s_\cE) \mbox{ mod }\cF\}\in\Hom(\cF,\cF_{\cE}/\cF)\to \Ext^1(\cF,\cF).$$
If we project instead to $\Ext^1(\cF,\cF(-E_\infty))$ we likewise
obtain deformations of framed spectral sheaves.

\begin{lemma}\label{CM algebroid} The sheaf $\underline{\End}_{\cE}$ over $\CM_n(E,\framing)$
of endomorphisms of the microlocalization of the universal sheaf
(i.e. the sheaf whose fiber at $\cF$ is $\End_{\cA_\cE}(\cF_\cE)$)
has the structure of Lie algebroid, the {\em CM algebroid}, given by
the tweaking action on spectral sheaves. The anchor map sends an
endomorphism to the corresponding $\Ext^1$ class defined above.
\end{lemma}

This algebroid is the spectral sheaf analog of the (transitive) Lie
algebroid on the moduli stack $Bun_G(X)$ of $G$-bundles on a curve
associated to the loop algebra $L{\mathfrak g}$ at a point $x\in X$,
whose fiber at a bundle $\cP$ is the adjoint twist of the loop
algebra, $(L{\mathfrak g})_{\cP}$.

 This
Lie algebroid structure becomes very simple for the typical case of
{\em rank one} spectral sheaves. Namely if $\cF$ is a rank one
torsion-free sheaf on its support, then its endomorphisms are
simply given by functions on the support: $\End(\cF_\cE)$ is the
direct sum of Laurent series on each component of the support of
$\cF_\cE$, i.e. each branch of the spectral curve passing through
$E_\infty$. In this case the action of the algebroid is simply given
by tweaking $\cF$ by meromorphic function germs,  and the resulting
flows may be easily written explicitly as Hitchin hamiltonian flows
following \cite{spin}. This is guaranteed if $\framing=\bigoplus
\Oo_{x_i}$ is a direct sum of skyscrapers at {\em distinct} points
$x_i\in E$, for example in the standard spinless case
$\framing=\Oo_b$.

For general spectral sheaves this produces a noncommutative family
of flows on the moduli space. Namely, the algebroid at $\cF$ is
(noncanonically) isomorphic to a sum $$\End(\cF_\cE)\simeq \bigoplus
L\gl_{k_i}$$ of loop algebras, where for each component $\Sigma_i$
of the spectral curve $\on{Supp}\cF$ near $E_\infty$, $k_i$ is the
rank of $\cF$ on $\Sigma_i$.

To put the resulting flows in the more familiar form of a hierarchy
of flows labeled by natural numbers, we can instead look at spectral
sheaves together with a ``Higgs" structure, breaking the symmetry
down to an abelian family of flows (in analogy with the abelian Lie
algebroid on the moduli space of Higgs bundles, giving the Hitchin
integrable system):

\begin{defn} Fix $\del\in\Gamma(\cT)$. An {\em Higgsed} $\framing$-framed CM spectral sheaf
is a pair $(\cF,\xi)$ consisting of a $\framing$-framed CM spectral
sheaf $\cF$ and a germ $\xi$ of a section of $\End(\cF)(E_\infty)$
at $b_\infty\in E_{\infty}$, with $\xi|_{E_{\infty}}=\del$ as a
section of $\Oo_{\Enatbar}(E_\infty)|_{E_{\infty}}=\cT_E$.
\end{defn}

We may then define the Higgsed CM hierarchy on enhanced spectral
sheaves as the action of the ring $\C[\del]$, where $P(\del)$ acts
on $(\cF,\xi)$ by deforming $\cF$ by $P(\xi)$. Note that the CM
flows defined above are the restriction of the CM hierarchy to the
locus of sheaves Higgsed by $\underline{t}$.

\subsection{The  KP Hierarchy}\label{KP section} Recall that the {\em
Kadomtsev-Petviashvili} (or {\em KP}) {\em equation} is the
following partial differential equation for a function
$u=u(t,x,y)$:\footnote{For compatibility with our later choice of
notation, we have permuted the usual labeling of variables in this
equation.}
\begin{equation}\label{KP eq}
\frac{3}{4}u_{xx}=(u_y-\frac{1}{4}(6uu_t+u_{ttt}))_t,
\end{equation}
which first arose in connection with the study of shallow water
waves. Overviews of this equation and its algebro-geometric
significance, as well as bibliographies, may be found in \cite{Mu,
JM, Ar}. In this section we review the multicomponent KP hierarchies
and their formulation using the Sato Grassmannian.

\subsubsection{Microdifferential Operators}
We begin by recalling the basic algebraic objects, the algebras of
differential and microdifferential operators over the ring of formal
power series: for the rest of this section, $\D\subset \cE$ will
denote the algebras of differential operators and of formal
microdifferential operators (also known as pseudodifferential
symbols) with coefficients in $\C[\![t]\!]$. Specifically, an
element of $\cE$ is a Laurent series
\begin{equation}\label{microdiff}
M=\sum_{N\ll \infty}a_N\del^N,\hspace{3em} a_i\in\C[\![t]\!]
\end{equation} in the
formal inverse $\del\inv$ of the derivation $\del=\del_t$ of
$\C[\![t]\!]$.  Such an element $M$ lies in $\D$ if $a_N=0$ for
$N<0$.

The symbol $\del\inv$ does not, of course, make sense as an operator
on $\C[\![t]\!]$; however, it is possible to give the set $\cE$ of
microdifferential operators an algebra product. The composition in
$\cE$ is determined by the Leibniz rule,
\begin{equation}\label{Leibniz}
\del^n \cdot f=\sum_{i\geq 0}\binom{n}{i} f^{(i)}\del^{n-i},
\end{equation}
where $\binom{n}{i}$ is defined for $n<0$ by taking \bd \binom{n}{i}
= \frac{n(n-1)\cdots (n-i+1)}{i(i-1)\cdots 2\cdot 1}. \ed This
product structure makes $\cE$ into a filtered algebra, with the
$k$th term in the filtration\footnote{We use subscripts to denote
the filtration degree in $\cE$; while this is nonstandard, we hope
it will prevent confusion between filtration degree and the
superscripts we will use to denote the rank of a free module.}
 given by
\bd \cE_k = \left\{a_k\del^k + a_{k-1}\del^{k-1} + \dots \right\}.
\ed The induced product on $\D$ is the usual one, and $\D$ becomes a
filtered subalgebra (again, with the filtration induced from $\cE$
agreeing with the usual order filtration on differential operators).

\subsubsection{The Multi-Component KP
Hierarchy}\label{multicomp}

We now review the Lax formulation of the (multi-component) KP
hierarchy; for more information,
 see \cite{KvdL,Plaza}.

The algebra $\volt_n=\gl_n(\cE)$ of matrix-valued microdifferential
operators consists of series as in \eqref{microdiff}, but with
$a_N\in\gl_n[\![t]\!]$ with its usual multiplication. We also have a
subalgebra $\volt_n^{\C}=\gl_n(\!(\del\inv)\!)\subset\cE$ of {\em
constant coefficient} matrix microdifferential operators, and an
abelian subalgebra $\Gamma=\C(\!(\del\inv)\!)\subset\volt_n^\C$ of
constant coefficient scalar operators.

\begin{defn}
\mbox{}
\begin{enumerate}
\item A {\em matrix KP Lax operator} is a microdifferential operator
of the form \bd L=\Id
\partial+u_1\partial\inv+u_2\partial^{-2}+\cdots\in\gl_n(\cE), \ed
where  $u_i\in\gl_n[\![t]\!]$ for all $i$. We let $\cL_n$ denote the
set of all such $n\times n$ matrix KP Lax operators.
\item A {\em matrix KP wave operator} is a
microdifferential operator of the form
\begin{equation}\label{volterra op}
W=\Id+w_1\partial\inv+w_2\partial^{-2}+\cdots,
\end{equation}
where $w_i\in\gl_n[\![t]\!]$ for all $i$.  We refer to the group (in
fact group scheme) of all matrix KP wave operators as the {\em
$n$-component Volterra group} and denote it by $\Volt$.
\end{enumerate}
\end{defn}

We let $\Volt^{\C}\subset \Volt$ denote the multiplicative group of
constant coefficient matrix microdifferential operators, that is,
expressions of the form \ref{volterra op} with $w_i\in\gl_n(\C).$ We
also let $\Gap$ denote the abelian Lie algebra $\C[\del]\subset
\gl_n(\cE)$.

\begin{lemma}\label{Lax and Volterra} The set $\cL_n$ of matrix KP Lax
operators forms an infinite-dimensional affine space. The
$n$-component Volterra group scheme
 $\Volt$ acts transitively (on the left) on $\cL_n$ by $(W,L)\mapsto WLW^{-1}$.
Under this action, the stabilizer of $\del$ is the sub-group-scheme
$\Volt^{\C}$, and thus $\cL_n$ is naturally identified with
$\Volt^{\C}\backslash\Volt$ as a scheme.
\end{lemma}

\begin{defn}
The {\em multicomponent KP hierarchy} is the collection of
compatible evolution equations on a Lax operator $L$ defined as
follows:
\begin{equation}\label{KP flows}
\frac{\partial L}{\partial t_n}=[L,(L^n)_+],
\end{equation} where
$(M)_+=\sum_{N\geq 0} a_N \del^N\in\gl_n(\D)\subset\gl_n(\cE)$
denotes the differential part of a matrix microdifferential operator
$M$ as in \eqref{microdiff}.  That is, we let the operator
$L=L(t,t_1,t_2,\dots)$ depend on the infinitely many time variables
$t_n$ and then require that the dependence of $L$ on $t_n$ (i.e. its
``evolution along the $n$th time'') satisfies \eqref{KP flows}.
\end{defn}

The KP hierarchy can be written as an action of the abelian Lie
 algebra $\Gap=\C[\del]$ on $\cL_n$, i.e. a collection of commuting
 vector fields $\del_n$ on the affine space $\cL_n$ corresponding to the
 action of $\del^n\in \Gap$.  We define the vector field $\del_n$ on
 the affine space $\cL_n$ by taking its value at $L$ to be the
 commutator $[L,(L_n)_+]$.  A solution $L$ of the equations \eqref{KP
 flows} of the KP hierarchy is then just an operator
 $L(t,t_1,t_2,\dots)$ that gives (formal) integral curves of all these
 vector fields simultaneously.  Note that the first KP time $t_1$ is
 naturally identified with translation along the original variable
 $t$. In the case $n=1$, the compatibility of the
 equations \eqref{KP flows} in $x=t_2$ and $y=t_3$ (i.e. the fact that
 the corresponding vector fields on the space of Lax operators
 commute) implies that $u=u_1$ satisfies Equation \eqref{KP eq}.

\subsubsection{Sato Grassmannian}
Sato's formulation of the KP hierarchy begins with the introduction
of an infinite-dimensional Grassmannian.

Consider the vector space $\C(\!(z\inv)\!)^n$ of $n$-component
Laurent series in a parameter $z\inv$. This is a topological vector
space when equipped with a basis of open neighborhoods of $0$ given
by the subspaces $z^k\C[\![z\inv]\!]^n$ for $k\in\Z$.\footnote{More
precisely, this makes $\C(\!(z\inv)\!)^n$ into a {\em Tate vector
space}, or locally linearly compact space \cite{Hecke,Dr}. This is a
topological vector space which can be written as a direct sum of a
discrete vector space and a linearly compact vector space (the
topological dual to a discrete vector space).}  The space
$\C(\!(z\inv)\!)^n$  is the direct sum of the linearly compact
vector space $z\inv\C[\![z\inv]\!]^n$ and the discrete vector space
$\C[z]^n$.

A $c$-{\em lattice} (or simply {\em lattice}) in $\C(\!(z\inv)\!)^n$
 is a compact open subspace; one can check that, equivalently,
a vector subspace $B\subset \boldc(\!(z\inv)\!)^n$ is a c-lattice if
there exist integers $k$ and $\ell$ such that $\big(\boldc
[\![z\inv]\!] z^k\big)^n \subseteq B \subseteq \big(\boldc
[\![z\inv]\!]z^\ell\big)^n$. A $d$-{\em lattice} in
$\C(\!(z\inv)\!)^n$ is a discrete subspace that
 is complementary to a $c$-lattice.
\begin{defn}
The {\em Sato Grassmannian} $\GR_n=\GR(\C(\!(z\inv)\!)^n)$ is the
set of $d$-lattices in $\C(\!(z\inv)\!)^n$.
\end{defn}
\noindent Equivalently, one has the following well-known
description.
\begin{lemma}
The Sato Grassmannian parametrizes subspaces $B\subset
\boldc(\!(z\inv)\!)^n$ whose projections on
$\boldc(\!(z\inv)\!)^n/\C[\![z\inv]\!]^n$
 have finite dimensional
kernel and cokernel.
\end{lemma}
\noindent The index of the projection map in the lemma, also known
as the {\em index of the subspace $B$}, gives a numerical invariant
of subspaces $B$. The {\em big cell} $\GR_n^{\circ}\subset\GR_n$
consists of subspaces $B$ which project {\em isomorphically} onto
$\boldc(\!(z\inv)\!)^n/z\inv\C[\![z\inv]\!]^n\cong \C[z]^n$; such
$d$-lattices are said to be {\em generic}.

The Sato Grassmannian can be given the structure of an
infinite-dimensional scheme: see Section \ref{d-lattice section}.
With this structure, each connected component consists of exactly
the subspaces of index $k$ for a fixed $k$; we call the index $k$
component $\Gr_n^k$.

\subsubsection{KP Flows Via the Sato Grassmannian}
Sato's purpose in introducing the Grassmannian $\GR_n$ (see
\cite{Sa}) was as follows.

Given a matrix KP wave operator $W$, one obtains a free right
$\D$-submodule $W\cdot \D^n \subset \cE^n$; since $W\cdot \cE^n_{-1}
= \cE^n_{-1}$, the submodule $W\cdot\D^n$ is a generic $\D$-lattice,
and thus $W\cdot(\D^n/\D^nt)$ is a generic $d$-lattice.

The identification $\C(\!(z\inv)\!)^n=\cE^n/\cE^n t$ gives rise to
an action of the group $GL_n(\cE)$ (on the left) and of
 the Lie algebra
$\gl_n(\cE)$ on $\C(\!(z\inv)\!)^n$ and therefore, one may check,
 on the Sato Grassmannian $\GR_n$. In particular, we obtain actions on $\GR_n$ of the Volterra group
$\Volt$ and its subgroup $\Volt^{\C}$ of constant coefficient
negative microdifferential operators, as well as an action by vector
fields of the Lie algebra $\Gamma_n$.  The abelian subalgebra
$\Gamma$ thereby gives rise to an infinite family of commuting
vector fields on $\GR_n$. Sato's approach to the KP hierarchy is
then encapsulated in the following theorem (see \cite{Mu}):

\begin{thm}[Sato]
\mbox{}
\begin{enumerate}
\item The group $\Volt$ preserves the big cell $\GRon$, on which it acts simply transitively.
Thus, every generic $d$-lattice is of the form  $W\cdot(\D^n/\D^nt)$
 for a unique $W\in \Volt$.

\item The isomorphism \bd \Volt^{\C}\backslash\GRon \xleftarrow{\sim}
\Volt^{\C}\backslash\Volt \xrightarrow{\sim} \cL_n, \hskip.3in W
\cdot (\D^n/\D^nt)\longleftrightarrow W(\Id\del) W\inv \ed between
the quotient of the big cell and the space of Lax operators
identifies the infinitesimal action of
$\del^n\in\C[\del]\subset\Gamma$ on $\Volt^{\C}\backslash\GRon$ with
the $n$th KP flow $\dfrac{\partial}{\partial t_n}$ on Lax operators.
\end{enumerate}
\end{thm}
\noindent

\section{$\D$-Bundles}\label{D bundles}
In this section we study $\D$-bundles and use them to relate the KP
and CM integrable systems. In Section \ref{D bund} we define
$\D$-bundles and describe their main properties. Rank one
$\D$-bundles are, roughly, nothing more than right ideals in the
algebra of differential operators on a curve (with the caveat that
they can only be embedded in this algebra locally on the curve).
 In Section \ref{D bund and particles} we describe the
identification between moduli spaces of $\D$-bundles and moduli
spaces of CM spectral sheaves. In Section \ref{adelic Gr} we relate
the moduli spaces of $\D$-bundles to Wilson's ad\`elic Grassmannian.
In Section \ref{Sato} we give a $\D$-bundle description of the Sato
Grassmannian, which we use in Section \ref{microopers} to interpret
KP Lax operators as enhanced $\D$-bundles, the {\em micro-opers}.
Finally in Section \ref{compatibility} we explain the compatibility
between the KP hierarchy on $\D$-bundles on cubic curves and the
corresponding CM integrable systems.

\subsection{$\D$-Bundles on Curves}\label{D bund}

Let $X$ denote either a smooth quasiprojective curve or a
Weierstrass cubic curve. We denote by $\cT_X$ the tangent sheaf,
$T^*X$ the cotangent bundle, $\D=\D_X$ the sheaf of differential
operators on $X$, with the convention that in the singular cubic
case these notations are defined as in Section \ref{cubics} (i.e. as
``log" versions).

\begin{defn}[See \cite{chiral}]
A {\em $\D$-bundle} $M$ on $X$ is a locally projective
 coherent right $\D_X$-module---if $X$ is singular, we require in addition
that $M$ be isomorphic to $\D_X^n$ (for some $n$) in a neighborhood
of the singular locus.
\end{defn}
On an affine curve, one obtains a large number of examples of
$\D$-bundles of rank one by taking finitely generated right ideals
of $\D(X)$. A typical rank $1$ $\D$-bundle $M$ on a curve $X$ is not
locally free, but only {\em generically} locally free:
 away from finitely many points of $X$, $M$ is locally isomorphic to
$\D$---note the similarity to the behavior of the ideal sheaf of a
collection of points on an algebraic surface.
   For example, taking $X={\mathbf A}^1$ and thus
$\D(X) = \C[x,\del]$, the right ideal $I$ in $\D({\mathbf A}^1)$
 generated by $x^2$ and $1-x\del$
is projective (indeed, the ring $\D({\mathbf A}^1)$ is hereditary),
but there is no $f\in\C[x]$ such that $M_f$ is a free
$\D(X)$-module.
 See \cite{CH cusps, BW ideals, cusps, solitons} for more discussion.

Fix a coherent torsion sheaf $\Vc$ on $E$ supported on the smooth
locus of $E$; under the Fourier-Mukai transform (see Section
\ref{Fourier}), this determines a vector bundle $V$ on $E$.
\begin{defn}\label{V-framed D-bun}
A {\em $V$-framed $\D$-bundle} is a $\D$-bundle $M$ equipped with a
$\D$-module filtration $\{M_k\}$ and an isomorphism \bd \phi:
\oplus_{k\geq k'}\on{gr}_k(M)\rightarrow V\otimes\on{gr}_{\geq
k'}(\D) \ed of $\on{gr}(\D)$-modules
 for
some $k'$ with the following properties:
\begin{enumerate}
\item[(i)] If $E$ is singular, there is an open neighborhood $U$ of $\infty$ such
that $M|_U$ is isomorphic to $V\otimes\Dlog|_U$ compatibly with the
isomorphism $\phi$.
\item[(ii)] The canonical filtration $\{M_k\}$ of $M$ (see Definition
\ref{canonical filt}) satisfies
\begin{displaymath}
\on{rk}(M_k) =\begin{cases} \on{rk}(V)(k+1)& \text{for $k\geq 0$},\\
0& \text{if $k<0$.}\end{cases}
\end{displaymath}
\end{enumerate}
An isomorphism of $V$-framed $\D$-bundles $(M, \{M_k\}, \phi)$ and
$(M', \{M_k'\}, \phi')$ is a $\D$-module isomorphism $\psi: M
\rightarrow M'$ such that $\psi(M_k) = M_k'$ for all $k\gg 0$ and
such that $\phi'\circ\on{gr}(\psi) = \phi$ (when $k'$ is large
enough so that both sides are defined).
\end{defn}
\begin{remark}
Conditions (i) and (ii) of the definition are not necessary for a
reasonable theory of $\D$-bundles; however, they are essential for
applications to multicomponent KP. Condition (i) simply tells us, in
light of the relationship (which we will explain shortly) between
singularities of $\D$-bundles and singularities of meromorphic KP
solutions, that the KP solution is regular ``at infinity.''
 If we replace filtered modules over $\D$ by filtered
modules over $\on{Sym}(\cT_E \oplus \theo_E)$, then Condition (ii)
becomes the condition that the vector bundle on $S = {\mathbf
P}(\cT_E\oplus\theo_E)$ corresponding to the graded module $\cR(M)$
is trivial upon restriction to a generic fiber of the projection map
$S\rightarrow E$---in particular, it is an open condition.
\end{remark}

We let $\BunDP(X,V)$ denote the moduli stack of $V$-framed
$\D$-bundles (for a precise definition see Section \ref{isomorphism
theorem}). For a projective curve $X$, $\BunDP(X,V)$ is an algebraic
stack.

\subsection{Moduli of $\D$-Bundles and Calogero-Moser
Spectral Sheaves}\label{D bund and particles} We now specialize to
the case of a cubic curve $X=E$. Suppose $V$ is a semistable vector
bundle of degree zero (with trivial pullback to the normalization,
in case $E$ is singular), and let $V^\vee$ denote the finite length
sheaf on $\G\subset E$ defined as the Fourier-Mukai transform of
$V$ (see Section \ref{Fourier}). We will prove the following theorem
(Theorem \ref{KP/CM stack equiv}):

\begin{thm}\label{state KP/CM equivalence}
There is an isomorphism of stacks \bd F: \BunDP(E,V)\longrightarrow
\SSh(E,\Vc), \ed given by the $\D$-module Fourier-Mukai transform
(Theorem \ref{cubic Laumon-Rothstein}). In particular there is an
isomorphism between trivially framed rank $k$ $\D$-bundles
$\BunDP(E,\Oo^k)$ and the union over $n$ of the $k$-spin
$n$-particle Calogero-Moser spaces $\CM_n^k(E)$.
\end{thm}
\noindent Note that the following invariant of $\D$-bundles appears
implicitly in the theorem. Let $M$ be a $V$-framed $\D$-bundle on
$E$; the framing induces a canonical inclusion
$\on{gr}(M)\hookrightarrow V\otimes\on{gr}(\D)$ of finite colength
(where we have used the canonical filtration on $M$, see Definition
\ref{canonical filt}).
\begin{defn}\label{local c2}
  The {\em local second Chern
class} $c_2(M)$ is the numerical invariant \bd c_2(M) =
\on{length}(V\otimes\on{gr}(\D)/\on{gr}(M)). \ed
\end{defn}
This local second Chern class is the noncommutative analog of the
numerical invariant $\on{length}(\cF^{**}/\cF)$ for a torsion-free
coherent sheaf $\cF$ on a nonsingular surface, which, if the surface
is projective, is exactly $c_2(\cF)-c_2(\cF^{**})$: thus, it is the
part of the second Chern class that measures ``how far $\cF$ is from
being a vector bundle.'' With this definition, we then have that the
component of $\BunDP(E,V)$ parametrizing $\D$-bundles with $c_2=n$
corresponds to $\CM_n(E,\Vc)$.

Let $M$ denote a $V$-framed $\D$-bundle and $\cF$ the corresponding
$V^\vee$-framed spectral sheaf. Let $W=\pi_*\cF(-E_\infty)$ denote
the associated semistable degree $0$ vector bundle on $E$ (part (i)
of Definition \ref{CM spectral sheaves def}), so that its Fourier
transform $W^\vee$ is the torsion sheaf on $E$ measuring the
position of the Calogero-Moser particles.
\begin{prop}\label{precise generic triv} The $V$-framed $\D$-bundle $M$
is canonically identified with $V\otimes\D$ away from the location
of the corresponding CM particles. More precisely, there is a short
exact sequence
$$
0\to M\to \wt{V}\ot\D\to Q \to 0,$$ where the vector bundle $\wt{V}$
is identified with $V$ away from the support of the torsion sheaf
$W^\vee$ and $\on{supp}(W^\vee)=\on{supp}(Q)$.
\end{prop}
\begin{proof}
Give $M$ the canonical filtration.  There is then a canonical
injective
 homomorphism
$M_0=\on{gr}_0(M)\rightarrow \gr_0{V\otimes \D}=V$.
 Let $D=\on{supp}(V/M_0)$ and $U=X\setminus D$;
it is immediate that $U$ is the largest open set over which $M\cong
V\otimes \D$ and that the framing of $M$ determines an inclusion
$M\hookrightarrow (V\otimes\D)|_U$.
  Since $M$ is finitely generated,
one may choose a coherent sheaf $\wt{V}$ so that
  $V\subseteq \wt{V}\subset V|_U$ and
$M\subseteq \wt{V}\otimes\D$.

 It remains to check that
$D=\on{supp}(W^\vee)$. For this, we use the exact sequence \eqref{eq
for M} of Proposition \ref{WIT properties}; it is then immediate
from \eqref{eq for F} and the discussion following that, in the
notation used there,  $D= \on{supp}(F^0(Q_k)) = \on{supp}(W^\vee)$.
\end{proof}

\begin{remark}\label{generic trivializations}
Proposition \ref{precise generic triv} is an analog of an easy
statement about vector bundles (or torsion-free sheaves) on a ruled
surface. Namely, the condition that the vector bundle be trivial on
a fiber is open in the base (since the trivial bundle is open in
bundles on $\pline$), and over that open set, the trivialization is
fixed by choosing a trivialization of the bundle along a section.

It is not hard to see that a $V$-framed $\D$-bundle on any
nonsingular curve is canonically identified with $V\ot\D$ outside of
a finite subscheme, the {\em cusps} of the $\D$-bundle
\cite{cusps}, where the local $c_2$ is supported; moreover,
this is true for families of $\D$-bundles when properly formulated.
 This follows for
example from the local description of $\D$-bundles (Proposition
\ref{genericity criterion}): the canonical trivialization is defined
wherever the local data of the $\D$-bundle is in the big cell of
the Sato Grassmannian. In fact, there is a
more precise version of the proposition that holds on any
nonsingular curve: namely, if we choose $\wt{V}$ minimal, then
 $Q$ will have length $n=c_2(M)$.  This exact sequence will prove its value
 in
\cite{BGN}.
\end{remark}

\subsection{The Ad\`elic Grassmannian}\label{adelic Gr}

In this section we let $X$ denote an arbitrary smooth projective
curve or a Weierstrass cubic.

\begin{defn} The rank $k$ ad\`elic Grassmannian $\Gr^{ad}_k(X)$ is the set of isomorphism
classes of (unframed) $\D$-bundles $M$ equipped with a generic
trivialization $M\ot K(X)\cong \D_X^k\otimes K(X)$. When $k=1$ we
denote $\Gr^{ad}_1(X)$ by $\Gr^{ad}(X)$.
\end{defn}

\begin{remark}
The ad\`elic Grassmannian is not, in any reasonable way, a variety
or, more generally, ind-scheme.  However,
 it does have a reasonable algebro-geometric structure if
we keep track of the poles of the trivialization. More precisely,
for a finite subset ${\mathbf x}=(x_1,\dots,x_r)\in X^r$ the
ad\`elic Grassmannian $\Gr^{ad}_k(X,{\mathbf x})$ at $\mathbf x$ is
the moduli of $\D$-bundles $M$ equipped with an isomorphism
$M|_{X\sm {\mathbf x}}\to \D_X^k|_{X\sm {\mathbf x}}$. The set
$\Gr^{ad}_k(X,{\mathbf x})$ is the set of points of an ind-scheme of
ind-finite type over $X^r$, and the full ad\`elic Grassmannian is
the inductive limit of these ind-schemes under the (unfiltered)
directed system of the $X^r$ under the action of permutations and
diagonal maps.  The directed system of ind-schemes itself is what is
known as a {\em factorization space} \cite{chiral}---see \cite{W}.
\end{remark}

\begin{remark}
The ad\'elic Grassmannian we have defined
 is not quite the same as the ad\'elic Grassmannian of
Wilson \cite{Wilson CM}.  However, the difference is (mostly)
harmless: Wilson considered the case in which $X={\mathbf A}^1$ (or,
more precisely, $X={\mathbf P}^1$ with the canonical marked point
$\infty$). The triviality of $\on{Pic}(X)$ in this case (up to a
twist at the marked point $\infty$) shows that Wilson's ad\'elic
Grassmannian is the quotient of our $\Gr^{ad}(\pline)_{\infty}$ by
the action of multiplication by $K(\pline)^{\times}=\C(z)^\times$ on
the framing, together with modification at $\infty$ ($\Oo\mapsto
\Oo(k\infty)$). We prefer this definition since it is representable
by a factorization ind-scheme.
\end{remark}

Fixing $\infty\in X$, let $\BunDP(X,V)_\infty\subset \BunDP(X,V)$
denote the open subspace of $\D$-bundles which are locally free at
$\infty$ and let $\Gr^{ad}_k(X)_\infty\subset \Gr^{ad}_k(X)$ the
subset of $\D$-bundles with a generic trivialization defined at
$\infty$ (which has a structure of ${\mathbb C}$-space). For example
if $X$ is a singular cubic, we may take for $\infty$ the singular
point, where all our $\D$-bundles are already required to be
locally free.

 Let $\BunDP(X,k)_\infty\rightarrow \on{Bun}_k(X)$
denote the moduli stack for pairs $(V,M)$ consisting of a rank $k$
vector bundle $V$ on $X$ (of any degree) trivialized at $\infty$ and
a $V$-framed $\D$-bundle $M$ on $X$. Let $\K(X)_\infty$ denote
rational functions on $X$ regular at $\infty$. We state now an
algebraic description of Wilson's decomposition of the (rank one
rational) ad\`elic Grassmannian into Calogero-Moser spaces (which
we interpret as moduli spaces of $\D$-bundles):

\begin{corollary}
\mbox{}
\begin{enumerate}
\item  There is a
 morphism of spaces $\BunDP(X,\Oo^k)_\infty\to
\Gr^{ad}_k(X)_\infty$ sending a framed $\D$-bundle to its canonical
trivialization on an open $\infty\in U\subset X$.

\item  Let $GL_k(\K(X)_\infty)$ act on $\Gr^{ad}_k(X)$ by changing the
trivialization.   There is a bijection on the level of field-valued
points between $\BunDP(X,k)_\infty$ and
$\Gr^{ad}_k(X)_\infty/GL_k(\K(X)_\infty)$.
\end{enumerate}
\end{corollary}

\begin{proof}
The first assertion is a consequence of Proposition \ref{precise
generic triv} in light of Remark \ref{generic trivializations}.

Given a vector bundle $V$ and a $V$-framed $\D$-bundle $M$ as in the definition
of $\BunDP(X,k)_\infty$, choose a trivialization of $V$ in a neighborhood of
infinity compatible with the given trivialization at $\infty$.  By
Proposition \ref{precise generic triv} and
Remark \ref{generic trivializations}, we obtain an object of
$\Gr^{ad}_k(X)_\infty$.  In fact, this is easily seen to give a morphism
from a $GL_k(\K(X)_\infty)$-torsor over $\BunDP(X,k)_\infty$ to
$\Gr^{ad}_k(X)_\infty$.

To get an inverse on the level of field-valued points, we proceed as
follows.  Let $\xi=\on{Spec}(K)$.
Given an object $M \cong \D_U^k$ of $\Gr^{ad}_k(X)_\infty$ for some nonempty
open subset $U$ of $X$ that contains $\infty$, we may give $M$ the filtration
induced from that of $\D_U^k$---this makes $M$ into a $V$-framed $\D$-bundle
where $V\otimes T^\ell = \on{gr}_\ell(M)$ for $\ell\gg 0$, and moreover
$V$ is trivialized near $\infty$ by the canonical trivialization of
$\on{gr}_0(\D^k)$ there.  We thus get a function
$\Gr^{ad}_k(X)_\infty(\xi)\rightarrow \BunDP(X,k)_\infty(\xi)$.  This function
is certainly invariant under the action of $GL_k(\K(X)_\infty)$, so induces
a function
\bd
\Gr^{ad}_k(X)_\infty/GL_k(\K(X)_\infty)(\xi)\xrightarrow{p} \BunDP(X,k)_\infty(\xi).
\ed
These functions are easily checked to be inverses to each other.
\end{proof}

The ad\`elic Grassmannian was originally studied by Wilson as a
parameter space for certain algebro-geometric (finite gap)
solutions of the KP hierarchy. Namely, $\Gr^{ad}(X)$ parametrizes
all Krichever data for the curve $X$ and its {\em cusp quotients}
$X\to Y$, curves whose normalization is identified with $X$ and for
which the normalization map is a bijection. To see this we use the
identification of $\D$-bundles with torsion-free sheaves on cusp
quotients discovered by Cannings and Holland \cite{CH cusps} (this
identification is given by identifying the structure of the deRham
cohomology of $\D$-bundles with generic trivialization, as
explained in \cite{chiral}). The following proposition is due to
\cite{CH cusps} in rank one; a more general result, implying in
particular the general case, appears in \cite{cusps}.

\begin{prop}\cite{CH cusps,cusps}
$\Gr^{ad}_k(X)_\infty$ is isomorphic to the direct limit over cusp
quotients $X\to Y$ smooth at $\infty$ of the set of isomorphism
classes of rank $k$ torsion-free $\OY$-modules equipped with a
generic trivialization defined at $\infty$.
\end{prop}
\noindent
Note that the same statement holds, if we remove the restrictions
on $\infty$.

This interpretation of $\D$-bundles gives rise to a construction of
solutions to the multicomponent KP hierarchy (orbits of the action
of $\C[z]$ on $\GR^k/GL_k(\C[[z\inv]])$) from $\D$-bundles on $X$,
smooth at $\infty$, considered as defining Krichever data on cusp
quotients $Y$ of $X$:

\begin{corollary} Fix a local coordinate $z\inv$ on $X$ at $\infty$. Then there is a
canonical Krichever map $\Gr^{ad}_k(X)_\infty\to \GR^k$ from the
ad\`elic Grassmannian to the Sato Grassmannian, sending a rank $k$
torsion-free sheaf $\cV$ on a curve $Y$ with trivialization near
$\infty$ to the subspace of $\C(\!(z\inv)\!)^k$ defined by sections
of $\cV$ on $Y\sm \infty$. Quotienting by change of trivialization
gives a composite
$$\BunDP(X,k)_\infty\to \Gr^{ad}_k(X)_\infty/GL_k(\K(X)_\infty) \to
\GR^k/GL_k(\C[[z\inv]]).$$
\end{corollary}

It is also easy to describe the corresponding flows directly in
terms of $\D$-bundles: the flows modify a $\D$-bundle $M$ inside
its localization $M\ot\C(\!(z\inv)\!)$ by multiplication by powers
of $z$. Below we will introduce a ``Fourier dual" relation between
$\D$-bundles on cubic curves and KP solutions, where we modify $M$
inside its microlocalization by powers of a vector field $\del$.

\subsection{Local $\D$-Bundles and the Sato Grassmannian}\label{Sato}
In this section we study $\D$-bundles on the disc $D=\spec
\C[\![t]\!]$: so, we let $\D$ and $\cE$ denote the rings of
differential and microdifferential operators with coefficients in
${\mathbf C}[\![t]\!]$, as in Section \ref{KP section}.
\footnote{It is important to note that we use the disc
rather than the {\em formal} disc $\wh{D}=\on{Spf}\C[\![t]\!]$, on
which the theory of $\D$-modules is trivial.}
 We introduce $\D$-lattices and state our
generalization of Sato's $\D$-module description of the Sato
Grassmannian.

\subsubsection{$\D$-Lattices}

Consider the fiber $\cE|_{t=0}\overset{\on{def}}{=}\cE/\cE\cdot t$
of the right
 $\C[\![t]\!]$-module $\cE$.
Under the identification $z\longleftrightarrow\del$, we see that
this fiber is isomorphic to the vector space $\C(\!(z\inv)\!)$; its
realization as $\cE/\cE t$, however, gives it a natural structure of
left $\cE$-module. Likewise, the fiber $\cE^n|_{t=0}$ is identified
with the vector space $\C(\!(z\inv)\!)^n$. Taking the fibers at zero
of the right $\C[\![t]\!]$-submodules $\D^n$ and $\cE_{-1}^n$ of
$\cE^n$, we obtain the subspaces $\C[z]^n$ and
$z\inv\C[\![z\inv]\!]^n$, respectively, of $\C(\!(z\inv)\!)^n$.

Generalizing the example of $\D^n\subset\cE^n$, we have the
following.

\begin{defn}\label{pointwise D-lattice}
A {\em $\D$-lattice} in $\cE^n$ is a finitely generated right
$\D$-submodule $M\subset \cE^n$ such that $M\cdot\cE=\cE^n$
(equivalently, the natural map $M\otimes_\D \cE \rightarrow \cE^n$
is an isomorphism).

A $\D$-lattice is said to be {\em generic} if it is transversal to
$\cE_{-1}^n$, that is, $\cE^n = M\oplus \cE^n_{-1}$.
\end{defn}

A $\D$-lattice $M\subset \cE^n$ is a $\D$-bundle on the disc with
a microlocal trivialization (i.e. near $\del=\infty$). Namely, $M$
inherits the structure of $\Oo^n$-framed $\D$-bundle from its
embedding in $\cE^n$. Thus we may
reformulate the notion of $\D$-lattice as an $\Oo^n$-framed
$\D$-bundle on the disc, equipped with a {\em filtered} isomorphism
$M\ot_\D\cE\to \cE^n$.

\begin{lemma}
Let $M\subset\cE^n$ be a finitely generated $\D$-submodule. Then $M$
is a $\D$-lattice if and only if the following conditions are
satisfied:
\begin{enumerate}
\item There exists an integer $k$ such that $M\rightarrow \cE^n/\cE^n_k$ is
surjective.
\item There exists an integer $\ell$ such that $M\rightarrow \cE^n/\cE^n_\ell$
is injective with cokernel a finitely generated projective
$\C[\![z]\!]$-module.
\end{enumerate}
\end{lemma}
\noindent The proof uses standard techniques.
In fact we will later adopt this equivalent reformulation, since it
is easier to work with in families than Definition \ref{pointwise
D-lattice}: see Section \ref{D and d}.

We also have the following characterization of generic
$\D$-lattices:
\begin{prop}\label{genericity criterion}
A $\D$-lattice $M\subset \cE^n$ is generic if and only if $M$ (with
its induced filtration, Section \ref{microloc section}) is
isomorphic to $\D^n$ (with its standard filtration) as filtered
$\D$-modules.
\end{prop}
\begin{proof}
$M$ is generic if and only if, with the induced filtration, the
induced
 homomorphism
\bd \psi: \on{gr}(M)\rightarrow \on{gr}(\cE^n/\cE_{-1}^n)\cong
\on{gr}(\D^n) \ed of $\C[\![t]\!]$-modules is an isomorphism.
Observe that
 $\psi$ has a kernel if and
only if $M_{-1}\neq 0$.  Furthermore, if $\on{ker}(\psi) = 0$, then
$\psi$ is an isomorphism if and only if $\on{gr}(M) \cong
M_0\otimes\on{gr}(\D)$.  A standard argument shows that this last
equality holds if and only if the induced $\D$-module map
$M_0\otimes\D \rightarrow M$ is an isomorphism of filtered
$\D$-modules.
\end{proof}

\subsubsection{Sato Grassmannian as Moduli Space}
We now state our general Sato theorem (restated and proved below as
Theorem \ref{Sato-type thm}), relating $\D$-lattices to the
Grassmannian.

To any $\D$-lattice $M\subset\cE^n$ we may assign its fiber $M/Mt$
at $0$. This is a subspace (with no additional module structure in
general) of $\cE^n/\cE^nt = \C(\!(z\inv)\!)^n$, and in fact always
gives a $d$-lattice: see Proposition \ref{D nat trans}.

We then have the following theorem:
\begin{thm}[See Theorem \ref{Sato-type thm}]\label{state Sato-type}
The Sato Grassmannian $\GR_n$ is isomorphic to the moduli space for
$\D$-lattices, under the map taking a $\D$-lattice $M\subset \cE^n$
to the $d$-lattice $M/Mt \subset \cE^n/\cE^n=\C(\!(z\inv)\!)^n$.
Under this map, generic $\D$-lattices are identified with
$d$-lattices in the big cell.
\end{thm}
The set-theoretic identification of the set of {\em generic} (i.e.
free) $\D$-lattices and the set of {\em generic} $d$-lattices (the
big cell $\GR_n^\circ$) was proven by Sato; see \cite{Mu}.  We will
postpone the proof until Section \ref{D and d}.

\subsubsection{Interpretation: Noncommutative Krichever Data}
Recall that a (rank one) Krichever datum (defining an
algebro-geometric solution for KP \cite{Mu}) associated to a curve
$X$ with a marked smooth point $x$ and local coordinate $z$ at $x$
is a torsion-free sheaf $\cL$ on $X$, equipped with a
trivialization near $x$. Equivalently we have a torsion-free sheaf
$\cL$ on the affine curve $X\sm x$ equipped with an isomorphism
$\cL\ot_{\Oo_{X\sm x}}\C(\!(z)\!)$ (allowing us to glue $\cL$ on the
punctured curve to the trivial bundle on the disc at $x$). Theorem
\ref{state Sato-type} claims that the (rank one) Sato Grassmannian
is identified with the moduli space of (rank one) $\D$-bundles on
the disc, equipped with a trivialization near infinity (i.e. of the
associated filtered $\cE$-module). Thus the entire Grassmannian
parametrizes a kind of Krichever data for the noncommutative
$\pline$-bundle $\boldp(\D)$. The latter is characterized as the
space on which coherent sheaves are filtered $\D$-modules (if we
adopt the definition of noncommutative variety as Grothendieck
category or as differential graded category), and is obtained by
gluing the noncommutative affine bundle $\spec \D$ to a ``disc
bundle at infinity" $\spec \cE_-$ along $\spec\cE$. Classical
Krichever solutions correspond to commutative subrings of $\D$
(necessarily the affine rings of curves), in other words to maps
from this noncommutative curve to ordinary curves. Thus classical
algebro-geometric solutions correspond to $\D$-lattices which are
pulled back from torsion-free sheaves $\cL$ on curves $X$ under
maps $\boldp(\cD)\to X$.

\subsection{Lax Operators and Micro-Opers}\label{microopers}

As we have seen, the
$\D$-lattices provide a $\D$-module interpretation for points
in the Sato Grassmannian, the big cell of which parametrizes
 KP wave operators. In this section we
introduce {\em micro-opers}, which give a similar interpretation for
points in the quotient $\Volt^{\C}\backslash\GR_n$ and thus also for
 KP Lax operators. They are
less rigid structures than $\D$-lattices and hence turn out to be
better suited for our goal of ``globalizing'' the Sato dictionary and
geometrically
interpreting meromorphic Lax operators on a (cubic) curve.

Given a right $\D$-module $M$, we write $M_\cE = M\otimes_\D \cE$.

\begin{defn}
A $V$-framed {\em micro-oper} on a differential curve (Section
\ref{cubics}) $(X,\del)$ is a $V$-framed $\D$-bundle $M$ on $X$
equipped with a left $\cE$-module endomorphism $\del_M$ of
$M_\cE=M\ot_\D\cE$ whose principal symbol is $\Id_V\otimes \del$
with respect to the induced filtration of $M_\cE$. In other words,
$\del_M$ has degree $1$ with respect to the filtration and induces
the isomorphism
$$\left(\on{gr}(\del_M)=\Id_V\otimes\del\right):\left(\on{gr}_n M_\cE\simeq V\ot \gr_n
\cE\right) \longrightarrow \left(\on{gr}_{n+1} M_\cE\simeq
V\ot\gr_{n+1}\cE\right).$$
\end{defn}

By a {\em rank $n$ local micro-oper} we denote an $\Oo^n$-framed
micro-oper on the disc. A local micro-oper is {\em generic} if the
underlying filtered $\D$-module is trivial, i.e. isomorphic to
$\D^n$ with its filtration. Local micro-opers are parametrized by a
quotient of the Sato Grassmannian:

\begin{prop}\label{microopers and GR}
The parameter space (moduli stack) for rank $n$ local micro-opers
is the quotient $\Volt^{\C}\backslash\GR_n$ of the Sato Grassmannian
by negative constant coefficient operators.
\end{prop}

\begin{proof}
A $\D$-lattice $M\hookrightarrow \cE^n$ defines a micro-oper
structure on $M$ by remembering only the natural filtration, framing
and $\del$ action on $\cM_E=\cE^n$. This defines a map from $\GR^n$
to micro-opers. Conversely, given a micro-oper structure on $M$ we
may pick a $\D$-lattice structure (i.e. an isomorphism
$M\ot_\cD\cE\to \cE^n$) on $M$ compatible with filtration and
framing. Since $\del_M$ is acting by right $\cE$-module
endomorphisms of $M_\cE$, the lattice structure identifies $\del_M$
with an endomorphism of $\cE^n$ given by left multiplication by an
operator of the form $\del\otimes\Id+a_0+a_1\del\inv+\cdots$, where
$a_i$ are matrices over formal power series. Changing the
$\D$-lattice structure by the left action of $\on{Id} +
\gl_n(\cE_{-1})$ on $\cE^n$, we can conjugate the image of $\del_M$
to $\del\ot\Id$.  Moreover, we can do so uniquely up to the
centralizer of $\del\ot\Id$ in $\gl_n(\cE_<)$, namely $\Volt^{\C}$.
This identifies the $\D$-lattice structures on $M$ inducing the
given micro-oper structure on $M$ with a $\Volt^{\C}$-orbit on
$\GR^n$, as desired.

The above proof extends immediately to $S$-families of micro-opers
(note that $\D$-lattice structures on $M$ will exist only locally
on the parameter scheme $S$),
proving the stronger assertion of the proposition.
\end{proof}
\begin{corollary}
The space of generic local micro-opers is isomorphic to the affine
space of matrix KP Lax operators.
\end{corollary}
\begin{proof}
The genericity condition precisely characterizes micro-opers in the
image of $\Volt^{\C}\backslash\GRon$, since $\D$-lattice structures
on such a framed $\D$-module are automatically in the big cell by
Proposition \ref{genericity criterion}. Note also that $\D^n$ has no
automorphisms as a framed $\D$-module. Thus a trivial framed
$\D$-module has a unique trivialization. Therefore $\del_M$ defines
a right $\cE$-module endomorphism of $\cE^n$, hence an element $L$
of $\gl_n(\cE)$ acting from the left. The micro-oper conditions
guarantee that this is indeed a matrix KP Lax operator, i.e. has the
form $\del\otimes\Id$ plus lower order terms.
\end{proof}

On a global differential curve $(X,\del)$ we have the following
``twisted'' analog of a matrix Lax operator:
\begin{defn}
A {\em $V$-twisted Lax operator} over an open set $U\subseteq X$ is
a filtration-preserving right $\cE$-module endomorphism \bd L:
V\otimes \cE\rightarrow V\otimes\cE(1) \ed with principal symbol
$\on{Id}\otimes\del$.
\end{defn}
\noindent Note that ``filtration-preserving'' in this case means
that \bd L(V\otimes\cE_k)\subset (V\otimes\cE(1))_k =
V\otimes\cE_{k+1} \ed for all $k$.  An easy argument shows that, in
the case $V=\theo^n$, such $L$ are exactly the matrix Lax operators
(more precisely, their initial values for fixed KP times).

Let us fix a point $x\in X$ at which $\del\neq 0$
 and use $\del$ to identify $\widehat{\theo}_{X,x}$
with $\C[\![t]\!]$. Then the restriction of any micro-oper on $X$
near $x$, together with a trivialization of $V$ near $x$, define a
micro-oper on the disc, i.e. a point of $\Volt^{\C}\backslash\GR^n$.
If this local micro-oper is generic then this is equivalent to the
data of a matrix Lax operator on the disc. In fact, by Proposition
\ref{precise generic triv} there is an open set $U$ on which $M_\cE$
is canonically identified with $V\ot\cE$, so that the $\cE$-module
endomorphism $\del_M$ gives rise to a ($V$-twisted) matrix
microdifferential operator of degree one and principal symbol the
identity:

\begin{corollary}\label{rational Lax}
Let $(M,\del_M)$ denote a $V$-framed micro-oper, and $U\subset X$
the open subset where $M$ is trivial as a $V$-framed $\D$-bundle
(Proposition \ref{precise generic triv}). Then the micro-oper
$M|_U$ determines, and is determined by, a $V$-twisted Lax operator.
\end{corollary}

\subsubsection{Micro-Opers and Opers}
We see that micro-opers on a curve are global geometric analogs of
matrix Lax operators on the disc. Thus they are the matrix KP
analogs of opers, introduced by Beilinson and Drinfeld \cite{BD
opers} following Drinfeld and Sokolov \cite{DS} (henceforth we only
consider $GL_n$ opers). Opers are special connections defined on any
smooth curve, while opers on the disc are identified with KdV Lax
operators, i.e. $n$th order differential operators
$L=\del^n+u_1\del^{n-1}+\cdots+u_n\in\D$.

More precisely, micro-opers are the analogs of affine opers,
introduced in \cite{BF} for the geometric study of Drinfeld-Sokolov
hierarchies. A ($GL_n$-)affine oper on a curve $X$ is a vector
bundles $V$ on $X\times\pline$, equipped with a connection $\nabla$
along $X$ on the bundle of sections $V_{\aline}$ of $V$ on
$\aline\times X$, and a flag $V_\infty^{\cdot}$ on the fiber
$V_\infty$ of $V$ at $X\times\infty$. The connection is required to
have a first order pole at $\infty$ and to satisfy a strict form of
Griffiths transversality with respect to the flag at $\infty$. The
open set of {\em generic} affine opers, for which $V$ is trivial
along $\pline$, is identified with opers. Affine opers on the disc
are identified with (a quotient of) the loop Grassmannian for
$GL_n$, while opers form the corresponding big cell.

The $GL_n$ loop Grassmannian is embedded in the Sato Grassmannian,
reflecting the inclusion of the KdV (and Gelfand-Dickey)
hierarchies in KP. This corresponds to the identification of affine
opers with special micro-opers. Namely, affine opers (on any curve
differential $(X,\del)$) are identified with $\Oo$-framed
micro-opers $(M,\del_M)$ for which $\del_M^n$ preserves the
submodule $M\subset M_\cE$. The identification preserves big cells:
$M$ is locally free if and only if the corresponding affine oper is
generic. To define the vector bundle $V$ on $X\times \pline$ we
consider $M_\cE$ as a $\Oo_X(\!(z\inv)\!)$-module via the
endomorphism $z=\del_{M}^n$. The extension to $\pline\times X$ and
the flag at $\infty$ are constructed from the filtration, while the
affine oper connection comes from the (left version of the) right
$\D$-module structure on $M_\cE$. (See \cite{solitons} for a more
leisurely discussion.)

\subsection{Flows of the KP and CM Hierarchies}\label{compatibility
of flows} The flows of the multicomponent KP hierarchies on the Sato
Grassmannian and the space of Lax operators have natural
formulations in terms of $\D$-lattices and micro-opers,
respectively. First recall that the KP flows on $\GR_n$ are given by
the action of the subalgebra
$\Gamma_+=\C[\del]\on{Id}\subset\gl_n(\cE)$ of constant coefficient
scalar differential operators. The action of this Lie algebra has
the following simple description on $\D$-lattices. Given a
polynomial $P(\del)\in \C[\del]$, its infinitesimal action on a
$\D$-lattice $M\subset\cE^n$ is given by translating the submodule
$M$ by the action of the (right $\D$-module) endomorphism of
$\cE^n$ given by left multiplication by $P(\del)$. This means we
deform elements $m\in M\subset \cE^n$ to $m+\epsilon P(\del)\cdot
m$. On the level of tangent spaces, write $B = M/Mt$ for the
corresponding $d$-lattice, and make the identification \bd
\on{Hom}\left(B, \C(\!(z\inv)\!)^n/B\right) = \on{Hom}_{\D}(M,
\cE^n/M). \ed Then the tangent space to $\GR_n$ at $M$ has the form
\bd T_M\GR_n = \on{Hom}\left(B, \C(\!(z\inv)\!)^n/B\right). \ed
The vector field on $\GR_n$ given by $P(\del)$ then
 has value at $M$ equal to the composite \bd M \hookrightarrow \cE^n
\xrightarrow{P(\del)\cdot -} \cE^n \twoheadrightarrow \cE^n/M; \ed
this is typically nontrivial since $M$ is a right, but not left,
$\D$-submodule of $\cE^n$.

The action of $\Gamma_+$ descends from $\GR_n$ to its quotient
$\Volt^{\C}\sm \GR_n$ (whose big cell parametrizes Lax operators).
By Proposition \ref{microopers and GR} this is the parameter space
for local micro-opers. It is easy to describe this action directly
on micro-opers:

\begin{prop} The matrix KP flow associated to a polynomial $P(\del)\in \C[\del]$
on a local micro-oper $(M,\del_M:M_\cE\to M_\cE$ is given by
translating the $\D$-module $M\subset M_\cE$ by the action of the
(right $\D$-module) endomorphism $P(\del_M)$ (that is we deform
elements $m\in M\subset M_\cE$ to $m+\epsilon\; P(\del_M)\cdot m$).
\end{prop}

Micro-opers on a global differential curve $(X,\del)$ (i.e.
meromorphic twisted Lax operators) also carry an action of the
abelian Lie algebra $\C[\del]$, defined just as in the local
setting.\footnote{Note that we don't in fact need the full structure
of differential curve, only the choice of $\del$ up to scalar, which
is canonical in the cubic case.} Namely, the action of $P(\del)\in
\C[\del]$ on a micro-oper $(M,\del_M)$ is simply given by
translating the $\D$-submodule $M\subset M_\cE$ by the action of
the (right $\D$-module) endomorphism $P(\del_M)$, and preserving
the filtration, framing and endomorphism $\del_M$ on $M_\cE$. In
particular, these flows act on the associated rational Lax operator
of Corollary \ref{rational Lax} by the standard multi-component KP
flows. More precisely, if we trivialize $V|_U$ then $P$ acts on the
Lax operator $L_M\in \gl_n(\cE)$ by commutator with $P(\del)$.

\subsubsection{The KP Algebroid}
We may also describe the KP flows on micro-opers in terms of
deformations of the underlying $\D$-bundles. More precisely, we
introduce, by analogy with Section \ref{Lie algebroids}, a Lie
algebroid on $\D$-bundles, which describes all deformations of a
$\D$-bundle $M$ coming from its microlocal endomorphisms (such as
$\del_M$), acting by moving $M$ inside $M_\cE$:

\begin{lemma} The sheaf $\underline{\End}_{\cE}$ over $\BunDP(E,V)$
of endomorphisms of the microlocalization of the universal sheaf
(i.e. the sheaf whose fiber at $M$ is $\End_{\cE}(M_\cE)$) has the
structure of Lie algebroid, the {\em KP algebroid}, whose action on
a $\D$-bundle $M$ deforms $M\subset M_\cE$ by left multiplication.
\end{lemma}

On the other hand, from Definition \ref{CM flows definition} we have
an action of $\C[\del]$ on $\Vc$-framed CM spectral sheaves
$\CM_n(E,\Vc)$. More generally by Lemma \ref{CM algebroid} we have a
natural algebroid describing all deformations of a CM spectral sheaf
along the curve $E_\infty$. The definitions of the CM and KP
algebroids are precisely analogous, as are the structures of Higgsed
spectral sheaf and micro-oper. It is then an immediate consequence
of the Fourier-Mukai construction of the isomorphism between
$\D$-bundles and spectral sheaves (Theorem \ref{state KP/CM
equivalence}) (and of the compatibility of the Fourier-Mukai
transform with microlocalization \cite{PRo}) that the corresponding
algebroids and hierarchies are identified:

\begin{thm}\label{compatibility}
Let $\bF:\BunDP(E,V)\to\CM_n(E,\Vc)$ denote the Fourier-Mukai
isomorphism of framed $\D$-bundles and framed CM spectral sheaves.
\begin{enumerate}
\item $\bF$ identifies the KP and CM algebroids.

\item $\bF$ lifts to an isomorphism of the moduli stack of $V$-framed micro-opers
with the moduli stack of $\Vc$-framed Higgsed CM spectral sheaves,
identifying the multicomponent KP hierarchy on micro-opers with the
Higgsed CM hierarchy (i.e. intertwining the two
$\C[\del]$-actions).
\end{enumerate}
\end{thm}

We defer the proof to Section \ref{proof of compatibility}, after
the relevant facts about the Fourier-Mukai transform have been
established.

\subsubsection{Bispectrality}\label{bispectrality}
We have now defined two relations between $\D$-bundles and KP
solutions: the construction of rational Lax operators on $X$ from
micro-opers on $X$, and the construction of Krichever data from the
ad\`elic Grassmannian, i.e. $\D$-bundles with generic
trivialization. Alternatively, we have define two Lie algebroids on
the moduli $\BunDP(E,V)$ of $\D$-bundles, describing deformations
of a $\D$-bundle $M$ by endomorphisms of $M\ot_\D\cE$ and of
$M\ot_\D\D(\!(z\inv)\!)$.   For
$\D$-bundles on a general curve, there is no obvious relation
between the two constructions. In particular note that the
constructions of KP solutions land in different copies of the Sato
Grassmannian: Krichever data give rise to subspaces of
$\C(\!(z\inv)\!)$, where $z\inv$ is a local coordinate at a point
$\infty\in X$, while micro-opers give rise to subspaces of
$\C(\!(\del\inv)\!)$, where $\del$ is a nonvanishing vector field on
$X$.

It is natural to expect that the relation between the two
constructions is a sort of Fourier transform, exchanging microlocal
and local trivializations. This can be made precise in the rational
case $X=\pline$, using the geometric Fourier transform.\footnote{
Note that the categories of $\D$-modules on the rational cubic and
on $\pline$ are equivalent, and $\D$-bundles on $\pline$ trivial at
$\infty$ and $\D$-bundles on the cuspidal cubic curve are
canonically identified \cite{SS,cusps}. } This is an autoequivalence
of the category of $\D$-modules on $\aline$, and exchanges
$\C(\!(z\inv)\!)$ and $\C(\!(\del\inv)\!)$. Let us then consider the
moduli problem for $\D$-bundles on $\pline$, with framing (i.e.
trivialization at $\del\inv=0$) and trivialization at $z\inv=0$,
i.e. the noncommutative version of sheaves on the quadric
$\pline\times \pline$, framed along $\pline \times\infty \cup
\infty\times\pline$.  This moduli stack acquires an automorphism
from the geometric Fourier transform, interchanging the KP algebroid
(deformations at $\del\inv=0$) with the algebroid of deformations at
$z\inv=0$ (given by endomorphisms of $\D$-bundles restricted to
$\spec \C(\!(z\inv)\!)$). In other words, the automorphism exchanges
the spectral parameter $\del$ and the differential parameter $z$,
and is easily seen to give (on the level of points) the bispectral
involution \cite{Wilson bispectral,Wsurv2} in the rank one case.
Specifically this is a geometric reformulation of the picture of
Berest and Wilson \cite{BW automorphisms}, and the two algebroids in
the rank one case are the two commutative subalgebras of the
automorphisms of the Weyl algebra studied in \cite{BW
automorphisms}.

\section{Fourier Transform for Cubic Curves}\label{Fourier}
In this section we explain the Fourier-Mukai autoequivalence of the
derived category of a cubic curve, extend it to give an analog of the
theorem of Laumon and Rothstein concerning $\D$-modules on abelian
varieties, and give some fundamental calculations for the
Fourier-Mukai transform of torsion-free sheaves.

\subsection{Fourier Equivalence for Weierstrass Cubics}
Fix a cubic curve $E$ as before. In this section we describe
 the Fourier-Mukai
autoequivalence of the derived category of coherent sheaves on $E$.

 Recall that the generalized Jacobian $\Jac(E)$ (the
group of line bundles of degree $0$ on $E$) is isomorphic to the
smooth locus $\G\subset E$ via the map $e\mapsto \theo(e-b)$. The
compactified Jacobian of $E$, denoted by $\overline{\Jac}(E)$, is
the moduli space for rank 1, torsion-free sheaves on $E$ of degree
(equivalently, Euler
characteristic) $0$; it is isomorphic to $E$, with the additional
point (when $E$ is singular) coming from the unique rank 1
torsion-free sheaf of degree $0$ that is not locally free (which is
isomorphic, in a neighborhood of the singular point, with the
ideal of the singular point).

Let $\Delta\subset E\times E$ denote the diagonally embedded copy of
$E$. We define a sheaf $\cP^\vee$ on $E\times \overline{\Jac}(E) =
E\times E$ by $\cP^\vee = I_{\Delta}\otimes
(\theo(b)\boxtimes\theo(b))$. $\cP^\vee$ is flat over both factors:
indeed, we have a short exact sequence
\begin{equation}\label{Pcheck}
0 \rightarrow \cP^\vee \rightarrow \theo(b)\boxtimes\theo(b)
\rightarrow \theo_{\Delta}(2b)\rightarrow 0
\end{equation}
with the second and third terms flat over both factors, implying
that the kernel $\cP^\vee$ is as well.

\begin{defn}
Let $\pi_i:E\times E\rightarrow E$ denote projection on the $i$th
factor.  The {\em Fourier functor} $F$ from the coherent derived
category $D^b_{\on{coh}}(E)$ to itself is defined by \bd F(M) =
{\mathbf R}{\pi_2}_*(\cP^\vee\overset{{\mathbf L}}{\otimes} \pi_1^*
M). \ed
\end{defn}
\noindent Note that ${\mathbf L}\pi_i^* = \pi_i^*$ since the
projection is flat.  Note also that our curve $E$ is Gorenstein, so
the dualizing sheaf $\omega = \omega_E$ is a line bundle---in fact,
for our cubic curves, $\omega_E \simeq \theo_E$.

To the sheaf $\cP^\vee$ we may also associate its {\em derived
dual}, the object $(\cP^\vee)^\vee = {\mathbf
R}\uHom(\cP^\vee,\theo_{E\times E})$ of the derived category of
$E\times E$.  We write $\cP = (\cP^\vee)^\vee$; this notation is
consistent, since by taking the derived dual twice one obtains a
complex quasi-isomorphic to the original complex $\cP^\vee$.
 In fact, $\cP$ is
a sheaf that fits in an exact sequence
\begin{equation}\label{Pform}
0\rightarrow \theo(-b)\boxtimes\theo(-b)\rightarrow \cP\rightarrow
\theo_{\Delta}(-2b)\rightarrow 0
\end{equation}
and, like $\cP^\vee$, $\cP$ is a flat family of torsion-free
sheaves of rank $1$ and degree $0$ over both factors.
 Moreover, let $(-1): E\rightarrow E$ denote
the involution on $E$ induced by the inverse for the group structure
of $\G$.  Then
\begin{equation}\label{P formula}
\cP = (\on{id}\times (-1))^*\cP^\vee.
\end{equation}

Consider the functor
\begin{equation}\label{Fbar}
\overline{F}(N) ={\mathbf R}{\pi_1}_* \big(\cP\overset{{\mathbf
L}}{\otimes} \pi_2^* N   \big)[1].
\end{equation}

\begin{thm}[See \cite{BuK}, Theorem 2.12]\label{FM equiv}
The Fourier functor $F: D^b_{\on{coh}}(E)\rightarrow
D^b_{\on{coh}}(E)$ is an exact equivalence of triangulated
categories, with quasi-inverse given by \eqref{Fbar}.
\end{thm}
\begin{remark}
In \cite[Theorem 5.2(1)]{solitons}, we announced a proof of this theorem
(a proof using Bridgeland's criterion appeared in an early uncirculated
draft of the present paper).
The theorem cited in \cite{BuK} is actually considerably more
general, but a special case of that general theorem gives the fact we need.
\end{remark}

\begin{notation}
We let $\Fin = \overline{F}[-1]$.
\end{notation}

\begin{corollary}\label{FM equiv2}
The Fourier functor $F: D^b_{\on{qcoh}}(E)\rightarrow
D^b_{\on{qcoh}}(E)$ is an exact equivalence of triangulated
categories.
\end{corollary}
\begin{proof}
A complex $M$ of quasicoherent sheaves is the colimit of coherent
complexes $M_i$.  By Remark 2.2 and Lemma 4.1 of
\cite{Bokstedt-Neeman}, we then have \bd \overline{F}\circ F(M) =
\underset{\longrightarrow}{\lim}\, \overline{F} F(M_i) \simeq M
\;\;\text{and}\;\; F\circ \overline{F}(M) =
\underset{\longrightarrow}{\lim}\, F\overline{F}(M_i)\simeq M \ed as
desired.
\end{proof}
\begin{corollary}\label{circle prod}
We have
\begin{equation}
{\mathbf R}(p_{13})_*\left(p_{12}^*\cP\otimes
p_{23}^*\cP^\vee\right) \simeq \Delta_*\theo_E[-1]
\;\;\text{and}\;\; {\mathbf
R}(p_{13})_*\left(p_{12}^*\cP^\vee\otimes p_{23}^*\cP\right) \simeq
\Delta_*\theo_E[-1]
\end{equation}
in the derived category of coherent $\theo_{E\times E}$-modules.
\end{corollary}
\begin{proof}
This is immediate from Corollary \ref{FM equiv2} and \cite{Toen}.
\end{proof}

\subsection{Fourier Transform of Torsion-Free Sheaves}
Our goal in this subsection is to prove analogs, for singular
Weierstrass cubics, of the usual characterizations of torsion-free
sheaves on an elliptic curve whose Fourier transforms are again
sheaves (that is, have cohomology in only a single degree).  The
reader who is familiar with the standard techniques for these
problems may consult the statements of Propositions
\ref{WIT0} and \ref{WIT1} for the expected facts; note, however,
that while the proof of Proposition \ref{WIT1} is standard,
the proof of Proposition \ref{WIT0} is {\em not} entirely standard.

Recall that, if $M$ is a torsion-free coherent sheaf on $E$, the
{\em slope} of $M$ is
\begin{equation}
\mu(M) \overset{\on{def}}{=}\frac{\on{deg}(M)}{\on{rk}(M)}.
\end{equation}
On a Weierstrass cubic curve $E$, the slope of $M$ satisfies $\mu(M)
= \chi(M)/\on{rk}(M)$ where $\chi(M)$ is the Euler characteristic of
$M$; see Section 0.2 of \cite{FMW} for this and other basic facts
about torsion-free sheaves and semistability on Weierstrass cubics.

Recall the {\em Harder-Narasimhan filtration} of a torsion-free
coherent sheaf $M$ on $E$: this is the unique decreasing filtration
\begin{equation}\label{HN filtration}
F_{\infty}(M) = 0 \subset F_{\mu_1}(M) \subset F_{\mu_2}(M) \subset
\dots \subset F_{\mu_r}(M) = M
\end{equation}
(that is, $\mu_1 > \mu_2 > \dots > \mu_r$) such that each term
$\on{gr}_{\mu_i}(M) = F_{\mu_i}(M)/F_{\mu_{i-1}}(M)$ of the
associated graded sheaf is torsion-free and semistable of slope
$\mu_i$.

\begin{prop}\label{WIT0}
Suppose $M$ is a torsion-free coherent sheaf on $E$. In the notation
of \eqref{HN filtration}, $F(M)$ and $\overline{F}(M)$ are
concentrated in cohomological degree $0$ if and only if $\mu_r >0$.
\end{prop}
\begin{proof}
We give the proof only for $\overline{F}(M)$; the proof for $F(M)$
is nearly identical in light of Equation \eqref{P formula}.

Use the exact sequence \bd 0\rightarrow \cP^{\vee} =
I_{\Delta}\otimes (\theo(b)\boxtimes\theo(b)) \rightarrow
\theo(b)\boxtimes\theo(b) \rightarrow \theo_{\Delta}(2b) \rightarrow
0 \ed on $E\times E$.  Tensoring with $p_1^*M$ and applying
${\mathbf R}{p_2}_*$ gives an exact sequence \bd 0\rightarrow
\overline{F}^0(M)\rightarrow \theo(b)\otimes H^0(M(b)) \rightarrow
M(2b) \rightarrow \overline{F}^1(M) \rightarrow \theo(b)\otimes
H^1(M(b)) \rightarrow 0 \ed on $E$.  So, $\overline{F}(M)$ is a
sheaf in degree $0$ if and only if
\begin{enumerate}
\item[(a)] $H^1(M(b))=0$ and
\item[(b)] $M(b)$ is globally generated.
\end{enumerate}
First, then, suppose that $\mu_r >0$; we will prove that $M$
satisfies (a) and (b).  An inductive argument using the long exact
cohomology sequence shows that, if $H^1(T(b)) =0$ for all
torsion-free semistable $T$ of slope $\mu(T)>0$ then also $H^1(M(b))
= 0$. But for such $T$, Serre duality gives $H^1(T(b))^* =
\Hom(T(b),\theo) = \Hom(T,\theo(-b)) = 0$ by the semistability of
$T$.  This proves (a).

For (b), we suppose for the moment that
\begin{equation}\label{gg}
\text{$T(b)$ is globally generated for all semistable $T$ with
$\mu(T)>0$}
\end{equation}
and let $Q_{\mu}$ denote $M/F_{\mu}(M)$.  Consider the commutative
diagram \bd \xymatrix{0\ar[r] & \theo\otimes
H^0\big(\on{gr}_{\mu_i}(M)(b)\big) \ar[r]\ar[d] & \theo\otimes
H^0\big(Q_{\mu_{i-1}}(b)\big) \ar[r]\ar[d] & \theo\otimes
H^0\big(Q_{\mu_i}(b)\big) \ar[r]\ar[d]
& 0\\
0\ar[r] & \on{gr}_{\mu_i}(M)(b)\ar[r] & Q_{\mu_{i-1}}(b) \ar[r] &
Q_{\mu_i}(b)\ar[r] & 0.} \ed Part (a) above shows that the top row
is exact.  By \eqref{gg}, the left-hand vertical arrow is
surjective.  Hence $M(b)/F_{\mu_{i-1}}(M)(b)$ is globally generated
if $M(b)/F_{\mu_i}(M)(b)$ is, and an induction then proves part (b).
So it suffices to prove \eqref{gg}.

To prove \eqref{gg}, suppose that $T$ is torsion-free and semistable
and $\mu(T)>0$.  If $T(b)$ is not globally generated, then there is
some $p\in E$ such that \bd H^0(T(b))\rightarrow H^0(T(b)\otimes
\theo_E/m_p) = T(b)\otimes \theo_E/m_p \ed is not surjective.  It
follows that there is a quotient $T(b)\xrightarrow{\phi}\theo_E/m_p$
such that \bd H^0(\phi): H^0(T(b))\rightarrow H^0(\theo_E/m_p) \cong
{\mathbf C} \ed is zero.  The long exact cohomology sequence then
implies that $\Hom(\on{ker}(\phi),\theo)^* \cong H^1(\on{ker}(\phi))
\neq 0$. Choosing a nonzero map
$\on{ker}(\phi)\xrightarrow{a}\theo$, we get a pushout diagram \bd
\xymatrix{0\ar[r] & \on{ker}(\phi)\ar[r]\ar[d]^{a} &
T(b)\ar[r]^{\phi}\ar[d]
& \theo_E/m_p\ar[r]\ar[d]^{\simeq} & 0\\
0 \ar[r] & \theo \ar[r] & \xi \ar[r] & \theo_E/m_p\ar[r] & 0.} \ed
Letting $\on{tors}(\xi)$ denote the (possibly zero) torsion subsheaf
of $\xi$, we get a nonzero map $T(b)\rightarrow \xi/\on{tors}(\xi)$,
where $\xi/\on{tors}(\xi)$ is torsion-free of slope at most $1$. But
this contradicts semistability of $T(b)$ (since $\mu(T(b))>1$),
hence it contradicts semistability of $T$.  This proves (b).

It remains to show that if $M$ satisfies (a) and (b) then $\mu_r>0$.
So suppose that $M$ satisfies (a) and (b).  It then follows from the
long exact cohomology sequence that
\begin{enumerate}
\item[(a')] $H^1(\on{gr}_{\mu_r}(M)(b)) = 0$ and
\item[(b')] $\on{gr}_{\mu_r}(M)(b)$ is globally generated
\end{enumerate}
hold as well.

Suppose that $\mu_r<0$.  Write $T= \gr_{\mu_r}(M)$.  By paragraph
0.2 of \cite{FMW}, \bd \on{deg}(T(b)) = \chi(T(b)) = \on{rk}(T) +
\on{deg}(T) \leq \on{rk}(T). \ed Since $H^1(T(b))= 0$ by (a'), we
get $h^0(T(b))\leq \on{rk}(T)$.  By (b'), we have a surjective map
$\theo\otimes H^0(T(b))\rightarrow T(b)$ of torsion-free sheaves,
necessarily of the same rank.  This map is, consequently, an
isomorphism.  But then \bd H^1(T(b)) = H^1\big(\theo\otimes
H^0(T(b))\big) \neq 0, \ed contradicting (a').  So $\mu_r>0$ after
all.  This completes the proof.
\end{proof}

\begin{prop}\label{WIT1}
Let $M$ be a torsion-free coherent sheaf on $E$.  Then
$\overline{F}^0(M)=0$ if and only if $\mu_1\leq 0$ (in the notation
of \eqref{HN filtration}).
\end{prop}
The proof is standard.

\subsection{Fourier Equivalence for $\Dlog$-Modules}\label{FT for D}

In general, given a quasicoherent sheaf of rings $\cA$ on $E$, let  $\cA
-\on{mod}$ ($\on{mod}-\cA$) denote the category of left (right)
$\cA$-modules that are quasicoherent as $\theo_E$-modules.

Recall the sheaf of rings $\D=\Dlog$ of differential operators on a
cubic curve $E$ (for a singular curve these are the log differential
operators defined in Section \ref{Dlog introduced} and studied in
Section \ref{D algebras}). We wish to calculate the Fourier
transform of $\D$ as an $\Oo_E$-bimodule (in fact as an algebra
object in the derived category). Consider $\D$ as a
quasicoherent sheaf on $E\times E$, and let $\cA = (p_{14})_*
\left(\cP\otimes p_{23}^*\D\otimes\cP^\vee\right)$ be its bimodule
Fourier transform; this is a special $D$-algebra by Corollary
\ref{special D} and
 \cite[Props. 6.2 and 6.3]{PRo}. We also let
$\cR(\D)$ and $\cR(\cA)$ denote the {\em Rees algebras} of $\D$ and
$\cA$, that is, the graded algebras given by \bd \cR(\D) =
\bigoplus_{k\geq 0} \D^k\cdot t^k \;\;\text{and}\;\; \cR(\cA) =
\bigoplus_{k\geq 0}\cA_k\cdot t^k \ed where $t$ is a formal variable
that keeps track of the grading.  Note that $\cR(\D)/t\cdot\cR(\D) =
\on{gr}(\D)$ and similarly for $\cA$. For a quasicoherent sheaf of
graded rings $\cR$ we let $\cR-\on{mod}$ denote the category of
finitely generated {\em graded} $\cR$-modules that are quasicoherent
as $\theo_E$-modules.
\begin{thm}\label{Dlog equiv}
\mbox{}
\begin{enumerate}
\item Let $(\cS_\D, \cS_\cA)$ be one of the following pairs:
\begin{enumerate}
\item $(\cS_\D, \cS_\cA) = (\D,\cA)$.
\item $(\cS_\D, \cS_\cA) = (\cR(\D), \cR(\cA))$.
\item $(\cS_\D, \cS_\cA) = (\cR(\D)/t^k\cR(\D), \cR(\cA)/t^k\cR(\cA))$ for
some $k>0$.
\end{enumerate}
Then the Fourier functor $F$ of Corollary \ref{FM equiv2} may be refined
to a functor \bd F: D^b(\cS_\D-\on{mod})\rightarrow D^b(\cS_\cA-\on{mod})
\ed that is an exact equivalence of triangulated categories.
\item Fix $k>0$.  Then the change-of-rings functors
\bd \cR(\D)-\on{mod} \rightarrow (\cR(\D)/t^k\cR(\D))-\on{mod}
\;\;\text{and}\;\; \cR(\cA)-\on{mod} \rightarrow
(\cR(\cA)/t^k\cR(\cA))-\on{mod} \ed given by $\cN\mapsto \cN/t^k\cN$
induce triangulated functors of the derived categories; moreover,
the following diagram commutes: \bd \xymatrix{
D^b(\cR(\D)-\on{mod}) \ar[r]^{F}\ar[d] & D^b(\cR(\cA)-\on{mod})\ar[d]\\
D^b((\cR(\D)/t^k\cR(\D))-\on{mod})\ar[r]^{F} &
D^b((\cR(\cA)/t^k\cR(\cA))-\on{mod}). } \ed
\end{enumerate}
\end{thm}
\begin{proof}
This follows from \cite[Theorem~6.5]{PRo} by Corollaries \ref{circle
prod} and \ref{special D}.
\end{proof}

\subsection{The Twisted Log Cotangent Bundle}\label{twisted log}

In this section we describe the Fourier transform algebra $\cA$ of
$\D$ from Theorem \ref{Dlog equiv}, and relate it to the twisted log
cotangent bundle $\Enat$ of $E$.

Let $\mu:\G\times E$ denote the action of the group $\G$ on $E$. For
a torsion sheaf $T$ on $\G$, we define the convolution action
$D^b(E)\to D^b(E)$, $\cF\mapsto T\ast \cF$ by
$$T\ast\cF=\mu_*(T\boxt \cF).$$
(Note that $\mu$ is an affine morphism and is proper on the support
of $T\boxt \cF$ for any $\cF$.)

\begin{lemma} The
Fourier transform $F(T\ast \cF)$ of convolution is canonically
identified with the tensor product $F(T)\otimes F(\cF)$.
\end{lemma}

The proof of the lemma is identical to the analogous convolution
statement for abelian varieties \cite{Mukai}, once we note the
character property of the Poincar\'e sheaf, namely the canonical
isomorphism on $\G\times E\times E$ of $(\mu\times 1)^*\cP$ and
$\cP_{12}\ot \cP_{23}$ (where $\cP_{12}$ denotes the restriction of
$\cP$ to $\G\times E$, pulled back to $\G\times E\times E$).

Next we note that the sheaf $\D$ of log differential operators on
$E$ is the $D$--algebra generated by the action of the group $\G$ on
$E$ (or specifically from that of its enveloping algebra).  Let
$U^k$ denote the $k$th filtered piece of the enveloping algebra of
$\G$, which is (as $\Oo_{\G}$--module) the $\C$--dual of functions
on the $k$th order neighborhood of the identity in $\G$ (i.e.
$U^k=\Omega_E((k+1)b)/\Omega_E$ via the residue pairing). We obtain
the following description of $\D$ as a bimodule (as for an arbitrary
group action on a variety):

\begin{lemma} Let $p_1,p_2$ denote the projections of $\G\times E$ on the factors
and $\mu$ the multiplication map to $E$. Then as an
$\Oo_E$--bimodule (sheaf on $E\times E$), $\D^k$ for any $k$ is
identified with $(\mu\times p_2)_*p_1^*U^k$.
\end{lemma}

We will now identify $\cA$ in terms of the twisted cotangent bundle
$\Enat$ of $E$ (Section \ref{about Enat}). Recall that $A$ denotes
the Atiyah extension, and $\Enatbar \overset{\on{def}}{=}
\bproj(\Sym A)$ the associated ruled surface. We let $p:\Enatbar
\rightarrow E$ denote the projection map, and $E_\infty\subset
\Enatbar$ the section at infinity.
 Note that
every other section of $\Enatbar$ has nonempty intersection with
$E_\infty$ since the Atiyah sequence is nonsplit.

\begin{thm}\label{cubic Laumon-Rothstein}\mbox{}
\begin{enumerate}
\item The bimodule Fourier transform $\cA_k$ of $\D^k$ is scheme-theoretically
supported on the diagonal and is canonically
identified with the Fourier transform $F(U^k)$ on $E$.

\item The commutative algebra $\cA$ is canonically isomorphic (as filtered $\Oo_E$--algebra)
to $p_*\theo_{\Enat}$ (functions on the twisted cotangent bundle),
inducing isomorphisms $\Enatbar = \bproj(\cR(\cA))$ and $\Enat =
\bspec(\cA)$.

\item  The Fourier transform induces an exact equivalence
of the bounded derived category of coherent $\D$-modules and the
bounded derived category of coherent sheaves on $\Enat$.

\end{enumerate}
\end{thm}

\begin{proof}
It follows from the above lemmas that $$F(\D^k\ot M)=F(U^k\ast
M)=F(U^k)\ot M$$ for any complex $M$, from which the description as
a bimodule follows. Since the algebra $\cA$ is completely determined
by the bimodule $\cA_1$ with the inclusion $\Oo_E\subset\cA_1$ as an
enveloping algebroid it follows that it is commutative (this is
immediate from the fact that the associated graded algebra is
isomorphic to $\Sym(\theo_E)$). Next note that the Fourier transform
of the nonsplit extension $U^1=\Omega_E(2b)/\Omega_E$ of $\cT_b$ by
$\Oo_b$ is isomorphic to the Atiyah extension $A$ of $\cT_E$ by
$\Oo_E$. It follows that $\cA$ is isomorphic to $p_*\theo_{\Enat}$
compatibly with the inclusions of $\theo_E$ in $\cA$ and
$p_*\theo_{\Enat}$. The isomorphisms on $\bspec$ and $\bproj$
follow. In order to fix these isomorphisms canonically, we note that
the fiber of $\bspec(\cA)$ over $b\in E$ has a canonical basepoint.
This point $d$ is characterized by the statement that the
Fourier-Mukai transform of the skyscraper at $d$ (as $\cA$-module)
is the trivial line bundle $\Oo_E$ (which is guaranteed by the
isomorphism $\Oo_d=\Oo_b$ as plain $\Oo_E$-modules) with its
canonical $\D$-module structure. On the other hand from the
definition of $A$ we find that the fiber $F_b=p\inv(b)\subset\Enat$
has a canonical identification with the cotangent fiber to $E$ at
$b$. There is now a unique isomorphism $\bspec(\cA)\to\Enat$
identifying these basepoints, fixing the algebra isomorphism above
uniquely. The final assertion follows from Theorem \ref{Dlog equiv},
noting that the coherence condition for $\D$-modules and
$\cA$-modules is identified by the Fourier transform.
\end{proof}

\subsection{Microlocalization}
We have, in addition to the algebras discussed above,
{\em microlocalizations} of the algebras $\D$ and $\cA$: these are
obtained by inverting elements that have
 invertible principal symbol and completing
with respect to the given filtration. If we microlocalize $\D$, we
obtain the sheaf $\cE$ of microdifferential operators associated to
$\D$.  Similarly, if we microlocalize $\cA$ we obtain a sheaf of
filtered algebras $\cA_\cE$ that is the sheaf of ``functions on the
punctured formal neighborhood of $E_\infty$''---that is, if we take
the structure sheaf of the formal completion of $\Enatbar$ along
$E_\infty$ and invert a local defining function for the closed
subscheme $E_\infty$ in this formal completion, the resulting sheaf
of functions is exactly $\cA_\cE$.

By \cite{AVV}, the microlocalizations are obtained as follows.  Writing
$\cR$ for the Rees ring of either algebra, we form the graded rings
$\on{gr}_{(n)} = \cR/t^n\cR$.  Each of these is a nilpotent extension of
$\on{gr} = \cR/t\cR$, and so any local lift to $\on{gr}_{(n)}$ of an invertible
element of $\on{gr}$ is invertible.
Our ring $\on{gr}_{(n)}$ is
generated by a single element, which we denote by $\del$.  We lift this
locally and invert to obtain localized rings $(\on{gr}_{(n)})_\partial$ which
form an inverse system of graded rings.  We take the inverse limit and
``reverse the formation of the Rees algebra'' to obtain a filtered ring,
the microlocalization, whose Rees algebra is this inverse limit.  Since the
graded rings at each stage are quasicoherent sheaves, one can make sense
of the Fourier transform of the microlocalization (which is not itself
quasicoherent) as the ``de-Reesed inverse limit'' of the Fourier-transformed
inverse system.
\begin{corollary}\label{microlocal rings}
We have $F(\cE) = \cA_\cE$ and $\Fin(\cA_\cE) = \cE$.
\end{corollary}
\begin{proof}
This was essentially proven in the pre-publication (arXiv) version
 of \cite{PRo},
but we sketch the proof for completeness.

The Fourier transform of the localization
$F\left((\on{gr}_{(n)}(\D))_{\partial}\right)$ is a nilpotent
extension of $F(\on{gr}(\D)_\partial) = \on{gr}(\cA)_\partial$, and
consequently any local lift of the element $\partial$ to
$F\left((\on{gr}_{(n)}(\D))_{\partial}\right)$ is invertible. The
universal property of localizations then gives us an induced
homomorphism \bd (\cR(\cA)/t^k\cR(\cA))_{\partial} =
F(\cR(\D)/t^k\cR(\D))_\partial \rightarrow
F\left((\cR(\D)/t^k\cR(\D))_{\partial}\right). \ed A standard
argument using the filtration shows that this is an isomorphism and,
consequently, we obtain an identification of the two inverse
systems. The result now follows from the construction of the
microlocalization in \cite{AVV}.
\end{proof}

In fact, we will also want slightly more: given a filtered
$\D$-module $M$ or $\cA$-module $\cF$, we may similarly form an
inverse system of graded modules $\on{gr}_{(n)}(M)$ (or similarly
for $\cF$) over the system of graded rings $\on{gr}_{(n)}(\cR)$ used
in Corollary \ref{microlocal rings}, invert a local lift of
$\partial$, and take the inverse limit.  If $M$ or $\cF$ is equipped
with a good filtration, this procedure gives us $M_\cE$ or
$\cF_{\cE}$, respectively. We can thus make sense of $F(M_\cE)$ or
$\Fin(\F_{\cE})$. With this in mind, we have the following.
\begin{corollary}\label{microlocal FM}
Suppose (for simplicity) that $M$ is a $\D$-module with good filtration
whose Fourier dual is the $\cA$-module $F$ (in some cohomological degree)
 with induced good filtration.
Then the microlocalizations $M_\cE$ and $\cF_{\cE}$
 satisfy $F(M_\cE) = \cF_{\cE}$ and
$\Fin (\cF_{\cE}) = M_\cE$.
\end{corollary}
\begin{proof}
The proof follows the same argument as in Corollary \ref{microlocal rings}.
\end{proof}

\section{Isomorphism Theorem for Moduli Spaces of $\D$-Bundles}\label{isomorphism of
moduli}

In this section we give our principal application of the Fourier-Mukai
transform.

\subsection{Moduli Stacks and Isomorphism Theorem}\label{isomorphism
theorem} We begin by defining the moduli stacks for $\Vc$-framed
spectral sheaves and $V$-framed $\D$-bundles.

As above, let $\Vc$ denote a nonzero coherent torsion sheaf on $E$ supported
on the smooth locus
$\G$.  We will let $\SSh(E,\Vc)$ denote the {\em moduli stack of
$\Vc$-framed spectral sheaves} on $\Enatbar$: its objects are pairs
$(\cF,\psi)$ consisting of
\begin{enumerate}
\item a finitely presented $S$-flat family $\cF$ of sheaves of pure dimension
$1$ on $\Enatbar$ parametrized by a scheme $S$, and
\item an isomorphism $\psi:\cF|_{E_{\infty}\times S}\rightarrow \Vc$,
\end{enumerate}
such that the restriction to $\Enatbar\times\{s\}$ is a $\Vc$-framed
spectral sheaf on $\Enatbar$ for each $s\in S$ (Definition \ref{CM
spectral sheaves def}). An isomorphism from $(\cF,\psi)$ to
$(\cF',\psi')$ consists of an isomorphism $\cF\rightarrow \cF'$ of
coherent sheaves on $\Enatbar\times S$ that is compatible with the
framings $\psi,\psi'$.

Standard arguments prove that $\SSh(E,\Vc)$ is an algebraic stack,
locally of finite type over ${\mathbf C}$.  Note that $\SSh(E,\Vc)$
will have many components, since we have not fixed one discrete
invariant of $\cF$, namely $\on{rk}(p_*\cF(-E_\infty))$.

Let $V$ denote the semistable vector bundle of degree $0$ on $E$
that is the Fourier transform of $\Vc$. Let $\BunDP(E,V)$ denote the
{\em moduli stack of $V$-framed $\D$-bundles on $E$}: its objects
are triples $(M, \{M_k\},\phi)$ consisting of
\begin{enumerate}
\item a finitely generated $S$-flat family $M$ of
locally projective $\D_S$-modules on $E\times S$,
\item a $\D$-module filtration $\{M_k\}$ of $M$ such that $M_k$ is a vector
bundle over $E\times S$ for $k\gg 0$, and
\item an isomorphism $\phi: \oplus_{k\geq N} \on{gr}_k(M)\rightarrow
\oplus_{k\geq N}\on{gr}_k(\D)_S\otimes V$ as $\on{gr}(\D)_S$-modules
for $N\gg 0$,
\end{enumerate}
such that the restriction to $E\times\{s\}$ is a $V$-framed $\D$-bundle
for each $s\in S$.
 We will refer to such a triple
as an {\em $S$-flat family of  $V$-framed $\D$-bundles} on $E$, or,
when $S$ is understood, simply as a $V$-framed $\D$-bundle.

An isomorphism between two such triples $(M,\{M_k\}, \phi)$ and
$(M',\{M_k'\},\phi')$ is an isomorphism $\psi:M\rightarrow M'$ of
$\D_S$-modules such that $\psi(M_k) = M_k'$ for all $k\gg 0$ and
$\on{gr}\psi \circ \phi|_{\on{gr}_k(M)} = \phi'\circ
\on{gr}\psi|_{\on{gr}_k(M)}$ for all $k\gg 0$.

The principal result of this section is the following isomorphism of
moduli stacks.

\begin{thm}\label{KP/CM stack equiv}
The Fourier functor $F$ induces an isomorphism of stacks \bd F:
\BunDP(E,V)\longrightarrow \SSh(E,\Vc). \ed
\end{thm}

\begin{remark} This is not, in fact, the most general such theorem possible,
in the following sense.  One may enlarge the class of $\D$-bundles and
the class of spectral sheaves in such a way that the Fourier-Mukai transform
induces an isomorphism between the moduli stacks for these more general
classes of objects.  Since these more general classes do not seem to be useful
for studying solutions of the KP or CM systems, we confine ourselves to
formulating the general theorem.  The proof uses the same techniques as
we use for Theorem \ref{KP/CM stack equiv}.
\end{remark}
\begin{thm}
Given a coherent sheaf $\cF$ on $\Enatbar$ of pure dimension one,
define the quotient sheaf $\cF_{\leq 0}$ as in
Definition \ref{F-}. Then the Fourier-Mukai
transform induces an isomorphism of the moduli stacks for the
following objects:
\begin{enumerate}
\item Pairs $(\cF,\phi)$ consisting of a coherent sheaf $\cF$ of pure dimension
one on $\Enatbar$ and an isomorphism $\phi: \cF|_{E_\infty}\rightarrow \Vc$
that satisfies:
\begin{enumerate}
\item $\cF$ has no nonzero subsheaves $\cG\subset\cF$ that satisfy both
$\cF|_{E_\infty}=0$ and $\mu(p_*\cG)>0$.
\item The quotient sheaf $\cF_{\leq 0}$ of $\cF$ is zero.
\end{enumerate}
\item Pairs $(M,\psi)$ consisting of a torsion-free $\D$-module $M$ equipped
with a good filtration and an isomorphism
\bd
\psi: \oplus_{k\geq k'}\on{gr}_k(M)\rightarrow V\otimes (\oplus_{k\geq k'}\on{gr}_k(\D))
\ed
for some $k'$ sufficiently large.
\end{enumerate}
\end{thm}


\subsection{Framed $\D$-Bundles and Framed Spectral
Sheaves} We begin by proving Theorem \ref{KP/CM stack equiv} at the
level of points: the Fourier transform identifies $V$-framed
$\D$-bundles (Definition \ref{V-framed D-bun}) with $\Vc$-framed
spectral sheaves (Definition \ref{CM spectral sheaves def}). We
continue to fix a coherent torsion sheaf $\Vc$ on $E$ that is
supported on $\G$ and $V = \Fin(\Vc)$.

We first recall a basic result on sheaves on projective bundles (see
Sections II.5 and III.8 of \cite{Hartshorne}).
 Suppose $F = \bigoplus_{k\geq n} F_k$ is a finitely
generated graded $\cR(\cA)$-module, and let $\widetilde{F}$ denote
the corresponding sheaf on $\Enatbar = \bproj(\cR(\cA))$.
Conversely, given a coherent sheaf $\cF$ on $\Enatbar$, let
$\Gamma\cF = \bigoplus_{k \geq 0} \Gamma_k\cF \overset{\on{def}}{=}
\bigoplus p_*\cF(kE_\infty)$ denote the associated
$\cR(\cA)$-module.
\begin{prop}\label{Serre theorem}
Suppose $S$ is a noetherian scheme, and let $\cA_S =
\cA\boxtimes\theo_S$ on $E\times S$.  Then:
\begin{enumerate}
\item If $F$ is a finitely generated $\cR(\cA_S)$-module, then there is a map
$F\rightarrow \Gamma\widetilde{F}$ of graded $\cR(\cA_S)$-modules
that is functorial in $F$ and is an isomorphism in sufficiently high
degrees.
\item If $\cF$ is a coherent sheaf on $\Enatbar\times S$, then there is an isomorphism
$\widetilde{\Gamma\cF}\rightarrow \cF$ of coherent sheaves that is
functorial in $\cF$.
\item The functor $N\mapsto \widetilde{N}$ is exact and commutes with tensor product
over $\theo_S$.
\item The functor $\Gamma$ takes exact sequences of coherent sheaves on $\Enatbar\times S$
to complexes of graded $\cR(\cA_S)$-modules that are exact in all
sufficiently high graded degrees.
\end{enumerate}
\end{prop}

Suppose that $M$ is a $V$-framed $\D$-bundle; in particular, it
comes equipped with a canonical structure of filtered $\D$-module
(see Definition \ref{canonical filt}).  Let $\cR(M)$ denote the Rees
module of $M$.

We will need the following basic fact relating $M$ and $\cR(M)$:
\begin{prop}\label{van o thm}
\mbox{}
\begin{enumerate}
\item Suppose that $M$ is a $V$-framed $(\D)_S$-bundle on $E$ for
some scheme $S$. Then the Rees module $\cR(M)$ is an $S$-flat family
of
 torsion-free graded $\cR(\D)$-modules,
and the framing becomes an isomorphism $\cR(M)/t\cR(M)
\xrightarrow{\phi} \on{gr}(\D)_S\otimes V$ in high degree as
$\on{gr}(\D)_S$-modules.
\item Conversely, if $N = \oplus N_k$ is an $S$-flat family of
torsion-free graded $\cR(\D)$-modules equipped with a $V$-framing
$N/tN \xrightarrow{\phi} \on{gr}(\D)_S\otimes V$, then
$\underset{\longrightarrow}{\lim} N_k$ together with the induced
filtration, filtered $\D$-module structure and framing is an
$S$-flat family of
 $V$-framed $\D$-bundles.
\end{enumerate}
These constructions give an equivalence between the groupoid of
$V$-framed $\D$-bundles and the groupoid of $V$-framed torsion-free
$\cR(\cD)$-modules.
\end{prop}
\begin{proof}
See \cite{LvO}.
\end{proof}

We will first check that $M$, together with its filtration, transforms
to a filtered sheaf on $\Enatbar$.
\begin{lemma}\label{D bundle to sheaf}
Suppose that $M$ is a $V$-framed $\D$-bundle.  Then:
\begin{enumerate}
\item $F(\cR(M))$ is a sheaf in cohomological degree $1$.
\item The element $t\in\cR(\cA)$ is a non-zero-divisor on $F(\cR(M))$.
\end{enumerate}
\end{lemma}
\begin{proof}
Suppose $(M, M_k,\phi)$ is a $V$-framed $\D$-bundle. Part (i) of
Proposition \ref{can filt} implies
that, for $m$ sufficiently large, the canonical filtration $\Theta$
on $M$ satisfies $\Theta_k(M) = M_k$ for all $k\geq m$.  By parts
(a) and (i) of Proposition \ref{can filt}, we find that \bd
\Theta_k(M) \subseteq \Theta_k\big(M_{\cE}\big)/
\Theta_{-1}\big(M_{\cE}\big) \ed for all $k$.   Part (iii) of
Proposition \ref{can filt} implies that
$\Theta_k\big(M_{\cE}\big)/\Theta_{-1}\big(M_{\cE}\big)$ is an
iterated extension of copies of $V$; since $F^0(V) = 0$, we also
have $F^0(M_k)=0$.

Since $F^0(V)=0$, we find that $F(M_k) = F^1(M_k)\rightarrow
F^1(M_{k+1}) = F(M_{k+1})$ is injective for all $k$ in the great
range.  This implies (2).
\end{proof}

\begin{corollary}
The exact sequence \bd 0\rightarrow t\cR(M)\rightarrow
\cR(M)\rightarrow \on{gr}(\D) \otimes V \rightarrow 0 \ed of graded
$\cR(\D)$-modules in sufficiently high degree has as its Fourier
transform a short exact sequence
\begin{equation}\label{dual framing seq}
0\rightarrow tF(\cR(M))\rightarrow F(\cR(M))\rightarrow
\on{gr}(\cA)\otimes \Vc\rightarrow 0
\end{equation}
of graded $\cR(\cA)$-modules in cohomological degree $1$.
\end{corollary}
\begin{proof}
It follows from Lemma \ref{D bundle to sheaf} that the exact
sequence \bd 0\rightarrow t\cR(M)\rightarrow \cR(M) \rightarrow
\cR(M)/t\cR(M) = \on{gr}(\D)\otimes V \rightarrow 0 \ed is taken by
the Fourier functor to an exact sequence \bd 0\rightarrow
tF(\cR(M))\rightarrow F(\cR(M)) \rightarrow F(\cR(M)/t\cR(M)) =
F(\on{gr}(\D)\otimes V)\rightarrow 0 \ed of sheaves in cohomological
degree $1$.  So it suffices to check that $F(\on{gr}(\D)\otimes V) =
\on{gr}(\cA)\otimes\Vc$ as $\gr(\cA)$-modules. This is immediate
from part (3) of Theorem \ref{Dlog equiv}.
\end{proof}

Write $N_k = F(M_k)$ and $N = F(\cR(\cM))$, and let $\cF =
\widetilde{F(\cR(M))} = \widetilde{N}$ (see \ref{Serre theorem}).
Taking the sequence of coherent sheaves on $\Enatbar$ associated to
the sequence \eqref{dual framing seq} gives an exact sequence \bd
0\rightarrow \cF(-E_\infty) \rightarrow \cF \rightarrow
\Vc\rightarrow 0 \ed on $\Enatbar$; in other words, $\cF$ is a
$\Vc$-framed coherent sheaf on $\Enatbar$.

\begin{prop}\label{WIT properties}
If $M$ is a $V$-framed $\D$-bundle, then $\cF =
\widetilde{F(\cR(M))}$ (with its induced $\Vc$-framing) is a
$\Vc$-framed CM spectral sheaf.
\end{prop}
\begin{proof}
Observe first that $\cF$ has one-dimensional support: indeed,
 since
$M_k \subset M_{\cE,k}/M_{\cE,-1}$ is torsion-free of negative degree, a
standard Fourier-Mukai computation shows that $F(M_k)$ has support equal
to $E$, and thus
$\on{supp}(\cF)$ has dimension at least one.

Suppose $\cF'$ is a subsheaf of $\cF$ of dimension $0$. By part (1)
of Proposition \ref{Serre theorem} and part (2) of Lemma \ref{D
bundle to sheaf}, the natural map $p_*\cF(kE_\infty)\rightarrow
p_*\cF((k+1)E_\infty)$ is injective for $k$ sufficiently large. On
the other hand, if the support of $\cF'$ had nontrivial intersection
with $E_\infty$,  the natural maps $p_*\cF'(kE_\infty)\rightarrow
p_*\cF'((k+1)E_\infty)$ would fail to be injective for all $k$.
Since $p_*\cF'(kE_\infty)\subseteq p_*\cF(kE_\infty)$
 for all $k$, it follows that $\on{supp}(\cF')\cap E_\infty =\emptyset$.
Then the maps $p_*\cF'(kE_\infty)\rightarrow p_*\cF'((k+1)E_\infty)$
are isomorphisms for all $k$, implying that $\cup_{k\geq
0}p_*\cF'(kE_\infty)$ is a finite-length $\theo_E$-submodule of
$F(M)$.  This transforms under $\Fin$ to an $\theo$-coherent
$\D$-submodule of $M$, a contradiction since $M$ is torsion-free
over $\D$. Thus $\cF$ is of pure dimension $1$.

It remains to prove the normalization conditions (i) and (ii) in
Definition \ref{CM spectral sheaves def}.  By hypothesis we have that
$\on{deg}(M_k) = -n$ and $\on{rk}(M_k) = (k+1)\on{rk}(V)$ for all
$k$ sufficiently large; a standard Fourier-Mukai computation then shows that
$F(M_k) = N_k$ has rank $n$ and degree $(k+1)\on{rk}(V)$ for $k\gg 0$.
Condition (i) then follows from Lemma \ref{CM equiv cond}.

To verify Condition (ii), we use the exact sequence
\begin{equation}\label{eq for M}
0 \rightarrow M_k \rightarrow M_{\cE,k}/M_{\cE,-1} \rightarrow Q_k \rightarrow 0\end{equation}
for $k\gg 0$
in which, by Definition \ref{V-framed D-bun} and
Proposition \ref{can filt props cor}, $Q_k$ is a torsion $\theo_E$-module
supported on the smooth locus of $E$.  Applying $F$, we get an exact
sequence
\begin{equation}\label{eq for F}
0\rightarrow F^0(Q_k)\rightarrow N_k \rightarrow N_k/N_{-1}\rightarrow 0.
\end{equation}
For $k\gg 0$, moreover, we have $N_k = p_*\cF(kE_\infty)$ and
$N_k/N_{-1} \cong p_*(\cF(kE_\infty)/\cF(-E_\infty))$ by Theorem
\ref{Dlog equiv}.  It follows that for $k\gg 0$,
$p_*\cF(-E_\infty) = F(Q_k)$, a vector bundle $W$ of the kind required by
Condition (ii) of Definition \ref{CM spectral sheaves def}.  This proves
the proposition.\end{proof}

We next study the transform of a CM spectral sheaf.
\begin{defn}\label{F-}
Let \bd F_{\leq 0} \overset{\on{def}}{=}
\frac{\big(\Gamma\cF/\on{tors}(\Gamma\cF)\big)}
{F_{>0}\big(\Gamma\cF/\on{tors}(\Gamma\cF)\big)}. \ed The quotient
sheaf $\cF\rightarrow \cF_{\leq 0}$ is defined by $\cF_{\leq 0} =
\widetilde{F_{\leq 0}}$.
\end{defn}
\begin{remark}
Note that, by Proposition \ref{Serre theorem}, we have $F_{\leq 0}
\rightarrow \Gamma\cF_{\leq 0}$ which is an isomorphism in high
degree.
\end{remark}

\begin{prop}\label{SS defs agree}
Every CM spectral sheaf $(\cF,\phi)$ has the following two
properties:
\begin{enumerate}
\item $\cF$ has no nonzero subsheaves $\cG\subset \cF$
satisfying both $\cG|_{E_\infty}=0$ and $\mu(p_*\cG)>0$.
\item The quotient sheaf $\cF_{\leq 0}$ of $\cF$ is zero.
\end{enumerate}
\end{prop}

\begin{proof}
Suppose that $\cG\subset \cF$ is a subsheaf satisfying
$\cG|_{E_\infty}=0$. Then $p_*\cG=p_*\cG(-E_\infty)\subset
p_*\cF(-E_\infty)$, implying, since the last term is a semistable
vector bundle of degree $0$, that $\on{deg}(p_*\cG)\leq 0$.  This
proves statement (1) above. For statement (2), we proceed as
follows. For every $k$, we have a map $\phi_k: F_{-1}\rightarrow
(F_k)_{\leq 0}$ from a semistable vector bundle of degree $0$ to a
sheaf whose Harder-Narasimhan subquotients all have nonpositive
degrees; a standard argument shows that the image is a semistable
torsion-free sheaf of degree $0$ with torsion-free quotient.
 From the short exact
sequence $0\rightarrow F_{-1}\rightarrow F_k\rightarrow
Q_k\rightarrow 0$ for all $k\geq 0$, where $Q_k$ is a torsion
$\theo$-module, it then follows that the induced map from $Q_k$ to
$(F_k)_{\leq 0}/\on{Im}(\phi_k)$ is zero.  Since $F_k$ surjects onto
$(F_k)_{\leq 0}$, however, it follows that $\phi_k$ must already be
surjective, i.e. $(F_k)_{\leq 0}$ is torsion-free semistable of
degree $0$, and $F_{-1}\rightarrow (F_k)_{\leq 0}$ is surjective.
The filtration on $(F_k)_{\leq 0}$ induced from the one on
$F_{-1}$ by part
(ii) of \cite[Lemma~1.2.5]{FriedMorgminuscule}
then has
subquotients that are line bundles of degree $0$.

Suppose, then, that $\cF_{\leq 0}$ is nonzero.  Then, for all $k$
sufficiently large, $p_*\cF_{\leq 0}(kE_\infty) = (F_k)_{\leq 0}$.
We claim that for all $k$ sufficiently large, $p_*\cF_{\leq
0}(kE_\infty) = p_*\cF_{\leq 0}((k+1)E_\infty)$: if not, then, since
for $k$ large $\on{rk}(p_*\cF_{\leq 0}(kE_\infty)) =
\on{rk}(p_*\cF_{\leq 0}((k+1)E_\infty))$, we find that
 $\on{deg}(F_k)_{\leq 0}$ is a strictly increasing function of $k$, a
contradicition since it is bounded above.  It follows that, for $k$
sufficiently large, the terms in the graded $\cA$-module $\oplus_m
(F_m)_{\leq 0}$ stabilize, and thus that there is a morphism---the
multiplication map---of the form $\At\otimes (F_k)_{\leq 0}
\rightarrow (F_k)_{\leq 0}$ that is the identity map on the subsheaf
$\theo\otimes (F_k)_{\leq 0}$---in other words, the extension given
by tensoring the Atiyah extension with $(F_k)_{\leq 0}$ is split.
But we have already seen that $(F_k)_{\leq 0}$ is a successive
extension of line bundles on $E$, and is in particular a vector
bundle on $E$.  Hence, taking the trace on $\on{Ext}^1((F_k)_{\leq
0}, (F_k)_{\leq 0})$ gives a splitting (up to scale)
 of the inclusion map
\begin{displaymath}
\on{Ext}^1(\theo,\theo) \rightarrow \on{Ext}^1((F_k)_{\leq 0},
(F_k)_{\leq 0})
\end{displaymath}
that tensors an extension of $\theo$ by $\theo$ with the vector
bundle $(F_k)_{\leq 0}$.  In particular, the Atiyah extension {\em
cannot} split when tensored with $(F_k)_{\leq 0}$ if the latter is
nonzero.  This completes the proof.
\end{proof}

\begin{prop}\label{WIT inv}
Suppose $\cF$ is a $\Vc$-framed spectral sheaf on $\Enatbar$. Then
$\Fin(\Gamma\cF)$ corresponds, under the equivalence of Proposition
\ref{van o thm}, to a  $V$-framed $\D$-bundle on $E$.
\end{prop}
\begin{proof}
Since $\cF$ is of pure dimension $1$ and $\cF|_{E_\infty} = \Vc$,
$\cF$ has no local sections supported set-theoretically on
$E_\infty$. Thus the map $\cF(kE_\infty)\rightarrow
\cF((k+1)E_\infty)$ is injective for all $k$, and the action of $t$
on $\Gamma\cF$ is regular. Taking the direct image of \bd
0\rightarrow \bigoplus\cF(kE_\infty)\rightarrow \bigoplus
\cF((k+1)E_\infty) \rightarrow
\bigoplus\cF((k+1)E_\infty)|_{E_\infty} = \Vc(kE_\infty)\rightarrow
0 \ed thus gives an exact sequence
\begin{equation}\label{Fseq}
0\rightarrow t\Gamma\cF \rightarrow \Gamma\cF \rightarrow \Gamma\cF/t\Gamma\cF
 = \on{gr}(\cA)\otimes \Vc \rightarrow 0
\end{equation}
in high degree.

Consider the torsion-free quotient
$\Gamma_k\cF/\on{tors}(\Gamma_k\cF)$. By Proposition \ref{SS defs agree},
 for $k\gg 0$ all terms in
its Harder-Narasimhan filtration have positive slope; applying
Proposition \ref{WIT0} and using the exact sequence
\begin{equation}\label{torsseq}
0\rightarrow \on{tors}(\Gamma_k\cF)\rightarrow \Gamma_k\cF
\rightarrow \Gamma_k\cF/\on{tors}(\Gamma_k\cF)\rightarrow 0,
\end{equation}
we find that $\Fin(\Gamma_k\cF)$ is a sheaf in cohomological degree
$0$ for $k\gg 0$.  It then follows that the same is true for every
term of \eqref{Fseq} for $k\gg 0$, and we obtain an exact sequence
\bd 0\rightarrow t\Fin(\Gamma\cF) \rightarrow \Fin(\Gamma\cF)
\rightarrow \on{gr}(\D)\otimes V\rightarrow 0 \ed in high degrees.
In particular, $M = \cup \Fin(\Gamma_k\cF)$ is a filtered
$\D$-module equipped with an isomorphism $\on{gr}(M) =
\on{gr}(\D)\otimes V$ of $\on{gr}(\D)$-modules in high degree.

We wish to prove next that $M$ is torsion-free; so suppose not. The
$V$-framing implies that the torsion submodule $\on{tors}(M)$ is
contained in $M_k =\Fin(\Gamma_k\cF)$ for some $k$ sufficiently
large.  In particular, $\on{tors}(M)$ is an $\theo$-coherent
$\D$-submodule of $M$.  Since $F(M) = p_*(\cF|_{\Enat})$ is a sheaf
concentrated in cohomological degree $1$, we find that
$F(\on{tors}(M))$ is an $\theo$-coherent sheaf concentrated in
degree $1$ that gives a subsheaf $\cG =
\widetilde{F(\cR(\on{tors}(M)))}$ of $\cF$ with the property (since
$\on{tors}(M)_k = \on{tors}(M)_{k+1}$ for $k \gg 0$) that
$\cG|_{E_\infty}=0$.  By Proposition \ref{Serre theorem}, $p_*(\cG)
= F(\on{tors}(M)_k)$ for $k\gg 0$; this sheaf is torsion-free since
it is the direct image of a sheaf of pure dimension $1$ that has
zero intersection with $E_\infty$.  Because $\Fin(p_*\cG) =
\on{tors}(M)_k$ in degree $0$, Proposition \ref{WIT0} then implies
that $p_*\cG$ has positive slope, contradicting part (1) of
Proposition \ref{SS defs agree}.  So $\on{tors}(M)= 0$.

By Condition (ii) of Definition \ref{CM spectral sheaves def}, we
get $\on{rk}(M_k) = (k+1)\on{rk}(V)$ for all $k\gg 0$.  Since $M$ is
$V$-framed and $\on{gr}(M)$ is torsion-free, we conclude that $M$ is
normalized if and only if $M_{-1}=0$.  So, suppose that $M_{-1}\neq
0$; it would then follow that $\on{rk}(M_k) < \on{rk}(V)$ for all
$k<0$.  Let $k$ be the smallest integer such that $M_k\neq 0$, and
consider the $\D$-submodule $N$ of $M$ generated by $M_k$.  This is
a $\D$-bundle of rank $\on{rk}(M_k)$. Then, under the induced
filtration from $M$, $\on{rk}(\on{gr}_\ell(N)) = \on{rk}(M_k)$ for
all $\ell\geq k$.  Thus, shifting the filtration by $k$, we find
that $N(k)$ is a normalized $\D$-bundle. Let $\cG$ be the Fourier
transform of $N$.  It follows from our earlier argument that
$\on{deg}(p_*(\cG((-1+k)E_{\infty})))=0$, hence
$\on{deg}(p_*(\cG(-E_\infty))) = (-k)\on{rk}(M_k)$.  Since $k<0$ and
$p_*\cG(-E_\infty)$ is a subsheaf of $\Gamma_{-1}(\cF)$, this
contradicts semistability of $\Gamma_{-1}(\cF)$.  So $M$ satisfies
Condition (ii) of Definition \ref{V-framed D-bun}.

To prove that Condition (i) of Definition \ref{V-framed D-bun} holds, we
use the exact sequence
\begin{displaymath}
0\rightarrow \Gamma_{-1}\cF \rightarrow \Gamma_k\cF \rightarrow
p_*(\cF(kE_\infty)/\cF(-E_\infty))\rightarrow 0.
\end{displaymath}
Using Corollary \ref{microlocal FM}, this transforms under
Fourier-Mukai to the exact sequence
\begin{displaymath}
M_k \rightarrow M_{\cE,k}/M_{\cE,-1} \rightarrow F(\Gamma_{-1}\cF)\rightarrow
0
\end{displaymath}
for $k\gg 0$.  By Condition (i) of Definition \ref{CM spectral sheaves def},
we have that $F(\Gamma_{-1}\cF)$ is a coherent $\theo_E$-torsion sheaf
supported on the smooth locus of $E$, from which it follows that
$M_k\cong M_{\cE,k}/M_{\cE,-1}$ as filtered sheaves in a neighborhood of
$\infty$.  The desired conclusion is then immediate.
\end{proof}

We may summarize the results of Propositions \ref{WIT properties}
and \ref{WIT inv} as follows.
\begin{corollary}\label{pointwise corr}
The construction $M\mapsto \widetilde{F(\cR(M))}$ gives a bijective
correspondence between $V$-framed $\D$-bundles on $E$ and
$\Vc$-framed spectral sheaves on $\Enatbar$.
\end{corollary}

\subsection{Proof of Theorem \ref{KP/CM stack equiv}}
Now that we have proven the isomorphism for ${\mathbf
C}$-points of the stacks (Corollary \ref{pointwise corr}),
 what remains is simply to extend this bijection to
families and check that it is functorial and compatible with base
change.

By standard limit arguments, we may assume that all parameter
schemes $S$ are noetherian; we will do so without comment below.

The proof will require the following facts from \cite{Bridgeland}.
\begin{lemma}\label{Br lemma}
Let $S$ and $T$ be schemes.
\begin{enumerate}
\item
Let $g:T\rightarrow S$ be a morphism and fix $\cE\in D(S\times E)$
of finite Tor-dimension over $S$.  Then there exists an isomorphism
\bd F\circ {\mathbf L}(g\times 1_E)^*\cE\simeq {\mathbf L}(g\times
1_E)^*\circ F(\cE). \ed
\item Let $\cE$ be a coherent sheaf on $S\times E$ that is flat over $S$.
Suppose that for each $s\in S$, $F(\cE_s)$ is a sheaf (in some
cohomological degree).  Then $F(\cE)$ is a sheaf on $S\times E$ that
is flat over $S$.
\end{enumerate}
\end{lemma}
\noindent
The isomorphism of Theorem \ref{KP/CM stack equiv}
 will follow from a certain collection of technical facts
(Lemmas \ref{F(M) is a sheaf} through \ref{claim needing no proof}).

Given a $V$-framed $\D$-bundle $(M,\{M_k\},\phi)$ on $E\times S$, we
obtain an object $F(\cR(M))$ of the derived category of
$\cR(\cA)$-modules.

\begin{lemma}\label{F(M) is a sheaf}
$F(\cR(M))$ is an $S$-flat sheaf of $\cR(\cA)$-modules concentrated
in degree $1$.
\end{lemma}
\begin{proof}
By assumption, $\cR(M)$ is $S$-flat.  By part (2) of Lemma \ref{Br
lemma} and part (1) of Lemma \ref{D bundle to sheaf}, the result
follows.
\end{proof}

It then follows from part (3) of Proposition \ref{Serre theorem}
that $\cF = \widetilde{F(\cR(M))}$ is an $S$-flat family of coherent
sheaves on $\Enatbar$.  We next prove:

\begin{lemma}\label{F(M) restriction}
$\widetilde{F\big(\cR(M)\big)}|_{\Enatbar_s} =
\widetilde{F\big(\cR(M|_{\Enatbar_s})\big)}$ for all $s\in S$.
\end{lemma}
\begin{proof}
By Lemma \ref{F(M) is a sheaf} and part (1) of Lemma \ref{Br lemma},
we have $F(\cR(M)_s) = F(\cR(M))_s$ for all $s\in S$.  The lemma
then follows by part (3) of Proposition \ref{Serre theorem}.
\end{proof}

\begin{lemma}\label{FM of R(M) seq}
The Fourier transform of the exact sequence
\begin{equation}\label{R(M) sequence}
0\rightarrow t\cR(M)\rightarrow \cR(M)\rightarrow \on{gr}(\D)\otimes
V\rightarrow 0
\end{equation}
in high graded degrees becomes a short exact sequence
\begin{equation}\label{F(R(M)) sequence}
0\rightarrow \cF(-E_\infty)\rightarrow \cF\rightarrow s_*\Vc
\rightarrow 0
\end{equation}
of sheaves on $\Enatbar$ (where $s: E\rightarrow \Enatbar$ is the
section at infinity).
\end{lemma}
In particular, it then follows by Corollary \ref{pointwise corr}
that $(\cF, F(\phi))$ is an $S$-flat family of $\Vc$-framed spectral
sheaves.

\begin{proof}[Proof of Lemma \ref{FM of R(M) seq}]
The Fourier transforms of the terms in \eqref{R(M) sequence} are all
sheaves in cohomological degree $1$ by Lemma \ref{F(M) is a sheaf}.
Hence, by the long exact cohomology sequence, $F$ applied to
\eqref{R(M) sequence} is an exact sequence of sheaves in degree $1$
in high graded degrees.  We get an exact sequence \bd 0\rightarrow
F(\cR(M))(-1)\rightarrow F(\cR(M))\rightarrow \on{gr}(\cA)\otimes
\Vc\rightarrow 0 \ed as a result.  Part (3) of Proposition
\ref{Serre theorem} then proves the lemma.
\end{proof}

It is clear from the construction that this takes isomorphisms of
$V$-framed $\D$-bundles to isomorphisms of $\Vc$-framed spectral
sheaves.  Moreover, by Lemma \ref{F(M) is a sheaf} and part (1) of
Lemma \ref{Br lemma}, this construction commutes with pull back
along morphisms $S'\rightarrow S$.  Hence it gives a morphism of
stacks $F: \BunDP(E,V)\rightarrow \SSh(E,\Vc)$.

In the other direction, we start with an $S$-flat family of
$\Vc$-framed spectral sheaves $(\cF,\psi)$.  We will prove:

\begin{lemma}\label{pushforward props}
Applying $\Gamma$ to
\begin{equation}\label{upstairs}
0\rightarrow \cF(-E_\infty)\rightarrow \cF \xrightarrow{\phi}
\Vc\rightarrow 0
\end{equation}
gives an exact sequence
\begin{equation}\label{Gammas}
0\rightarrow \Gamma\cF(-1) \rightarrow \Gamma\cF\rightarrow
\on{gr}(\cA)\otimes \Vc\rightarrow 0
\end{equation}
of $S$-flat graded $\cR(\cA)$-modules in all sufficiently high
graded degrees.  Moreover, for any morphism $S'\rightarrow S$ we
have $\Gamma_k(\cF_{S'}) = (\Gamma_k\cF)_{S'}$ for all $k$
sufficiently large.
\end{lemma}

\begin{proof}
The exactness is immediate from part (4) of Proposition \ref{Serre
theorem}; indeed, the sequence is exact in graded degree $k$ by
construction whenever \bd {\mathbf R}^1p_*\cF((k-1)E_\infty) = 0.
\ed

For flatness over $S$, we begin with the observation that it is
enough to check this flatness locally over $E$.  In particular, we
may restrict attention
 to
an open set of $E$ over which $\Enatbar$ is a trivial ${\mathbf
P}^1$-bundle and choose a \v{C}ech covering of ${\mathbf P}^1$ by
two open sets. We then obtain a two-term \v{C}ech complex
$C^{\bullet}$ over (an open set of) $\Enatbar\times S$; since $\cF$
is $S$-flat, so are the terms of the complex $\cF(kE_\infty)\otimes
C^{\bullet}$.  Thus the terms of
\begin{equation}\label{Cech complex}
p_*\cF(kE_\infty)\otimes C^0 \rightarrow p_*\cF(kE_\infty)\otimes
C^1
\end{equation}
are also $S$-flat.  Choosing $k$ sufficiently large that ${\mathbf
R}^1p_*\cF(kE_\infty) = 0$, we find that $\Gamma_k\cF =
p_*\cF(kE_\infty)$ is the kernel of the surjective map \eqref{Cech
complex} of $S$-flat sheaves, hence is itself $S$-flat. Moreover, if
$S'\rightarrow S$ is any morphism, the pullback of \eqref{Cech
complex} to $S'$ has $\Gamma_k(\cF_{S'})$ as its kernel (by
$S$-flatness of the terms of the exact sequence), completing the
proof.
\end{proof}

We next prove:

\begin{lemma}\label{Gamma exactness}
In all sufficiently high degrees, the Fourier transform of
\eqref{Gammas} is an exact sequence
\begin{equation}\label{seq of R(D)-mods}
0\rightarrow \Fin(\Gamma(\cF))(-1)\rightarrow
\Fin(\Gamma\cF)\rightarrow \on{gr}(\D)\otimes V\rightarrow 0
\end{equation}
 of $S$-flat families of $\cR(\D)$-modules concentrated in
cohomological degree $0$.
\end{lemma}
By Lemma \ref{pushforward props}, Lemma \ref{Gamma exactness}, and
part (1) of Lemma \ref{Br lemma}, the formation of the exact
sequence \eqref{seq of R(D)-mods} then commutes with arbitrary base
changes $S'\rightarrow S$.

\begin{proof}[Proof of Lemma \ref{Gamma exactness}]
By Lemma \ref{pushforward props}, part (2) of Lemma \ref{Br lemma},
and Corollary \ref{pointwise corr}, $\Fin(\Gamma\cF)$ is, in
sufficiently high graded degrees, an $S$-flat sheaf in cohomological
degree $0$.  Since this applies to $\Fin$ of all terms in
\eqref{Gammas}, it follows from the long exact cohomology sequence
that $\Fin$ applied to \eqref{Gammas} is an exact sequence of the
form \eqref{seq of R(D)-mods} in all sufficiently high graded
degrees.
\end{proof}

We obtain, as an immediate consequence,
 the following analog of Proposition \ref{van o thm}:

\begin{lemma}\label{claim needing no proof}
The natural map $\Fin(\iota_k):\Fin(\Gamma_k\cF)\rightarrow
\Fin(\Gamma_{k+1}\cF)$ is injective for all $k$ sufficiently large
and satisfies $\Fin(\iota_k)_{S'} = \Fin(\iota_{k,S'})$ for every
$S'\rightarrow S$.
\end{lemma}
It follows that
 $M = \cup_k \Fin(\Gamma_k\cF)$ is a $\D$-module satisfying
\begin{enumerate}
\item  $M_{S'} = \cup_k \Fin(\Gamma_k(\cF_{S'}))$ for every $S'\rightarrow S$,
\item  $M_k = \Fin(\Gamma\cF_k)$ defines a good filtration on $M$, and
\item $\bigoplus_{k\geq n}\on{gr}_k(M) = \bigoplus_{k\geq n}\on{gr}(\D)\otimes V$
as $\on{gr}(\D)$-modules for some $n$ sufficiently large.
\end{enumerate}
Combining these observations with Corollary \ref{pointwise corr}, we
find that $(M,\{M_k\}, \Fin(\psi))$ is an $S$-flat family of
$V$-framed $\D$-bundles on $E$.

It is again clear from the constructions that the map $(\cF,\psi)
\mapsto (M,\{M_k\},\Fin(\psi))$ takes isomorphisms of $\Vc$-framed
spectral sheaves to isomorphisms of $V$-framed $\D$-bundles.  We
have proven above the compatibility of this construction with base
change, so we obtain a morphism of stacks $\Fin:
\SSh(E,\Vc)\rightarrow \BunDP(E,V)$.

To see that these functors give isomorphisms of stacks, let $
(\cF,\psi)= F(M,\{M_k\},\phi)$. Applying $\Gamma$ to $\cF$ we obtain
an exact sequence of
 $\cR(\cA)$-modules of the form \eqref{Gammas} in which, by
Proposition \ref{Serre theorem}, the left-hand and middle terms are
isomorphic to $F(\cR(M))(-1)$ and $F(\cR(M))$ respectively in high
degrees.  Applying $\Fin$ to this exact sequence we then obtain, by
the construction and Theorem \ref{FM equiv}, an exact sequence
equipped with a canonical isomorphism to \eqref{R(M) sequence} in
high degrees; this reconstructs $(M,\{M_k\},\phi)$ if we forget the
terms $M_k$ in the filtration for small $k$, and hence is isomorphic
to $(M,\{M_k\},\phi)$ as a family of framed $\D$-bundles.  Thus
$\Fin\circ F \simeq 1$.

Similarly, starting from $(\cF,\psi)$ we apply $\Gamma$ to the
sequence \eqref{upstairs} to obtain \eqref{Gammas}; applying $\Fin$
to this we obtain a short exact sequence of the form \eqref{seq of
R(D)-mods} and thereby a flat family of $V$-framed $\D$-bundles
$(M,\{M_i\},\phi)$.  Applying $F$ to these data then returns the
complex \eqref{Gammas} in sufficiently high degrees by construction
and Theorem \ref{FM equiv}.  Finally, taking the associated coherent
sheaves on $\Enatbar\times S$ gives us \eqref{upstairs} by
Proposition \ref{Serre theorem}. Thus $F\circ\Fin \simeq 1$.

This completes the proof of Theorem \ref{KP/CM stack
equiv}.\hfill\qedsymbol

\subsection{Compatibility of KP and CM Hierarchies}\label{proof of
compatibility}

We now complete the proof of the compatibility between KP and CM
hierarchies, stated in Theorem \ref{compatibility}.

\begin{proof}[Proof of Theorem \ref{compatibility}]
We first establish that $\F$ identifies the KP algebroid
$\underline{\End}_\cE$ on $\BunDP(E,V)$ with the CM algebroid, also
denoted $\underline{\End}_\cE$, on $\CM(E,\Vc)$. Let $M,\cF$ denote
a $\D$-bundle and its Fourier transform spectral sheaf. By
Corollaries \ref{microlocal rings} and \ref{microlocal FM}, the
Fourier-Mukai transform is compatible with microlocalization,
 sending $M_\cE$ to $\cF_\cE$.
It follows that endomorphisms of each, with their Lie bracket, are
also identified. The definitions of the respective algebroid
structures (i.e. the actions deforming $M,\cF$ inside their
microlocalization by multiplication) are then also clearly
identified. Finally, the choice of a micro-oper structure $\del_M$
acting on $M_\cE$ corresponds to the choice of endomorphism $\xi$ of
$\cF_\cE$, and the restriction on the symbol of $\del_M$ is
precisely the condition making $\xi$ a Higgsing of $\cF$. Since the
KP flows on micro-opers and CM hierarchy on Higgsed spectral
sheaves are given by the action of the corresponding algebroids, the
compatibility of hierarchies is established.
\end{proof}

\section{$\D$-Lattices and d-Lattices}\label{D and d}
Our goal in this section is to prove an extension of Sato's
$\D$-module description of the big cell $\GRon$. Namely, we will
show that the entire Sato Grassmannian is a moduli space for certain
$\D$-submodules of $\cE^n$, the $\D$-lattices.

\subsection{d-Lattices}\label{d-lattice section}
We begin by recalling the scheme structure on $\GR_n$, following
\cite{AMP}.  In \cite{AMP}, it is proven that $\GR_n$ is an
infinite-dimensional scheme that represents a certain functor which
we now describe.

Let $z$ be a formal parameter.  Recall (Section \ref{Sato}) that a
vector subspace $B\subset \boldc(\!(z\inv)\!)^n$ is called a {\em
c-lattice} if there exist integers $k$ and $\ell$ such that
$\big(\boldc [\![z\inv]\!] z^k\big)^n \subseteq B \subseteq
\big(\boldc [\![z\inv]\!]z^\ell\big)^n$.

For a scheme $S$, a {\em discrete sub-$\theo_S$-module} (or {\em
d-lattice}) $L\subset \theo_S(\!(z\inv)\!)^n$ is a quasi-coherent
$\theo_S$-module $L$ equipped with an injection of $\theo_S$-modules
in $\theo_S(\!(z\inv)\!)^n$ such that
\begin{enumerate}
\item For every morphism $S'\xrightarrow{f} S$, the base-changed
(composite) morphism \bd f^*L \rightarrow
f^*(\theo_S(\!(z\inv)\!)^n) \rightarrow \theo_{S'}(\!(z\inv)\!)^n
\ed is injective.
\item For every $s\in S$, there exists an open neighborhood $U_s$ of $s$
and a c-lattice $B\subset \boldc(\!(z\inv)\!)^n$ such that
$L|_{U_s}\cap \theo_{U_s}\widehat{\otimes}B \subset
\theo_{U_s}(\!(z\inv)\!)^n$ is locally free of finite type over
$\theo_{U_s}$, and $L|_{U_s}\oplus \theo_{U_s}\widehat{\otimes} B
\rightarrow \theo_{U_s}(\!(z\inv)\!)^n$ is surjective.
\end{enumerate}

The functor $\GR_n: \underline{\on{Sch}}^{\on{op}}\rightarrow
\underline{\on{Sets}}$ takes, as its value on a scheme $S$, the set
\bd \GR_n (S) = \big\{ \text{d-lattices}\; L\subset
\theo_S(\!(z\inv)\!)^n\big\}. \ed

\begin{thm}[\cite{BS, AMP}]\label{Sato rep}
The functor $\GR_n$ is represented by a scheme.
\end{thm}

The scheme $\GR_n$ is called the {\em (rank $n$) Sato Grassmannian}.

It is convenient for us to use a slightly different, but equivalent,
definition of a d-lattice.

\begin{lemma}\label{whats a d-lattice}
A quasicoherent $\theo_S$-module $L$ with an injective map
$L\rightarrow \theo_S(\!(z\inv)\!)^n$ is a d-lattice if and only if
the following hold:
\begin{enumerate}
\item For every $s\in S$ there is an open neighborhood $U_s$ of $s$ and
an integer $k$ such that the map $L\rightarrow
\theo_S(\!(z\inv)\!)^n/\big(\theo_S[\![z\inv]\!]z^k\big)^n$ is
surjective with kernel a locally free $\theo_S$-module of finite
type.
\item For every $s\in S$ there is an open neighborhood $U_s$ of $s$ and
an integer $\ell$ such that the map $L\rightarrow
\theo_S(\!(z\inv)\!)^n/\big(\theo_S[\![z\inv]\!]z^\ell\big)^n$ is
injective with cokernel a locally free $\theo_S$-module of finite
type.
\end{enumerate}
\end{lemma}
\begin{proof}
First, we suppose $L$ is a d-lattice.  By the definition, for each
$s\in S$ there is an open set $U_s$ containing $s$ and a c-lattice
$B\subset \boldc(\!(z\inv)\!)^n$ giving an exact sequence \bd
0\rightarrow L|_{U_s}\cap \theo_{U_s}\widehat{\otimes}B \rightarrow
L|_{U_s} \rightarrow
\theo_{U_s}(\!(z\inv)\!)^n/\theo_{U_s}\widehat{\otimes}B \rightarrow
0 \ed over $U_s$ with kernel locally free of finite type.  Choosing
some $k$ such that $B\subseteq (\boldc[\![z\inv]\!]z^k)^n$, we get a
surjective map $f: L|_{U_s}\rightarrow
\theo_{U_s}(\!(z\inv)\!)^n/(\theo_{U_s}[\![z\inv]\!]z^k)^n$ whose
kernel sits in an exact sequence \bd 0\rightarrow L|_{U_s}\cap
\theo_{U_s}\widehat{\otimes}B \rightarrow \operatorname{ker}(f)
\rightarrow
\theo_S\widehat{\otimes}\big(\!(\boldc[\![z\inv]\!]z^k)^n/B\big)
\rightarrow 0. \ed This proves that (1) holds for $L$.

We next prove that there exists $\ell$ such that $L|_{U_s}\cap
(\theo_S[\![x\inv]\!]x^\ell)^n = 0$.  To see that such an $\ell$
exists, start with $U_s$ and $B$ such that
 $F = L|_{U_s}\cap (\theo_{U_s}\widehat{\otimes}B)$ is
finitely generated and locally free over $\theo_{U_s}$.  Choosing
$U_s$ smaller and $B$ larger if necessary, we may assume that $U_s =
\spec(R)$ is an affine scheme, that $B =
\big(\boldc[\![z\inv]\!]z^k\big)^n$, and that $F$ is a free
$\theo_{U_s}$-module of rank $N$.  So, suppose that $F =R^N$ and let
$P\subset R$ be the prime of $R$ corresponding to the point $s\in
S$.  Let $\phi: R^N\rightarrow R(\!(z\inv)\!)^n$ be the
corresponding map of $R$-modules; this map factors through
$\big(R[\![z\inv]\!]z^k\big)^n\subset R(\!(z\inv)\!)^n$ by
construction.
  By part (1) of the definition of d-lattice, the induced map
$\phi_{R_P/P}:(R_P/P)^N \rightarrow (R_P/P)(\!(z\inv)\!)^n$ is
injective.  Since $K = R_P/P$ is a field, for any $\ell$
sufficiently negative we find that the map $\phi_{K}: K^N\rightarrow
\big(K[\![z\inv]\!]z^k\big)^n/ \big(K[\![z\inv]\!]z^\ell\big)^n$ is
injective.  It follows from Nakayama's Lemma that, letting $r =
k-\ell-N$, there is an $R$-homomorphism $\psi:R^r \rightarrow
\big(R[\![z\inv]\!]z^k\big)^n/\big(R[\![z\inv]\!]z^\ell\big)^n$ such
that the localized sum \bd R_P^N\oplus R_P^r \xrightarrow{\phi
+\psi}
\big(R_P[\![z\inv]\!]z^k\big)^n/\big(R_P[\![z\inv]\!]z^\ell\big)^n
\ed is surjective, hence also an isomorphism.  Since the cokernel of
$\phi+\psi$ is finitely generated over $R$, it follows that there is
some element $f\in R\setminus P$ such that
\begin{equation}
(\phi + \psi)\otimes R_f: R_f^N \oplus R_f^r \rightarrow
\big(R_f[\![z\inv]\!]z^k\big)^n/\big(R_f[\![z\inv]\!]z^\ell\big)^n
\end{equation}
is an isomorphism.  Consequently, over $\spec(R_f)$, the map \bd
L|_{\spec(R_f)} \rightarrow
\theo_{\spec(R_f)}(\!(z\inv)\!)^n/\big(\theo_S[\![z\inv]\!]z^\ell\big)^n
\ed is injective with cokernel a locally free
$\theo_{\spec(R_f)}$-module of finite type (isomorphic to
$\theo_{\spec(R_f)}^r$), proving that (2) holds.

For the converse, part (1) of the definition of d-lattice is
automatic from statement (2) in the lemma, and part (2) of the
definition of d-lattice is immediate from statement (1) in the
lemma.
\end{proof}

\subsection{$\D$-Lattices}

We now turn to the description of the $\D$-modules that interest us.
For the remainder of the subsection, we let $\D =
\boldc[\![x]\!][\partial]$ where $\partial = \partial/\partial x$.
For a scheme $S$, we let $\D_S = \theo_S[\![x]\!][\partial]$; that
is, $\D_S$ is the sheaf on $S$ whose sections on an open set
$U\subset S$ are given by $\D_S(U) = \theo_S(U)[\![x]\!][\partial]$.
Similarly, we define a sheaf $\cE_S =
\theo_S[\![x]\!](\!(\partial\inv)\!)$ in the same way.  Note that
$\D_S$ is a quasicoherent sheaf on the formal scheme
$S\times\on{Spf}(\boldc[\![x]\!])$.  Note also that $\cE_S$ is not
quasicoherent over $S\times\on{Spf}(\boldc[\![x]\!])$. Nevertheless,
it is clear from the construction that, if $S'\xrightarrow{f}S$ is a
morphism of schemes, there is a homomorphism of sheaves of algebras
$f\inv \cE_S \rightarrow \cE_{S'}$ on $S'$.

\begin{note}
We use subscripts rather than superscripts to denote the filtration
on the sheaf $\cE$; this is nonstandard in the theory of
$\D$-modules, but makes our formulas more readable (especially in
this section).
\end{note}

For a scheme $S$, a {\em $\D$-lattice} $M\subset\cE_S^n$ is a
quasicoherent sheaf on $S\times\on{Spf}(\boldc [\![x]\!])$ with a
structure of right $\D_S$-module that is finitely presented as a
$\D_S$-module and comes equipped with an injective $\D_S$-module
homomorphism $M\rightarrow \cE_S^n$ that satisfies the following:
\begin{enumerate}
\item For all $s\in S$, there exist an open neighborhood $U_s$ of $s$ and
an integer $k$ such that $M|_{U_s}\rightarrow
\cE^n_{U_s}/\cE^n_{U_s, k}$ is surjective, with kernel a finitely
generated, locally projective $\theo_{U_s}[\![x]\!]$-module.
\item For all $s\in S$, there exist an open neighborhood $U_s$ of $s$ and
an integer $\ell$ such that $M|_{U_s}\rightarrow
\cE^n_{U_s}/\cE^n_{U_s, \ell}$ is injective, with cokernel a
finitely generated, locally projective
$\theo_{U_s}[\![x]\!]$-module.
\end{enumerate}

This definition is useful since, on ${\mathbf C}$-points, framed
$\D$-bundles give $\D$-lattices:

\begin{prop}\label{D-bundles give D-lattices}
Let $M$ be an $\theo^n$-framed $\D$-bundle.  Choose an isomorphism
$M_\cE \cong \cE^n$ compatible with the framing.  Then:
\begin{enumerate}
\item For any $k$ sufficiently
large, the map $M\rightarrow \cE^n/\cE_k^n$ is surjective, with
kernel a finitely generated projective ${\mathbf
C}[\![x]\!]$-module.
\item For any $\ell$ sufficiently small, the map
$M\rightarrow \cE^n/\cE^n_\ell$ is injective, with cokernel a
finitely generated projective ${\mathbf C}[\![x]\!]$-module.
\end{enumerate}
\end{prop}
\begin{proof}
Since we have chosen an isomorphism compatible with the framing, the
canonical filtration on $\cE^n$ induced from that on $M$
(Proposition \ref{can filt}) agrees with the standard one on
$\cE^n$. Corollary \ref{can filt props cor} implies that
$\on{gr}M\subset\on{gr}\cE^n$ is finitely generated and torsion-free
over $\on{gr}\D$ and that the inclusion is an isomorphism in all
sufficiently large graded degrees.  An inductive argument proves
that $M/M_k \cong \cE^n/\cE_k^n$ for all $k\gg 0$. Now each $M_k$ is
finitely generated over ${\mathbf C}[\![x]\!]$ by assumption and
torsion-free over ${\mathbf C}[\![x]\!]$ since $M$ is torsion-free
over $\D$.  Hence each $M_k$ is projective over ${\mathbf
C}[\![x]\!]$. This proves (1).

For (2), observe that Proposition \ref{can filt} yields that $M\cap
\cE_\ell^n=0$ for $\ell \ll 0$, so $M\rightarrow \cE^n/\cE^n_\ell$
is injective.  It follows from Corollary \ref{can filt props cor}
that $\cE^n/(\cE^n_\ell + M)$ is finitely generated over ${\mathbf
C}[\![x]\!]$, so to prove that it is projective, it suffices to
prove that it is torsion-free.

First, suppose that $\cE^n/M$ has nonzero ${\mathbf
C}[\![x]\!]$-torsion.  Then it has a $\D$-submodule $N$ that is
finitely generated over $\D$ and consists of ${\mathbf
C}[\![x]\!]$-torsion.  Taking the induced filtration on $\cE^n/M$
and thus on $N$, we get a filtration of $N$ with
$\on{gr}_\ell(N)\subset\on{gr}_\ell(\cE^n)$ for $\ell\ll 0$. But
this is impossible (since the right-hand module is ${\mathbf
C}[\![x]\!]$-torsion-free), so the filtration of $N$ must be bounded
below.  However, by Corollary \ref{can filt props cor} the
filtration is also bounded above, so $N$ is a finitely generated
torsion module over ${\mathbf C}[\![x]\!]$ equipped with a
$\D$-module structure, a contradiction.  Thus $\cE^n/M$ is ${\mathbf
C}[\![x]\!]$-torsion-free.

Now, suppose that $\cE^n/(\cE^n_\ell + M)$ has ${\mathbf
C}[\![x]\!]$-torsion submodule $T_\ell$.  For $\ell \ll 0$, we have
an exact sequence \bd 0\rightarrow \on{gr}_\ell\cE^n = {\mathbf
C}[\![x]\!]^n \rightarrow \cE^n/(M+\cE^n_{\ell-1}) \rightarrow
\cE^n/(M+\cE^n_\ell)\rightarrow 0, \ed so in particular we find that
the natural map $T_{\ell-1} \rightarrow T_\ell$ is injective.  But
$T_\ell$ is finitely generated over ${\mathbf C}[\![x]\!]$ for each
$\ell$; in particular, it is of finite length.  Now
$\underset{\longleftarrow}{\lim} \cE^n/(M+\cE^n_\ell) = \cE/M$ is
${\mathbf C}[\![x]\!]$-torsion-free, so
$\underset{\longleftarrow}{\lim} T_\ell =0$.  Since each $T_\ell$ is
of finite length, we get that $T_\ell =0$ for all sufficiently
negative $\ell$. This completes the proof.
\end{proof}

We remark that $\D$-lattices pull back well:

\begin{lemma}
If $S'\xrightarrow{f} S$ is a morphism and $M\subset \cE_S^n$ is a
$\D$-lattice, then $f^*M\rightarrow f^*\cE_S^n \rightarrow
\cE_{S'}^n$ makes $f^*M$ a $\D$-lattice.
\end{lemma}

We thus obtain a functor $\GR_{\D,n}: \underline{\on{Sch}}^{\on{op}}
\rightarrow \underline{\on{Sets}}$ by setting \bd \GR_{\D,n}(S) =
\big\{ \text{$\D$-lattices}\; M\subset \cE_S^n\big\}. \ed

The main theorem of this section is the following:

\begin{thm}\label{Sato-type thm}
The scheme $\GR_n$ represents the moduli functor $\GR_{\D,n}$ of
$\D$-lattices.
\end{thm}

\subsection{Proof of Theorem \ref{Sato-type thm}}
We begin the proof of Theorem \ref{Sato-type thm} by describing a
natural transformation between the two moduli functors.

Given a $\D$-lattice $M\subset\cE_S^n$, we get a map \bd
M\otimes\boldc[\![x]\!]/(x) \rightarrow
\cE_S^n\otimes\boldc[\![x]\!]/(x) \cong \theo_S(\!(z\inv)\!)^n, \ed
where we have identified $\partial\inv$ with $z\inv$.
\begin{prop}\label{D nat trans}
This construction defines a natural transformation of functors
$\GR_{\D,n}\rightarrow \GR_n$.
\end{prop}
\begin{proof}
Suppose that $M\subset \cE_S^n$ is a $\D$-lattice.  Over an open set
$U\subset S$, if we have $M|_U\twoheadrightarrow
\cE^n_U/\cE^n_{U,k}$ surjective with locally free kernel $K$, then
we get a surjective map $M|_U\otimes \boldc[\![x]\!]/(x)\rightarrow
(\cE^n_U/\cE^n_U x)/(\cE^n_{U,k}/\cE^n_{U,k}x) =
\theo_U(\!(z\inv)\!)^n/\big(\theo_U[\![z\inv]\!]z^k\big)^n$ with
kernel $K\otimes \boldc[\![x]\!]/(x)$, a locally free
$\theo_U$-module.

Similarly, over an open set $U\subset S$, if we take $\ell \ll 0$,
we get \bd 0\rightarrow M|_U \rightarrow
\cE_U^n/\cE^n_{U,\ell}\rightarrow W\rightarrow 0 \ed exact with $W$
finitely generated and locally free over $\theo_U[\![x]\!]$.
Tensoring with $\boldc[\![x]\!]/(x)$, we get \bd 0\rightarrow
(M|_U)/(M|_U x) \rightarrow
\theo_U(\!(z\inv)\!)^n/\big(\theo_U[\![z\inv]\!]z^k)^n \rightarrow
W\otimes \boldc[\![x]\!]/(x)\rightarrow 0 \ed exact, with $W\otimes
\boldc[\![x]\!]/(x)$ locally free over $\theo_U$. So $M/Mx
\rightarrow \cE^n/\cE^n x$ is a d-lattice.  We thus get a function
$\GR_{\D,n}(S)\rightarrow \GR_n(S)$ for every scheme $S$.
Furthermore, since tensor product over $\boldc[\![x]\!]$ commutes
with pullback, these functions commute with the pullback along
morphisms $S'\rightarrow S$, and thus define a natural
transformation, as desired.
\end{proof}

We want to prove that this natural transformation is actually an
isomorphism.  In particular, we first show that it is surjective:
given a d-lattice $L\subset \theo_S(\!(z\inv)\!)^n$, we want to
produce a $\D$-lattice $M\subset\cE_S^n$ whose fiber at $x=0$ is
$L$. Since both our functors are sheaves in the Zariski topology, it
suffices to do this over an open set $U_s$ containing $s$ for each
point $s$ of $S$.  We begin with the construction for those $L$ that
are {\em transverse} to some $\big(\theo_S[\![z\inv]\!]z^k\big)^n$.
\begin{lemma}\label{D lift}
If $L\subset\theo_S(\!(z\inv)\!)^n$ is a d-lattice satisfying
$L\oplus \big(\theo_S[\![z\inv]\!]z^k\big)^n =
\theo_S(\!(z\inv)\!)^n$, then there is an injection of
$\D_S$-modules $\D_S^n\xrightarrow{\iota}\cE_S^n$ satisfying
$\iota(\D_S^n) \oplus \cE_{S,k}^n = \cE_S^n$ and
$\iota(\D_S^n\otimes\boldc[\![x]\!]/(x)\!) = L$. Moreover, the image
$\iota(\D_S^n)$ consists of all sections $\theta$ of $\cE_S^n$ such
that $\theta\cdot D|_{x=0} \in L$ for every section $D$ of $\D_S$.
\end{lemma}
\begin{proof}
The first part follows from Theorem 7.4 of \cite{Mu} and Theorem 6.2
of \cite{LiMulasePrym}: indeed, those results construct such a $\D$-module
embedding of $\D^n$ in $\cE^n$ in the case in which the base scheme
is $\boldc$, but the proofs, which produce an element $\Psi$ of
$\GL_n(\cE_<)$ such that $\iota(m) = \Psi\cdot m$, work over any
base ring.

For the second part, let
\begin{equation}\label{multiplier module}
N = \big\{ \theta\in\cE_S^n \; \big| \; (\theta\cdot D)|_{x=0} \in L
\; \text{for all $D\in\D_S$}\big\}.
\end{equation}
Observe that $\iota(\D_S^n)\subseteq N$. Since $\iota(\D_S^n) = \Phi
\cdot \D_S^n$ for some invertible $\Phi\in\GL_n(\cE_S)$, after
multiplying by $\Phi\inv$ we may assume that $\iota(\D_S^n) =
\D_S^n$ and that $L = \theo_S[z]$.  Then any $\theta\in N$ satisfies
$\theta \cdot \theo_S[z]^n \subset \theo_S[z]^n$, which, by
Proposition 6.1 of \cite{LiMulasePrym},
implies that $\theta\in \D^n$ (again, the proof there is written for
the base $\boldc$ but the proof works equally well over any base
ring).
\end{proof}

Next, suppose that $L\subset\theo_S(\!(z\inv)\!)^n$ is any
d-lattice. Fix $s\in S$.  By restricting attention to an open
neighborhood of $s$ in $S$, we may assume, by Lemma \ref{whats a
d-lattice},
 that there are $k$ and $\ell$ such that we have
exact sequences \bd 0\rightarrow V \rightarrow L\xrightarrow{\alpha}
\theo_S(\!(z\inv)\!)^n/\big(\theo_S[\![z\inv]\!]z^k\big)^n\rightarrow
0 \ed and \bd 0\rightarrow L \rightarrow
\theo_S(\!(z\inv)\!)^n/\big(\theo_S[\![z\inv]\!]z^\ell\big)^n
\rightarrow W \rightarrow 0 \ed with $V$, $W$ vector bundles over
$\theo_S$.  Choosing a splitting $\alpha\inv$ of $\alpha$, we get an
$\theo_S$-submodule $L_- = \on{Im}(\alpha\inv)\subset L$ such that
$L_- \xrightarrow{\alpha}
\theo_S(\!(z\inv)\!)^n/\big(\theo_S[\![z\inv]\!]z^k\big)^n$ is an
isomorphism and $L/L_- \cong V$.  Similarly, shrinking $S$ further
if necessary, we may choose a splitting
$W\xrightarrow{\beta}\theo_S(\!(z\inv)\!)^n$ of
$\theo_S(\!(z\inv)\!)^n \rightarrow
\theo_S(\!(z\inv)\!)^n/\big(\theo_S[\![z\inv]\!]z^\ell\big)^n
\rightarrow W$. We then have $L\oplus\beta(W) \cong
\theo_S(\!(z\inv)\!)^n/\big(\theo_S[\![z\inv]\!]z^\ell\big)^n$. We
let $L_+ = L +\beta(W) \subset\theo_S(\!(z\inv)\!)^n$.

We now have $L_- \subset L \subset L_+$ with $L_- \oplus
\big(\theo_S[\![z\inv]\!]z^k\big)^n = \theo_S(\!(z\inv)\!)^n$, $L_+
\oplus \big(\theo_S[\![z\inv]\!]z^\ell\big)^n =
\theo_S(\!(z\inv)\!)^n$. By Lemma \ref{D lift}, we have injective
$\D_S$-homomorphisms $\D_S^n\xrightarrow{\iota_-}\cE_S^n$ and
$\D_S^n\xrightarrow{\iota_+}\cE_S^n$ such that
$\iota_-(\D_S^n)|_{x=0} = L_-$ and $\iota_+(\D_S^n)|_{x=0} = L_+$.
Lemma \ref{D lift} also implies that
$\iota_-(\D_S^n)\subset\iota_+(\D_S^n)$. By construction, moreover,
$\iota_+(\D_S^n)/\iota_-(\D_S^n) \cong \cE_{S,k}^n/\cE_{S,\ell}^n$
as $\theo_S[\![x]\!]$-modules.

\begin{prop}\label{lattice lift}
There is a finitely presented $\D_S$-module $M\subset\cE_S^n$
satisfying:
\begin{enumerate}
\item $\iota_-(\D_S^n)\subset M \subset \iota_+(\D_S^n)$.
\item $M\otimes \boldc[\![x]\!]/(x) = L$.
\item $M$ consists of all $\theta\in\cE_S^n$ such that
$\theta\cdot D|_{x=0} \in L$ for all sections $D$ of $\D_S$.
\end{enumerate}
\end{prop}
\begin{proof}
Consider the quotient $\D_S$-module $Q =
\iota_-(\D_S^n)/\iota_+(\D_S^n)$; by construction, this is a
finitely generated locally free $\theo_S[\![x]\!]$-module, hence,
working locally over $S$, choosing a basis, and changing to a left
$\D_S$-module, $Q$ corresponds to a vector bundle with flat
connection.  The usual construction of flat sections with power
series coefficients shows that there is a bijection between vector
subbundles of $Q/Qx$ and right $\D_S$-submodules of $Q$.  In
particular, choosing the vector subbundle $L/L_-$ of $L_+/L_-$,
there is a lift to a $\D_S$-module $M$, $\iota_-(\D_S^n)\subset M
\subset \iota_+(\D_S^n)$, with $M/Mx = L$.

For part (3), define a module $N$ as in \eqref{multiplier module},
and note that Lemma \ref{D lift} implies that we have
$\iota_-(\D_S^n)\subseteq M \subseteq N \subseteq \iota_+(\D_S^n)$.
Then $N/M$ defines a $\D_S$-submodule of $Q' = \iota_+(\D_S^n)/M$
all of whose sections $\theta$ satisfy $(\theta \cdot D)|_{x=0} =
0\in Q'/Q'x$; the standard results for flat connections then imply
that $N/M =0$, as desired.
\end{proof}

It follows that the natural transformation
$\GR_{\D,n}\rightarrow\GR_n$ is Zariski-locally surjective, hence
surjective.

Finally, suppose that $M'\subset\cE_S^n$ is a $\D$-lattice with
$M\otimes\boldc[\![x]\!]/(x) = L$.  Let $M\subset\cE_S^n$ be the
$\D$-lattice lifting $L$ that was constructed by Proposition
\ref{lattice lift}.  By part (3) of Proposition \ref{lattice lift},
$M'\subseteq M$. Since $M$ and $M'$ are $\D$-lattices, the quotient
$M/M'$ is a finitely generated
$\theo_S[\![x]\!]$-module: this follows by applying part (2) of the
Definition of $\D$-lattice and considering $M/M'\subset \cE^n/\cE^n_\ell$
for $\ell$ sufficiently small, to find that $M/M'$ is a submodule of a
finitely generated $\theo_S[\![x]\!]$-module, hence is itself finitely
generated (the ring $\theo_S[\![x]\!]$ is Noetherian).

  Tensoring the exact sequence \bd
M'\rightarrow M\rightarrow M/M'\rightarrow 0 \ed with
$\boldc[\![x]\!]/(x)$ and noting that the map from $M'$ to $M$
becomes the identity map of $L$, we find that
$(M/M')\otimes\boldc[\![x]\!]/(x) = 0$. It follows from Nakayama's
Lemma (Theorem 2.2 of Matsumura) that there is an element $a$ of
$\theo_S[\![x]\!]$ such that $a(M/M') = 0$ and $a\cong 1 \mod (x)$.
But the completeness of $\theo_S[\![x]\!]$ with respect to $(x)$
then implies that $a$ is a unit, and thus that $M/M'=0$. So the
natural transformation $\GR_{\D,n}\rightarrow\GR_n$ is injective.
This completes the proof of Theorem \ref{Sato-type thm}.\hfill\qedsymbol

\section{Appendix: $D$-Algebras and Log Differential Operators}\label{D
algebras}

In this appendix, we explain the construction of the sheaves $\Dlog$ on
singular cubics that we use in the body of the paper.  We also
explain
some facts about filtrations on $\D$-modules and microlocalization that
are needed for our main theorems.

\subsection{$D$-Algebras}
Let $X$ denote a reduced and irreducible quasiprojective complex
variety. A {\em differential $\theo_X$-bimodule} is a quasicoherent
sheaf $M$ of $\theo_X$-bimodules, such that, if we give $M$ the
filtration defined by $M_{-1}= 0$ and $M_i = \{m\in M \;|\; [r,m]\in
M_{i-1}\;\text{for all}\; r\in \theo_X\}$, then $M = \cup_i M_i$.

\begin{defn}[1.1.4 of \cite{BB}]
A {\em $D$-algebra} on $X$ is a quasicoherent sheaf $\cD$ of
associative algebras on $X$ equipped with an algebra homomorphism
$\theo_X\rightarrow \cD$ making $\cD$ a differential
$\theo_X$-bimodule.  As above, a $D$-algebra $\cD$ comes with a
canonical filtration; to be consistent with the usual notation for
the sheaf of rings of differential operators, we write $\cD^i$ for
the $i$th term in the filtration.

A $D$-algebra $\cA$ on $X$ is called a {\em special $D$-algebra}
(Section 6.2 of \cite{PRo}) if the terms $\on{gr}_i\cA$ in the
associated graded for the standard bimodule filtration are of the
form $\on{gr}_i(\cA) = \Delta_* \cF_i$ for some {\em free}
$\theo_X$-modules $\cF_i$ (i.e. the terms in the associated graded
are trivial vector bundles).
\end{defn}

The standard example of a $D$-algebra on a scheme $X$ is the sheaf
of rings of differential operators \cite{BB}.

\begin{lemma}\label{facts about D-algebras}
Let $X$ be a reduced, irreducible complex variety, and $U\subseteq
X$ a nonempty open subvariety.  Let $\D_X$ denote the sheaf of rings
of differential operators on $X$.  Then:
\begin{enumerate}
\item $\D_X(X)\subseteq \D_X(U)$.
\item $\D_X^k(X) = \D_X(X)\cap \D_X^k(U)$.
\item  Suppose
$\cA\subseteq\D_X$ is a subalgebra containing $\theo_X$.  Then $\cA$
is a $D$-algebra with canonical filtration given by $\cA_k = \cA\cap
\D_X^k$.
\end{enumerate}
\end{lemma}
\noindent
The proofs follow standard arguments.

Recall (Section \ref{Dlog introduced}) the definition of the sheaf
$\Dlog$ of log differential operators on a cubic curve, generated by
functions and the $\G$-invariant vector fields.

\begin{corollary}\label{special D}
Let $E$ be a Weierstrass cubic.  Then:
\begin{enumerate}
\item $\Dlog^k = \Dlog\cap \D^k_{\G}$.
\item $\Dlog$ is a special $D$-algebra with associated graded isomorphic
to $\Sym\theo_E$.
\end{enumerate}
\end{corollary}
\begin{proof}
(1) $\Dlog$ is a $D$-algebra by part (3) of Lemma \ref{facts about
D-algebras}.  We then have
\begin{align*}
\Dlog^k & = \Dlog\cap D_X^k \;\text{by part (3) of Lemma \ref{facts
about
D-algebras}}\\
& = \Dlog\cap D^k_{\G} \;\text{by part (2) of Lemma \ref{facts about
D-algebras}}.
\end{align*}
(2) By part (1) above, we have \bd \on{gr}(\Dlog) \subseteq
\on{gr}(\D_{\G}) = \Sym\theo_{\G}. \ed Since $\Dlog$ is generated
over $\theo_E$ by the vector field $X$, $\on{gr}(\Dlog)$ is
generated over $\theo_E$ by the image of the vector field $X$ in
$\D^1_{\G}/\D^0_{\G}$.  The $\theo_E$-submodule generated by this
element is exactly $\Sym\theo_E\subset \Sym\theo_{\G}$, as
desired.\end{proof}

\subsection{Microlocalization and Torsion-Free
$\D$-Modules}\label{microloc section}

For the remainder of this section, we fix a reduced and irreducible
quasiprojective complex curve $X$ and a
 $D$-algebra $\cD$ on $X$.

We will call the $D$-algebra $\D$ {\em symmetric} if it satisfies
$\on{gr}(\D) = \Sym(B)$ as algebras for some line bundle $B$ on $X$.

We will write $\cE$ for the  microlocalization of $\D$; for
background on the procedure of microlocalization and its properties,
see for example \cite{Kashi, vdE, AVV,Shap}. Note that, by our
assumption that $\D$ is symmetric, the associated graded ring
$\on{gr}(\cE)$ for the natural filtration is isomorphic to
$\oplus_{k\in{\mathbf Z}} B^k$ as a graded $\theo_X$-algebra.

By a {\em finitely generated (left or right) $\cD$-module} on $X$ we
mean a $\cD$-module that is quasicoherent as a sheaf of
$\OX$-modules and locally finitely generated over $\D$.

A symmetric $D$-algebra is, in particular, a sheaf of noetherian
domains by Proposition 2.6.1 of \cite{Bjork}. It follows by Theorem
1.15 of \cite{MR} that, for any affine open set $U\subseteq X$, the
ring $\D(U)$ has a quotient skew field, which we denote by $Q(\D)$
(it does not depend on $U$ since $X$ is irreducible). Consequently,
we may define a {\em torsion-free} $\D$-module to be a (left or
right) $\D$-module $M$ such that for every open set $U$ of $X$ and
all nonzero $s\in H^0(U,\D)$ and nonzero $m\in H^0(U,M)$ we have
$s\cdot m\neq 0$. If $M$ is a finitely generated and torsion-free
right $\D$-module and $U$ is an affine open set of $X$,
then $H^0(U,M)\subset H^0(U,M)\otimes Q(\cD) \cong Q(\cD)^r$ for
some integer $r$ which we refer to as the {\em rank} of $M$ (and
similarly for left $\D$-modules); note that this does not depend on
the choice of $U$.  Moreover, in this case there is
some vector bundle $W$ on $X$ such that $M$ embeds in the
$\cD$-module $\cD\otimes_{\OX} W$.
We then have the following standard fact about the microlocalization of
a finitely generated torsion-free $\D$-module.
\begin{lemma}
If $M$ is a torsion-free and finitely generated left (or right)
$\D$-module then the natural map $M\rightarrow \cE\otimes M$
(respectively $M\rightarrow M_\cE = M\otimes \cE$) is injective.
\end{lemma}

Recall that a {\em good filtration} of a $\D$-module $M$ is a
filtration $M = \cup F_i(M)$ that makes $M$ a filtered $\D$-module
and such that $\on{gr}_F(M)$ is a finitely generated
$\on{gr}(\cD)$-module.  We often write $M_i$ in place of $F_i(M)$
when there can be no confusion.

Suppose that $M$ is a $\D$-module with good filtration
$F_{\bullet}$. Then $M_{\cE} = M\otimes_{\cD}\cE$ has a canonical
sequence of filtrations induced from $F_{\bullet}$: set
$F_k(M_\cE,\ell) = F_\ell(M)\cdot\cE_{k-\ell}$; this is the $k$th
term in the $\ell$th filtration.  Note that for some choices of
$\ell$ it may not be the case that $\bigcup_{k\geq 0}F_k(M_\cE,\ell)
= M_\cE$.

We have the following easy lemma.

\begin{lemma}\label{incrfilt}
Suppose that $M$ is a $\D$-module with good filtration
$F_{\bullet}$.
\begin{enumerate}
\item $F_{\bullet}(M_\cE,\ell)$ makes $M_\cE$ into a filtered $\cE$-module
for each $\ell$.
\item
For each $\ell$ and $k$ we have $F_k(M_\cE,\ell-1) \subseteq
F_k(M_\cE,\ell)$.
\item
There exists $\ell\gg 0$ such that for all $n\geq 0$ and all $k$,
$F_k(M_\cE,\ell+n) = F_k(M_\cE,\ell)$.  Moreover, for any such
$\ell$, $\bigcup_k F_k(M_\cE,\ell) = M_\cE$.

\end{enumerate}
\end{lemma}

\begin{defn}\label{canonical filt}
We will refer to the filtration $F_{\bullet}(M_\cE,\ell)$ for
sufficiently large $\ell$ as in part (3) of Lemma \ref{incrfilt} as
the {\em induced filtration} or {\em canonical filtration}
 of $M_\cE$, and will simply write
$F_{\bullet}(M_\cE)$ or even $M_{\cE,\bullet}$ for its terms.  We also
call the induced filtration $M_k = M\cap M_{\cE,k}$ on $M$ the
{\em canonical filtration} on $M$.
\end{defn}
We
then have the following:

\begin{corollary}\label{can filt props cor}
Suppose that $M$ is a $\D$-module equipped with a good filtration
$F_{\bullet}$ such that $\bigoplus_{k\geq n}\on{gr}^F_k(M)$ is a
torsion-free $\on{gr}(\cD)$-module for some $n$.  Then:
\begin{enumerate}
\item $\operatorname{gr}_k(M_\cE) = \operatorname{gr}^F_k(M)$
for all $k$ sufficiently large.
\item   $F_{\bullet}(M_\cE)$ is a
filtration of $M_\cE$ with torsion-free associated graded.
\end{enumerate}
\end{corollary}
\begin{proof}
For $k\geq \ell\gg 0$, we have $F_k(M_\cE,\ell) = F_k(M_\cE,k) =
F_k(M)\cdot \cE_0$.  Also, for $k\geq\ell$, we have
 $F_{k-1}(M_\cE,\ell) = F_{k-1}(M_\cE,k) =  F_k(M)\cdot \cE_{-1}$.
Since $F_k(M)\cdot\cE_0 = F_k(M) + F_k(M)\cdot\cE_{-1}$, the map
$F_k(M)/F_{k-1}(M)\rightarrow F_k(M_\cE)/F_{k-1}(M_\cE)$ becomes \bd
\frac{F_k(M)}{F_{k-1}(M)}\longrightarrow
\frac{F_k(M)}{F_k(M)\cdot\cE_{-1}\cap M} \longrightarrow
\frac{F_k(M)+F_k(M)\cdot\cE_{-1}}{F_k(M)\cdot\cE_{-1}}, \ed which is
evidently surjective.  To prove that it is injective, it suffices,
by our assumption,
 to prove that its restriction to some nonempty
open set $U\subset X$ is injective.

By choosing $U$ sufficiently small and $n$ sufficiently large, we
may assume that $\on{gr}_n(M) = W$ is a vector bundle, that the
symbol map $F_n(M)\rightarrow \on{gr}_n(M)=W$ is a split surjection,
and, choosing a splitting $W\rightarrow F_n(M)$ as
$\theo_X$-modules, that the natural map $W\otimes\D\rightarrow
F_n(M)$ is an isomorphism on associated graded modules in degrees
greater than or equal to $n$. It follows by a standard inductive
argument that the
induced map of filtered $\cE$-modules $W\otimes \cE\rightarrow
M_\cE$ is an isomorphism.  This proves (1). Part (2) is then
immediate.
\end{proof}

Let $M$ be a finitely generated, torsion-free right $\D$-module equipped
with a good filtration such that $\on{gr}(M)$ is a torsion-free
$\on{gr}(\D)$-module in high degree---we will call $M$ a {\em $\D$-module
with torsion-free filtration} for simplicity.

\begin{prop}\label{can filt}
Suppose $M$ is a $\D$-module with torsion-free filtration on $X$.
Then the canonical
filtration on $M_\cE$ induced by the given filtration on $M$ is the
unique $\cE$-module filtration on $M_\cE$ that satisfies:
\begin{enumerate}
\item[(a)] The induced map $M\rightarrow M_\cE$ is a homomorphism of filtered
$\D$-modules.
\item[(b)] The induced map $\on{gr}_k(M)\rightarrow \on{gr}_k(M_\cE)$ is an
isomorphism of $\theo_X$-modules for $k\gg 0$.
\end{enumerate}
Furthermore, this filtration has the properties:
\begin{enumerate}
\item[(i)] $M_k= M\cap M_{\cE,k}$ for all $k\gg 0$.
\item[(ii)] $M_{\cE,k} = \cup_{n\geq 0} M_{k+n}\cdot \cE_{-n}$.
\item[(iii)] $\on{gr}(M_\cE) \cong V\otimes \on{gr}(\cE)$ as
graded $\on{gr}(\cE)$-modules.
\end{enumerate}
\end{prop}
\begin{proof}
The canonical filtration satisfies (a) and (b) by construction and
Corollary \ref{can filt props cor}.  If $\Theta_{\bullet}$ is any
other filtration of $M_\cE$ satisfying (a) and (b), then
$\Theta_k(M_\cE) = \Theta_k(M_\cE)\cE_0 \subseteq F_k(M)\cE_0
=F_k(M_\cE)$ for $k$ sufficiently large.  Since $\Theta_{k-n}(M_\cE)
= \Theta_k(M_\cE)\cE_{n-}$ and $F_{k-n}(M_\cE) = F_k(M_\cE)\cE_{-n}$
for all $n\geq 0$, we find that $F_k(M_\cE)\subseteq
\Theta_k(M_\cE)$ for all $k$.

The induced map $\on{gr}^F_k(M_\cE)\rightarrow
\on{gr}^{\Theta}_k(M_\cE)$ is an isomorphism for $k$ large by
Corollary \ref{can filt props cor} and assumption (b), hence for all
$k$ since it is the restriction of a homomorphism of torsion-free graded
$\on{gr}(\cE)$-modules.  It is immediate that $F_{k-1}(M_\cE) =
F_k(M_\cE) \cap \Theta_{k-1}(M_\cE)$ for all $k$. Suppose, by way of
inductive hypothesis, that $F_{k+n}(M_\cE)\cap\Theta_{k-1}(M_\cE) =
F_{k-1}(M_\cE)$.  Then \bd F_{k+n+1}\cap\Theta_{k-1} =
(F_{k+n+1}\cap \Theta_{k+n}) \cap \Theta_{k-1} = F_{k+n}\cap
\Theta_{k-1} = F_{k-1}. \ed It follows that $F_{k+n}\cap \Theta_k
=F_k$ for all $k$. Since $M_\cE = \cup_{k+n}F_{k+n}(M_\cE)$, we get
$\Theta_k(M_\cE) = F_k(M_\cE)$, proving the uniqueness statement.

Properties (ii) and (iii) are clear from the constructions. To prove
(i), let $F'_k(M) = M\cap F_k(M_\cE)$.  Then $F_k\subseteq F'_k$ for
all $k$.  The map $\on{gr}^{F'}_k(M)\rightarrow \on{gr}_k(M_\cE)$ is
injective by construction, hence by (b) is an isomorphism for $k$
sufficiently large.  It follows that the map
$\on{gr}_k^F(M)\rightarrow \on{gr}_k^{F'}(M)$ is an isomorphism for
$k$ sufficiently large, say $k\geq k_0$.  It follows that $F_k\cap
F'_{k-1} = F_{k-1}$ for $k\geq k_0$.  An inductive argument just
like that for $F$ and $\Theta$ above then proves that $F_k = F'_k$
for all $k\geq k_0$, as desired.
\end{proof}

Any finitely generated torsion-free $\D$-module can be equipped with
a torsion-free filtration locally on $X$, and consequently we have:
\begin{lemma}\label{no O coherent}
If $M$ is a finitely generated torsion-free $\D$-module,
then $M$ has no nonzero
$\theo$-coherent $\D$-submodules.
\end{lemma}

In fact, more is true: every finitely generated torsion-free $\D$-module
can be equipped with a torsion-free filtration globally on $X$, and in
rank $1$ this filtration is essentially unique.

\begin{prop}\label{uniqueness of filtration}
Suppose that $M$ is a finitely generated, torsion-free right $\D$-module
on $X$.
\begin{enumerate}
\item Let $n=\on{rk}(M)$.  Then $M$ admits a filtration such that $\on{gr}(M)$
is a torsion-free $\on{gr}(\D)$-module and $\on{rk}(\on{gr}_k(M))=n$ for
$k\geq 0$ and $0$ for $k<0$.
\item Suppose that, in addition, $M$ has rank $n=1$ and that $F^1$ and
$F^2$ are two $\cD$-module filtrations of $M$ such that for some
$n_0$ sufficiently large, $\bigoplus_{k\geq n_0} \gr_k(M,F^i)$ is a
torsion-free $\gr(\cD)$-module for $i=1,2$.  Then there exists
$l\in{\mathbf Z}$ such that for all $k\geq \min\{n_0, n_0+l\}$,
$F^1_{k+l} = F^2_k$.
\end{enumerate}
\end{prop}
\begin{proof}
(1) Since $M$ is torsion-free, there is a nonempty open subset $U\subseteq X$
such that $M|_U\cong \D^n_U$.  Let $D=X\setminus U$.  Since $M$ is
finitely generated, there exists $\ell$ such that the inclusion
$M\hookrightarrow \D^n_U$ factors through $\theo(\ell D)\otimes \D^n$; the
map $M\hookrightarrow \theo(\ell D)\otimes \D^n$ is also an isomorphism over
$U$.  The induced filtration on $M$ then has the desired property.

(2) The question is local on $X$, so we may restrict to an affine open subset
$U$ of $X$.  The claim then follows from Lemma 3.2 of \cite{NS}.
\end{proof}

\bibliographystyle{alpha}
\newcommand{\bl}{/afs/math.lsa.umich.edu/group/fac/nevins/private/math/bibtex/}

\bibliography{\bl baranovsky,\bl bradlow,\bl brosius,\bl
bialynicki-birula, \bl carrell,\bl donagi,\bl friedman,\bl
frankel,\bl ginzburg,\bl gomez,\bl griffiths, \bl grojnowski,\bl
hartshorne,\bl hitchin,\bl hungerford,\bl huybrechts,\bl jardim,\bl
kapranov,\bl kapustin,\bl kirwan, \bl kurke,\bl langton,\bl
laumon,\bl lehn,\bl matsumura,\bl milne, \bl nakajima,\bl nevins,\bl
segal,\bl sernesi,\bl thaddeus,\bl viehweg,\bl wilson}

\end{document}